\documentclass[opre]{arxiv_version}
\OneAndAHalfSpacedXI

\newcommand{\lr}{L^{\eta}}
\newcommand{\fr}{{\gamma}_{*}^{\eta}}

\newcommand\inner[2]{\left \langle #1, #2 \right \rangle}
\newcommand{\etal}{\eta\lambda^{*}}
\newcommand{\etam}{\eta\mu^{*}}
\newcommand{\con}{\gamma}
\newcommand{\olam}{\BFlambda^{*}}
\newcommand{\omu}{\BFmu^{*}}

\newcommand{\bth}{\BFtheta^{\eta}}

\newcommand{\bph}{\BFphi^{\eta}}
\usepackage{algorithm}
\usepackage{algorithmic}
\usepackage{soul}
\usepackage{color}

\usepackage{bbm}
\usepackage{enumitem}
\usepackage{graphics}
\usepackage{tikz}
\usetikzlibrary{shapes}
\usepackage{multirow}
{\theoremstyle{TH}\newtheorem{condition}{Condition}}
\usetikzlibrary{arrows}
\usetikzlibrary{automata,arrows,positioning,calc}
\usepackage{hyperref}
\usepackage{xcolor}
\newcommand{\E}{\mathbb{E}}
\renewcommand{\Pr}{\mathbb{P}}
\usepackage{natbib}
 \bibpunct[, ]{(}{)}{,}{a}{}{,}%
 \def\bibfont{\small}%
\TheoremsNumberedThrough     
\ECRepeatTheorems
\EquationsNumberedThrough    
\MANUSCRIPTNO{}
\begin{document}
\RUNAUTHOR{Varma et al.}
\RUNTITLE{Dynamic Pricing and Matching in Two-Sided Queues}
\TITLE{\Large Dynamic Pricing and Matching for Two-Sided Queues}
\ARTICLEAUTHORS{%
\AUTHOR{Sushil Mahavir Varma,\;\ Pornpawee Bumpensanti,\;\ Siva Theja Maguluri,\;\ He Wang}
\AFF{School of Industrial and Systems Engineering, Georgia Institute of Technology, Atlanta, GA 30332\\ 
\EMAIL{sushil@gatech.edu, pornpawee@gatech.edu, siva.theja@gatech.edu, he.wang@isye.gatech.edu}}
} 
\ABSTRACT{%
Motivated by applications from gig economy and online marketplaces, we study a two-sided queueing system under joint pricing and matching controls.  The queueing system is modeled by a bipartite graph, where the vertices represent customer or server types and the edges represent compatible customer-server pairs.
Both customers and servers sequentially arrive to the system and join separate queues according to their types. The arrival rates of different types depend on the prices set by the system operator and the expected waiting time. At any point in time, the system operator can choose certain customers to match with compatible servers. The objective is to maximize the long-run average profit for the system.
We first propose a fluid approximation based pricing and max-weight matching policy, which achieves an $O(\sqrt{\eta})$ optimality rate when all the arrival rates are scaled by $\eta$. We further show that a two-price and max-weight matching policy achieves an improved $O(\eta^{1/3})$ optimality rate. 
Under a broad class of pricing policies, we prove that any matching policy has an optimality rate that is lower bounded by $\Omega(\eta^{1/3})$. Thus, the latter policy achieves the optimal rate with respect to $\eta$. 
We also demonstrate the advantage of max-weight matching with respect to the number of server and customer types $n$. Under a complete resource pooling condition, we show that max-weight matching achieves $O(\sqrt{n})$ and $O(n^{1/3})$ optimality rates for static and two-price policies, respectively, and the latter matches the lower bound $\Omega(n^{1/3})$.  In comparison, the randomized matching policy may have an $\Omega(n)$ optimality rate.
}%
\KEYWORDS{queueing, dynamic pricing, max-weight matching, Markov decision process} 
\maketitle
%
\section{Introduction}
Most queueing models consider a fixed set of servers with sequentially arriving customers.
In this paper, however, we consider a two-sided queueing system where servers also arrive sequentially and then wait to be matched with customers. 
Various applications of online marketplaces and gig economy platforms can be modeled as two-sided queues---for example,
Uber and Lyft where passengers are matched with drivers, Uber Eats and DoorDash where customer orders are matched with meal delivery couriers, and crowdsourced workforce platforms such as TaskRabbit where tasks are matched with contributors.
Most of these platforms use both dynamic pricing and dynamic matching as levers to control the marketplaces. 
Motivated by these applications, we consider a canonical model of two-sided queues with multiple types of servers and customers. Each customer type is compatible with a subset of server types. For example, in the case of ride-hailing marketplaces, the types of servers (drivers) and customers are determined by the proximity of their current locations, as well as other factors such as the numbers of seats requested by passengers and the vehicle capacities.
Our model assumes a fairly general setting with arbitrary numbers of customer and server types, with their compatibility modeled by a bipartite graph.
At each point in time, the system operator posts a price for each customer and server type. Then, customers and servers who are willing to accept the quoted prices (after they factor in expected waiting costs) will enter the system. Those who entered will wait in queues separated by their types until they are matched to a compatible counterpart type. Once a customer-server pair is matched, the pair will leave the queueing system immediately to complete the service. The system operator earns a profit that is equal to the difference between the price charged to the customer and the price quoted to the server.
We formulate the above system as a Markov decision process (MDP) in the infinite time horizon. 
The operator can vary the prices for different customer and server types, as well as decide when to match and which customer-server pair to match.
The objective is to maximize the long-run average profit obtained by the system operator. 
There are several technical challenges to analyze such a stochastic system.
The first challenge is the \emph{curse of dimensionality} in solving and analyzing the MDP. As the number of customer or server types increases, the dimension of the state space increases exponentially (even when the buffer size of each queue is bounded). 
It is hence intractable to solve the exact MDP for large-scale systems with multiple types.
In this paper, we propose several {approximate} policies to obtain near optimal solutions for the MDP. 
The second challenge is that the stochastic behavior of the two-sided queueing system is complicated by
the interplay between pricing and matching decisions.
Our proposed policies use dynamic pricing to ensure the stability of the two-sided queue system, so the arrival rates of customers and servers vary with the queue lengths. As the queue lengths change, the matching decisions among different types are adjusted dynamically, which in turn affects the system state and pricing decisions.
As a result, the queue lengths of different types are intricately correlated. The system cannot be decomposed into a set of simple queue and the pricing and matching decisions cannot be decoupled and analyzed separately.
To solve this challenge, we use the Lyapunov drift method to 
analyze the stochastic system as a whole in order to bound the total queue length.
\subsection{Summary of Results}
We first present a fluid model for the two-sided queueing system and show that the profit obtained by the fluid model is an upper bound on the achievable profit under any policy. Based on the fluid model, we propose several pricing and matching policies.
In Section~\ref{sec:large-theta}, we analyze the proposed policies in a \emph{large-scale} regime in which the arrival rates of all types are scaled by a factor $\eta \to \infty$.
We consider a static pricing policy using the fluid solution combined with the max-weight matching algorithm.
We show that the profit loss of this policy from the fluid solution benchmark is $O(\sqrt{\eta})$ (Section~\ref{subsec: fluidsolution}). 
We then propose a generalization of the fluid pricing policy that uses two prices for each queue type \citep[see][]{kim2017value}.
For the two-price policy combined with max-weight matching policy, we show that the profit loss from the fluid solution benchmark is reduced to $O(\eta^{1/3})$ (Section~\ref{sec:two-price}). Furthermore, we prove that for a broad class of pricing policies, using any matching policy will result in a profit loss lower bounded by $\Omega(\eta^{1/3})$ (Section~\ref{sec:lower_bound}). 
In Section~\ref{sec:large-n}, we consider a \emph{large-system} regime in which both the number of server/customer types and the arrival rates are scaled ($n\to \infty, \eta \to \infty$).
We show that the max-weight algorithm is \emph{delay optimal}. In particular, 
max-weight matching minimizes the revenue loss under fluid pricing and two-price policies among all matching policies (Section~\ref{sec:delay-optimal}).
Under the complete resource pooling condition, we characterize the profit loss of max-weight matching:
the profit loss scales as $O(\sqrt{n\eta})$ for fluid pricing policy and $O((n\eta)^{1/3})$ for two-price policy (Section~\ref{sec:MW vs Rand}). 
Furthermore, we establish a lower bound showing that any pricing and matching will incur a profit loss of $\Omega((n\eta)^{1/3})$, so the two-price max-weight policy is asymptotically optimal (Section~\ref{sec:lower-bound-n}).
In contrast, if one directly applies the solution of the fluid model as a state-independent randomized matching policy, the profit loss scales as $O(n\eta^{1/2})$ for the fluid pricing policy and $O(n\eta^{1/3})$ for the two-price policy.
In Appendix~\ref{sec:mdp-appendix}, we further analyze the structure of the MDP model and propose approximate DP solutions. 
In some special cases, we are able to show structural properties of the optimal dynamic pricing policy. In addition, we present an LP-based approximation technique with a constraint generation algorithm to solve the MDP efficiently. 
\subsection{Literature Review}
\label{sec:liter_review}
\subsubsection*{Dynamic Matching.} Dynamic matching markets have numerous applications such as ride sharing \citep{banerjeeridehailing}, e-commerce marketplaces like Amazon.com or Ebay, kidney exchange \citep{roth2007kidney,yashkanoriabarter}, and payment processing networks \citep{spider}.
Below, we will discuss previous work involving dynamic matching in the context of two-sided queues.
 \cite{caldentey2009fcfs} and  \cite{adan2012exact} considered bipartite matching for two-sided queues on a first-come-first-served basis: each arriving customer is matched to a compatible server who has the earliest arrival time and has not been matched.
Under this matching rule, they analyzed steady-state matching rates between certain customer and server types.
 Furthermore, they deduced the necessary conditions on the frequency of arrivals for stability of the system and also derived the stationary distribution.
 \cite{matchingqueues} analyzed a general multi-sided queuing system, which includes the two-sided queueing system as a special case. Their objective was to minimize holding cost in a finite horizon. They presented a periodic review matching algorithm and showed asymptotic optimality as arrival rates become large. 
 \cite{dynamictypematchinghu} studied a two-sided matching system similar to ours. Their goal is to maximize the discounted reward obtained by matching customers and servers in a finite horizon, while accounting for the holding costs. They study conditions such that a priority rule is optimal. In addition, they present a matching algorithm based on fluid approximation and show that it is asymptotically optimal. 
The main distinction of \cite{dynamictypematchinghu} with our paper is that they do not consider dynamic pricing. In addition, while they use fluid approximation to generate static (open-loop) matching decisions, we use max-weight algorithm to generate (closed-loop) matching decisions that are adaptive to queue lengths.
\citet{chen2019pricing} studied a dynamic pricing and matching problem in which strategic customers and servers arrive dynamically and have heterogeneous waiting costs. Their paper assumes all customers and servers are compatible and considers a greedy matching policy on a first-come-first-served basis.
Dynamic matching problems were also studied in the context of kidney exchanges albeit in a non-two-sided setting in \cite{yashkanoriabarter,akbarpour2017thickness}. Pricing is usually forbidden in kidney exchanges due to ethical and legal reasons.
These papers study the value of ``batching'', i.e., holding compatible matching pairs in hope that better matching will arrive in future. However, both papers find that batching in general does not provide significant benefit.
\subsubsection*{Dynamic Pricing for Queues.} 
First we discuss the literature involving dynamic pricing in the context of single-sided queues and then review those involving two-sided queues. 
\cite{pricinglow1974} is one of the earlier works studying dynamic pricing in a single sided queue. The paper considered price dependent customer arrivals with a finite buffer; the rewards include the payment by customers and holding costs incurred by the operator. Monotonicity of the optimal pricing policy is showed. It was later extended to infinite buffer capacity in \cite{low1974optimal}.
 \cite{pricinginqueuing2001} considered a queuing model with customers who are sensitive to both waiting time and price. 
They presented structural properties on optimal pricing decisions and monotonicity of optimal bias function.
In the context of network services like call centers,  \cite{Tsitsiklis2000congestion} considered a system with finite total resource. They consider different types of price dependent customers arrivals which requests for a fraction of the resource. The objective is to find a pricing policy to maximize revenue. They show multiple structural properties like concavity of value function and monotonicity of optimal policy.
 \cite{kim2017value} considers a single server queuing system and studies the benefit of dynamic pricing over static pricing. They assume that the customers are delay sensitive and consider a revenue maximization objective. They present a static pricing policy and a two-price policy, and also provide the rate of convergence of these policies. 
 Our two-price policy considered in Section~\ref{sec:two-price} is motivated by the results from \cite{kim2017value}.
 The method of \cite{kim2017value} involves applying the Taylor series expansion to the revenue function and then bounding the expected steady-state queue length.
 The main distinction of \cite{kim2017value} with our paper is that they consider a single server queue, whereas we consider a network of two-sided queues with matching decisions. It is non-trivial to generalize the method presented in \cite{kim2017value} to a two-sided queueing network,
 as an exact analysis of the steady-state distribution is intractable due to the complex interaction among different queues.
 In addition, unlike the single server setting in \cite{kim2017value}, matching decisions play a critical role in our model and cannot be decoupled from the pricing decisions.
Aside from establishing asymptotic rates with large arrival rates, we also complement the result in \cite{kim2017value} by showing the advantage of two-price pricing (when combined with appropriate matching policies) for large network sizes.

The joint problem of dynamic pricing and matching was also studied by \cite{amywardpricingmatching} 
under the objective of maximizing the number of successful matches. They proposed an asymptotically optimal pricing and matching policy with large arrival rates.
 The differences with our work are that they proposed static policies based on the fluid model and analyzed the system for a finite time horizon.
A two-sided queueing model with both customer and server arrivals is studied by  \cite{ondemandservers}. They consider a setting where the arrival rate of the servers can be controlled.
However, the focus in \cite{ondemandservers} was to establish system stability and process level convergence, while the objective in our model is to maximize profit.
Several recent papers have studied dynamic pricing in the context of ride hailing systems \citep{korolko2018dynamic, besbes2018surge, hu2019surge}.  \cite{banerjeeridehailing,banerjee2018state} studied a closed queuing network, where the number of cars in the system is a constant and the customers abandon the system if they are not matched immediately. \cite{banerjeeridehailing} considered a state-independent pricing policy and prove the approximation ratio with respect to optimal pricing policy. \citep{banerjee2016dynamic} proposed a state-dependent pricing policy and argue that the benefit of dynamic pricing is in the robustness of the performance of the system.
In sum, most of the previous work on dynamic matching is either in the context of single-sided queues or not coupled with revenue optimization. Of the few that consider both of these, the matching policy considered is an open-loop policy. On the other hand, we consider all of these aspects and show the asymptotic optimality under closed-loop matching policies. 
\subsubsection*{Max-Weight Algorithm.}
 In this work, we apply a max-weight matching algorithm to two-sided queuing systems.
 This algorithm was first proposed by  \cite{tassiulas1990stability} in the context of communication networks. After that, the max-weight algorithm and the backpressure algorithm, which is a generalization of the max-weight algorithm, are studied intensively in the literature. The book by  \cite{srikantbook} provides an excellent summary. The performance of the max-weight algorithm in the context of a switch operating in heavy traffic has been studied by  \cite{maguluri2016heavy}. 
 The backpressure algorithm was used in the context of online ad matching in \cite{srikantad} and in the context of ride hailing in \cite{kanoria2019backpressure}. 

Heavy traffic analysis of the max-weight algorithm in the context of single-sided queue has a long line of literature. 
One analysis approach is based on fluid limits, diffusion limits and reflected Brownian motion (RBM) \citep{harrison2013brownian}. In this approach, the queueing process is studied under an appropriate scaling and the corresponding limiting fluid or diffusion process is shown to converge to a lower dimensional RBM. This phenomenon is called state space collapse (SSC). If the RBM is single dimensional, then it is called complete resource pooling (CRP). Examples on this line of work to study SSC under the max-weight algorithm in the context of single-sided queues are \citet{williams1998diffusion,stolyar2004maxweight,gamarnik2006validity}. In this paper, we employ another approach based on the Lyapunov drift method developed by \cite{atilla_srikant} and later used by \cite{siva_switch} for switch systems. We generalize the Lyapunov function for two-sided queues and analyzed the max-weight algorithm under the CRP condition similar to that in \cite{CRP_queueing,CRP_manufacturing}.

\subsection{Notation}
Throughout the paper, vectors are denoted by boldface letters. 
Functions applied on vectors are defined entrywise; e.g., $F(\BFlambda)$ is defined to be $(F(\lambda_1), \dots ,F(\lambda_m))$. 
For any two vectors $\BFa \in \mathbb{R}^n$ and $\BFb \in \mathbb{R}^m$, we denote the concatenated vector of dimension $n+m$ by $(\BFa,\BFb)$.  
We denote the $n$-dimensional vector with all 1's by $\BFone_n$, and the $n$-dimensional vector with all 0's by $\BFzero_n$; we omit the subscript $n$ if the sizes of these vectors are clear from the context. 
If $\BFx$ and $\BFy$ are of the same dimension, 
we use $\inner{\BFx}{\BFy}$ to denote the inner product, and $\BFx \circ \BFy$ to denote the Hadamard product (i.e., entrywise product). Any inequality $\BFx \leq \BFy$ is also defined entrywise.
We use the superscript $s$'' to denote variables related to servers and the superscript $c$'' for variables related to customers. We use $\BFe_j^{(c)}$ and $\BFe_i^{(s)}$ to represent unit vectors with a 1 for type $j$ customer and type $i$ server, respectively, and all 0's otherwise.
\section{Model} \label{sec: model}
We represent the types of customers and servers by a bipartite graph $G(N \cup M,E)$, where $N$ is the set of server types with $|N|=n$, $M$ is the set of customers type with $|M|=m$, and $E$ is the set of edges representing customer and server types that are compatible with each other (see Figure~\ref{fig:ex_multiple_link}). A pair $(i,j) \in E$ if and only if a type $j$ customer can be served by a type $i$ server. Each node in the bipartite graph is a queue of customers or servers waiting to be matched with any one of the compatible counterparts.  Our convention is to refer to incoming customers as \emph{demand} and incoming servers as \emph{supply}.
\begin{figure}[!htb]
    \FIGURE
    {\begin{tabular}{c}
        \begin{tikzpicture}[scale=0.7]
\draw[black, very thick] (0,0) -- (2,0) -- (2,1) -- (0,1);
\node[black,very thick] at (0.25,0.5) {2};
\draw[black, very thick] (0,1.5) -- (2,1.5) -- (2,2.5) -- (0,2.5);
\node[black,very thick] at (0.25,2) {1};
\fill[black] (1,-0.3) circle (0.05);
\fill[black] (1,-0.5) circle (0.05);
\fill[black] (1,-0.7) circle (0.05);
\draw[black, very thick] (0,-1) -- (2,-1) -- (2,-2) -- (0,-2);
\node[black,very thick] at (0.25,-1.5) {$n$};
\draw[black,very thick] (8,0) -- (6,0) -- (6,1) -- (8,1);
\node[black,very thick] at (7.75,0.5) {2};
\draw[black, very thick] (8,1.5) -- (6,1.5) -- (6,2.5) -- (8,2.5);
\node[black,very thick] at (7.75,2) {1};
\fill[black] (7,-0.3) circle (0.05);
\fill[black] (7,-0.5) circle (0.05);
\fill[black] (7,-0.7) circle (0.05);
\draw[black, very thick] (8,-1) -- (6,-1) -- (6,-2) -- (8,-2);
\node[black,very thick] at (7.75,-1.5) {$m$};
\draw[black,thick]  (2.75, 2.1) edge[<->]  (5.25, 2.1);
\draw[black,thick]  (2.75, 1.9) edge[<->]  (5.25, -1.4);
\draw[black,thick]  (2.75, 0.6) edge[<->]  (5.25, 1.9);
\draw[black,thick]  (2.75, 0.4) edge[<->]  (5.25, 0.4);
\fill[black] (4,-0.3) circle (0.05);
\fill[black] (4,-0.5) circle (0.05);
\fill[black] (4,-0.7) circle (0.05);
\draw[black,thick]  (2.75, -1.5) edge[<->]  (5.25, -1.5);
\node[black, align=center] at (1,3) {\footnotesize Servers};
\node[black, align=center] at (4,3) {\footnotesize \shortstack{Compatible\\Matchings}};
\node[black, align=center] at (7,3) {\footnotesize Customers};
\end{tikzpicture}
    \end{tabular}}
    {\centering{Bipartite graph representation for two-sided queues.}
    \label{fig:ex_multiple_link}}
    {}
\end{figure}
At each point in time, the system operator posts a price for each customer and server type.
Customers willing to pay the quoted prices, as well as servers who are willing to provide their service at the posted prices (i.e., wages), are admitted to the system. Thus, the system operator can vary the prices to control the arrival rates of customers and servers.
Customers and servers then wait in queues until they are matched. The first-come-first-serve (FCFS) discipline is employed for each queue separately, but it may \emph{not} hold among different types of customers and servers. 
Once a customer is matched with a compatible server, we assume that they depart from the system instantaneously to complete the service process.
The system operator's objective is to find a joint pricing and matching policy under which the system is stable (positive recurrent) and the long-run average profit is maximized. 
We assume that customers and servers arrive according to nonhomogeneous Poisson processes.
For each server type $i \in {N}$, there exists a supply curve $\mu_i: \mathbb{R}_+ \to \mathbb{R}_+$, such that if
the system operator sets a price $p_i^{(s)}$ and the expected waiting time is  $w_i^{(s)}$, the resulting arrival rate is $\mu_i\left(p_i^{(s)} - s^{(s)}_i w_i^{(s)}\right)$, where the constant $s^{(s)}_i > 0$ is the unit waiting cost of server type $i$.
Similarly, for each customer type $j \in {M}$, there exists a demand curve $\lambda_j: \mathbb{R}_+ \to \mathbb{R}_+$, such that if the system operator sets a price $p_j^{(c)}$ and the expected waiting time is $w_j^{(c)}$,
the resulting arrival rate is $\lambda_j\left(p_j^{(c)} + s^{(c)}_{j} w_j^{(c)}\right)$, where $s^{(c)}_{j} > 0$ is the unit waiting cost of customer type $j$.
We make the following assumption on the supply and demand curves.
\begin{assumption} \label{ass: monotonic}
The supply curves $\mu_i: \mathbb{R}_+ \to \mathbb{R}_+$ $(\forall i \in {N})$ are strictly increasing and twice continuously differentiable. The demand curves $\lambda_j: \mathbb{R}_+ \to \mathbb{R}_+$ $(\forall j \in {M})$ are strictly decreasing and twice continuously differentiable. 
\end{assumption}
Since $\lambda_j$  and $\mu_i$ 
are strictly monotone, their inverse functions exist, and we denote them by $F_j: \mathbb{R}_+ \to \mathbb{R}_+$ $(\forall j \in {M})$ and $G_i: \mathbb{R}_+ \to \mathbb{R}_+$ $(\forall i \in {N})$, respectively.
In addition, we define the revenue and cost functions as $r_j^{(c)}(\lambda_j)\overset{\Delta}{=}\lambda_j F_j(\lambda_j)$ for all $j \in {M}$ and $r_i^{(s)}(\mu_i)\overset{\Delta}{=}\mu_iG_i(\mu_i)$ for all $i \in {N}$.
We make the following assumption on the revenue and cost functions.
\begin{assumption} \label{ass: concavity}
The revenue function $r_j^{(c)}(\lambda_j)$ is strictly concave with $(r_j^{(c)})^{\prime \prime}(\lambda_j) < 0$ $(\forall j \in {M})$. The cost function
$r_i^{(s)}(\mu_i)$ is strictly convex with $(r_i^{(s)})^{\prime \prime}(\mu_i) > 0$ $(\forall i \in {N})$.
\end{assumption}
The concavity assumption on revenue function follows from the economic law of diminishing marginal return: 
as the system operator increases the customer arrival rate $\lambda_j$, the marginal revenue ${d r_j^{(c)}(\lambda_j)}/{d\lambda_j}$ decreases, which implies that the revenue function $r_j^{(c)}(\lambda_j)$ is concave. This assumption is often assumed in the revenue management literature \citep{gallego1994optimal,kim2017value}. We assume that the marginal cost ${d r_i^{(s)}(\mu_i)}/{d \mu_i}$ increases with $\mu_i$, since it becomes harder to recruit servers when we try to increase server arrival rate. This implies that the cost function $r_i^{(s)}$ is convex.
For those customers and servers waiting in queues, the system operator uses matching controls to govern the queueing process. At any given time, suppose $q_i^{(s)}$ is the number of type $i$ servers waiting in queue, and
$q_j^{(c)}$ is the number of type $j$ customers waiting in queue.
We denote the vector of all queue lengths by $\BFq = (q_j^{(c)}, \forall j\in M,\  q_i^{(s)}, \forall i\in N)$.
We denote the number of type $i$ servers to be matched to type $j$ customers by $y_{ij}$. The set of feasible matching decisions is
\[
    Y(\BFq) \ \overset{\Delta}{=}\ \left\{ \BFy \in \mathbb{Z}^{nm}_+ \ \Bigm|\ 
    \sum_{i=1}^n y_{ij} \leq q_j^{(c)} \ (\forall j \in {M}),\ 
    \sum_{j=1}^m y_{ij} \leq q_i^{(s)}  \ (\forall i \in {N}), \ 
    y_{ij}= 0 \ (\forall (i,j) \notin E)
    \right\}.
\]
We also define the projection of $Y(\BFq)$ to the queue length space as
\begin{equation}\label{eq:xq}
    X(\BFq) \ \overset{\Delta}{=}\ \left\{ \BFx \in \mathbb{Z}^{n+m}_+ \ \Bigm|\ 
    \exists\  \BFy \in Y(\BFq): \
    x^{(s)}_i = \sum_{j=1}^m y_{ij}
    \ (\forall i \in {N}),\ 
    x^{(c)}_j = \sum_{i=1}^n y_{ij}
    \ (\forall j \in {M})
    \right\}.
\end{equation}
When a pair of customer and server is matched by the system, they both depart from the system. Since a customer is only compatible to a subset of server types,
the system operator may have an incentive to hold some customers or servers in queue in order to achieve better matches in future.
\vspace{12pt}
\emph{Example: Ride Hailing.\ }
An application of the two-sided queueing model is in ride hailing systems. 
In such a system, the customer and server (drivers) types, as well as the matching compatibility graph, are determined by their geographical locations. A simple example with three regions is shown in Figure \ref{fig:ex_ride_hailing}. (Here, we ignore issues such as vehicle capacity and number of seats requested by customers, which can be accounted for by creating additional customer and server types.)
Based on the price and the waiting time quoted to customers, only a fraction of them who open the app will book a ride, which determines the customer arrival rate. Similarly, based on the price quoted to the drivers, they will choose whether or not to provide service. 
Thus, the arrival rates of customer and drivers depend on price and wait time and are governed by the demand and supply curve of each region. 
Once a customer confirms the price and books a ride, the system operator can determine which driver (from what region) should be matched to the customer. 
If a driver accepts the ride request, then they immediately become unavailable for any other ride requests (departing from the system). After the ride is complete, the car becomes available again, possibly in a different region. A simplifying assumption in our model is that we treat a driver who completes the service and re-enters the system the same as a new arrival.
\begin{figure}
  \TABLE{A ride hailing system with three regions.   \label{fig:ex_ride_hailing}}
  {\begin{tabular}[b]{c}
    \begin{tikzpicture}[scale=0.8, every node/.style={scale=0.8}, hexa/.style= {shape=regular polygon,
                                   regular polygon sides=6,
                                   minimum size=2cm, draw,
                                   inner sep=0,anchor=south}]
\node[hexa, thick] (1) at (0,{9*sin(60)}) {1};
\node[hexa, thick] (2) at ({2/2+2/4},{8*sin(60)}) {2};
\node[hexa, thick] (3) at ({4/2+4/4},{7*sin(60)}) {3};
\end{tikzpicture} \\
    \footnotesize (a) Geographical locations
  \end{tabular} \qquad
  \begin{tabular}[b]{c}
    \begin{tikzpicture}[scale=0.7]
\draw[black, thick] (0,0) -- (2,0) -- (2,1) -- (0,1);
\node[black, thick] at (0.25,0.5) {2};
\draw[black, thick] (0,1.5) -- (2,1.5) -- (2,2.5) -- (0,2.5);
\node[black,thick] at (0.25,2) {1};
\draw[black, thick] (0,-0.5) -- (2,-0.5) -- (2,-1.5) -- (0,-1.5);
\node[black, thick] at (0.25,-1) {3};
\draw[black, thick] (8,0) -- (6,0) -- (6,1) -- (8,1);
\node[black,thick] at (7.75,0.5) {2};
\draw[black,  thick] (8,1.5) -- (6,1.5) -- (6,2.5) -- (8,2.5);
\node[black, thick] at (7.75,2) {1};
\draw[black,  thick] (8,-0.5) -- (6,-0.5) -- (6,-1.5) -- (8,-1.5);
\node[black,thick] at (7.75,-1) {3};
\draw[black,thick]  (2.75, 2.1) edge[<->]  (5.25, 2.1);
\draw[black,thick]  (2.75, 1.9) edge[<->]  (5.25, 0.7);
\draw[black,thick]  (2.75, 0.7) edge[<->]  (5.25, 1.9);
\draw[black,thick]  (2.75, 0.5) edge[<->]  (5.25, 0.5);
\draw[black,thick]  (2.75, 0.3) edge[<->]  (5.25, -0.9);
\draw[black,thick]  (2.75, -0.9) edge[<->]  (5.25, 0.3);
\draw[black,thick]  (2.75, -1.1) edge[<->]  (5.25, -1.1);
\node[black, align=center] at (1,3) {\footnotesize Customer};
\node[black, align=center] at (4,3) {\footnotesize \shortstack{Compatible\\Matching}};
\node[black, align=center] at (7,3) {\footnotesize Server};
\end{tikzpicture} \\
    \footnotesize (b) Bipartite graph
  \end{tabular}}
 {We assume that riders can only be matched to cars in their own region or any neighboring regions. The two-sided system generated from the map is shown in subfigure (b).}
\end{figure}
\subsection{Continuous-Time MDP Formulation}
\label{sec:CTMDP}
We now formulate the system operator's decision problem as a continuous-time Markov decision process (CTMDP) and specify its states, actions, transition rates, and rewards.
The system state is represented by the queue lengths of all customer and server types $\BFq \in \mathbb{Z}_{+}^{n+m}$. 
The actions of the CTMDP include both pricing and matching decisions.
The matching decision must satisfy $\BFx \in X(\BFq)$ defined by Eq~\eqref{eq:xq}.
For the pricing decision, in order to leverage Assumption~\ref{ass: concavity}, it is more convenient to use arrival rates $(\BFlambda, \BFmu)$ rather than prices as the control variables.
In particular, for customer type $j \in M$, setting the arrival rate to $\lambda_j$ is equivalent to setting the price to $p^{(c)}_j = F_j(\lambda_j) - s^{(c)}_j w^{(c)}_j(\BFq)$. Similarly, for server type $i\in N$, setting the arrival rate to $\mu_i$ is equivalent to setting the price to $p^{(s)}_i = G_i(\mu_i) + s^{(s)}_i w^{(s)}_i(\BFq)$.
Thus, the action is a tuple $\BFz \overset{\Delta}{=} (\BFlambda,\BFmu,\BFx) \in \mathbb{R}^{2(m+n)}$.
Given this action, the transition rate from state $\BFq$ to state $\BFq + \BFe^{(c)}_j - \BFx$ (i.e., having a new arrival of type $j$ customer) is $\lambda_j$ ($\forall j \in M$), and a reward of $p^{(c)}_j$ is received upon the new arrival. 
The transition rate from state $\BFq$ to state $\BFq + \BFe^{(s)}_i - \BFx$ (i.e., having a new arrival of type $i$ server) is $\mu_i$ ($\forall i \in N$), and a cost of $p^{(s)}_i$ is paid upon the new arrival. 
The system operator's objective is to find a pricing and matching policy such that the long run average profit earned by the system operator is maximized. Note that as $\mathbf{x} \in X(\BFq)$, and we only make matches at the time of an arrival, the state never reaches $\BFzero_{n+m}$ and so we restrict the state space to be $\mathbb{Z}_+^{n+m} \backslash \{\BFzero_{n+m}\}$. We restrict our attention to policies that make the system stable in the long run, which is defined as follows.
\begin{definition}\label{def:stable_policies}
A joint pricing and matching policy is said to be stable, if the continuous-time Markov chain (CTMC) induced by this policy has a positive recurrent communicating class that contains the state $\BFq=e_j^{(c)}$ for some $j \in M$.
\end{definition}
\begin{remark}[Average waiting time.]
It is technically challenging to analyze the exact waiting time
$ w^{(s)}_i(\BFq)$ and $w^{(c)}_j(\BFq)$, since the waiting time of one type may depend on the queue lengths of all the types as well as the policy and matching policy used by the system operator. 
Additionally, in some applications,  real-time queue length information may not be visible to all market participants \citep{zohar2002adaptive}.
Therefore, we make a simplifying assumption that the waiting time \emph{perceived} by the customers and servers is the long-run average waiting time. That is, we assume 
\[
p^{(c)}_j = F_j(\lambda_j) - s^{(c)}_j \bar{w}^{(c)}_j \ \  \forall j\in {M},
\quad
p^{(s)}_i = G_i(\mu_i) + s^{(s)}_i \bar{w}^{(s)}_i \ \ \forall i\in {N},
\]
where $\bar{w}^{(c)}_j$ is the average of the waiting times of type-$j$ customers. Similarly, $\bar{w}^{(s)}_i$ is the average of the waiting times of type-$i$ servers. The scheme of announcing the long-run average waiting time to (impatient) customers is commonly assumed in the literature \citep{zohar2002adaptive,armony2009impact}.
Additionally, in the large scale setting that will be considered in the following sections, approximating real-time estimated waiting time with the long-run average waiting time will only result in a negligible error term of a higher order (see \cite{kim2017value}, Section 6.1 for a similar argument).
\end{remark}
\subsubsection{Equivalence to Holding Cost Models.}
The above model assumes that customers and servers are sensitive to both prices and waiting costs when they decide to enter the queueing system.
We now consider an alternative model, where customers and servers only react to prices, while the system operator pays  \emph{additional} holding costs for market participants waiting in queues. In particular, in this alternative model, the states, actions, and transition rates remain the same. Given a state $\BFq$ and an action $\BFz = (\BFlambda, \BFmu, \BFx)$, the reward function is defined as
\begin{equation}\label{revenue}
        \mathcal{R}(\BFq, \BFz)\overset{\Delta}{=} \sum_{j=1}^m \lambda_j F_j(\lambda_j)
    - \sum_{i=1}^n \mu_i G_i(\mu_i)
    - \sum_{j=1}^m s^{(c)}_j q^{(c)}_j
    - \sum_{i=1}^n s^{(s)}_i q^{(s)}_i,
\end{equation}
where $s^{(c)}_j$ and $s^{(s)}_i$ are the customers' and servers' impatience parameters introduced in the original model.
The following result shows that the two modelling approaches are indeed equivalent. 
\begin{proposition}
\label{prop: model-equivalence}
For any given stationary and ergodic policy with finite average waiting times and queue lengths,
the delay-sensitive model and the holding cost model have the same long-run average profit.
\end{proposition}
The proof of Proposition~\ref{prop: model-equivalence} follows an application of Little' Law and can be found in Appendix~\ref{sec:proof_of_equivalence}.
The advantage of considering the holding cost model is that the reward function $\mathcal{R}(\BFq, \BFz)$ does not explicitly depend on the waiting time.
Hence, we use the holding cost model in the rest of the paper.
\subsection{Discrete-Time MDP Formulation by Uniformization}
Instead of analyzing the CTMDP directly, we use the well-known uniformization technique \citep[e.g.][Chap.~11]{puterman2014markov} to obtain an equivalent discrete-time Markov decision process (DTMDP), which will simplify our analysis. 
The uniformized process works as follows.
We first choose a uniformization parameter $c$ defined below.
\begin{definition} \label{ass: upperboundedrates}
Suppose there exist constants $\BFlambda_{\max} \in \mathbb{R}^{m}_+$ and $\BFmu_{\max} \in \mathbb{R}^{n}_+$ such that for any price vector $\BFp \geq 0$ we have,
$   \BFlambda(\BFp) \leq \BFlambda_{\max}$, and $\BFmu(\BFp) \leq \BFmu_{\max}$.
Let $c$ be any constant such that 
$   c \geq \inner{\BFone_m}{\BFlambda_{\max}}+\inner{\BFone_n}{\BFmu_{\max}}. $
\end{definition}
The uniformized DTMDP is endowed with the same state $\BFq$ and action $\BFz = (\BFlambda, \BFmu, \BFx)$ as the CTMDP. Let $Z(\BFq) = [0, \BFlambda_{\max}] \times [0, \BFmu_{\max}] \times X(\BFq) $ be the set of feasible actions for queue length $\BFq \in \mathbb{Z}^{m+n}_+$. 
In the uniformized DTMDP, there is at most one customer arrival or one server arrival in each period: the state transitions from $\BFq$ to $\BFq + \BFe^{(c)}_j - \BFx$ with probability $\lambda_j / c$ ($\forall j \in M$); it transitions from $\BFq$ to $\BFq + \BFe^{(s)}_i - \BFx$ with probability $\mu_i / c$ ($\forall i \in N$); otherwise, no arrival happens in this period, and the state  remains at $\BFq$ with probability $1-(\inner{\BFone_m}{\BFlambda} + \inner{\BFone_n}{\BFmu})/c$.
The expected reward in one period is given by $\mathcal{R}(\BFq, \BFz)/c$.
Let $\BFq'$ be the state in the next period.
The Bellman equation of the DTMDP is
\begin{align}
    h(\BFq)+\frac{\gamma}{c}\ =\  \max_{\BFz \in Z(\BFq)} \left\{\frac{\mathcal{R}(\BFq,\BFz)}{c} + \E[h(\BFq')\mid \BFq,\BFz] \right\}, \quad \forall \BFq \in \mathbb{Z}^{n+m}_+ \backslash \{\BFzero_{n+m}\}, \label{MDP} 
    \end{align}
    where
    \begin{align}
    \E[h(\BFq')\mid \BFq,\BFz]\ =\ {}&\sum_{j=1}^m\frac{\lambda_j}{c}h(\BFq+\BFe_j^{(c)}-\BFx)
    +\sum_{i=1}^n\frac{\mu_i}{c}h(\BFq+\BFe_i^{(s)}-\BFx) \nonumber \\
    &\quad +(1-\sum_{j=1}^m \frac{\lambda_j}{c}-\sum_{i=1}^n \frac{\mu_i}{c})h(\BFq). \label{yqplus1} 
\end{align}
In the above equation, the solution $\gamma$ is the optimal long-run average profit, and $h(\BFq)$ is the bias function associated with state $\BFq$ ($\forall \BFq \geq 0$). The term $\E[h(\BFq')\mid \BFq,\BFz]$ is the expectation of the bias function $h$ after one transition in the uniformized process.
The expectation is taken with respect to the one-period transition probabilities conditional on the state $\BFq$ and the action $\BFz$.

We invoke \citet[Theorems 3.1 and 3.2, Lemma 3.1]{cavazos1989weak} to conclude that stationary optimal policy exists and the Bellman optimality equation \eqref{MDP} exhibits a solution. \citet[Theorems 3.1 and 3.2, Lemma 3.1]{cavazos1989weak} requires the following four conditions, which are all satisfied:
\begin{itemize}
    \item \citet[Assumption 1.1]{cavazos1989weak} holds as the transition matrix and the reward 
    are continuous functions of the actions $(\BFlambda, \BFmu, \BFx)$. Moreover, the action set is compact.
    \item \citet[Assumption 1.2]{cavazos1989weak} follows as one can consider an equivalent cost minimization formulation of our MDP with cost function $\mathcal{C}(\BFq, \BFz) = K - \mathcal{R}(\BFq, \BFz)$. With a large enough constant $K > 0$, $\mathcal{C}(\BFq, \BFz)$ is non-negative due to the compact action space. Moreover, $\mathcal{C}(\BFq, \BFz) \to \infty$ as $\|\BFq\|_1 \to \infty$ due to positive holding costs. 
    \item \citet[Assumption 1.3]{cavazos1989weak} requires existence of a stationary deterministic policy such that the corresponding Markov chain is irreducible on $\mathbb{Z}_+^{m+n} \backslash \{\BFzero_{n+m}\}$, positive recurrent, and has a finite long-run average reward. As we have full freedom to control the arrival rates, it is not too hard to construct such a policy. 
    
    Retain a spanning subgraph whose connected components are complete
bipartite, and let $Z_\ell$ denote the customer--server imbalance in
component $\ell$. Choose $0<a<b$. If $Z_\ell>0$, assign total customer
and server arrival rates $(a,b)$ in component $\ell$ (divided equally across all types); if $Z_\ell<0$,
assign $(b,a)$; and if $Z_\ell=0$ do the following:
\begin{itemize}
    \item If all queue lengths in component $\ell$ are zero, set arrival rates of all types to be $\epsilon > 0$.
    \item If all imbalances are zero ($Z_l = 0$ for all $l$), set the arrival rate of all types to be $\varepsilon>0$ and, upon the next arrival, match all previously queued items; the new
arrival remains, so the all-zero state is never reached
    \item If the queue lengths in component $\ell$ is non-zero, and the overall imbalance vector is nonzero ($\exists l^\prime$ s.t. $Z_{l^\prime} \neq 0$), set all rates in component $\ell$ equal to zero.
\end{itemize}
The reachability of any state $\mathbb{Z}_+^{m+n} \backslash \{\BFzero_{m+n}\}$ is easy to check, which ensures irreducibility. Moreover, each imbalance $Z_l$ has a negative drift as long as $Z_l \neq 0$. Once any of the imbalance reaches $0$ (with non-zero queue lengths), it stays at zero until the entire imbalance vector reaches $0$. Thus the imbalance vector has finite return times to zero, which is then followed by clearing out all the queues. This yields positive recurrence and finite
stationary expected queue length.
    \item The technical condition \citet[Assumption 3.1]{cavazos1989weak} holds by \citet[Lemma 2.1]{cavazos1989weak} and noting that the MDP has a finite number of possible transitions from any state.
\end{itemize}
Thus, the optimality equation formulation \eqref{MDP} is well-defined. 
In Appendix~\ref{sec:mdp-appendix}, we present additional analysis of the uniformized DTMDP. We  show the monotonicity structure of the optimal pricing policy in the single-link queueing system (i.e., $m=n=1$). Unfortunately, as the number of customer and server types becomes large, solving the DTMDP becomes intractable due to the curse of dimensionality. We propose two approximation methods to obtain near optimal solutions to the DTMDP. 
The first method is based on fluid approximation.
The remainder of the paper primarily focuses on this approach.
The second method uses value function approximation. We defer details of the second method to Appendix~\ref{sec:mdp-appendix}, as the remaining parts of the paper do not rely on it.
\subsection{Max-Weight Matching Policy}
\label{subsec:max-weight matching}
In the following sections, we will extensively use the \emph{max-weight} matching policy, so we provide its definition here. 
Suppose the system has state $\BFq$ and the set of feasible matches is $X(\BFq)$ (see Eq~\eqref{eq:xq}).
The policy chooses the matching decision $\BFx$ to be the solution of
\begin{equation}
    \argmax_{\BFx \in X(\BFq)}\left\{ 
     \inner{\BFq}{\BFx} 
      \right\}. \label{eq:max-weight}
\end{equation}
In other words, under the max-weight policy, when there is either a customer or a server arrival, a match will be made if any of the compatible types has a nonempty queue, and we will always match the arriving customer/server to the compatible type with the most number of customers/servers waiting in queue. Otherwise, if all the compatible counterparts' queues are empty, then the arrival is inserted into the queue of its own type.
The max-weight matching policy, originally proposed by \cite{tassiulas1990stability}, is extensively studied in the queueing literature. This literature is reviewed in Section~\ref{sec:liter_review}. 
Apart from the queueing literature, 
in the our model specifically,
there is also an alternative way to motivate the max-weight matching policy through quadratic value function approximation of the MDP. Suppose the bias function in Eq~\eqref{yqplus1} is approximated by $h(\BFq)\approx \inner{1}{\BFq^2}$, then the \emph{optimal} matching policy of the DTMDP will be very close to a max-weight policy define in Eq~\eqref{eq:max-weight}.
Appendix~\ref{sec: LP approx} contains a detailed discussion of the value function approximation method.
\begin{algorithm} 
\caption{Max-Weight Matching Policy} \label{alg: maxweight}
{\fontsize{10}{13}\selectfont
\begin{algorithmic}
\STATE \textbf{input:} current queue length $\BFq(k)$, new arrival $\BFa(k)$ \COMMENT{$k$ is a decision epoch}
\STATE \textbf{initialization:} $\BFy(k)=\BFzero$
    \FOR{$i \in {N}$}
    \IF {$a_i^{(s)}(k)=1$ \textbf{and} $\max_{j: (i,j)\in E} q_{j}^{(c)}>0$}
        \STATE choose $j^* \in \arg\max_{j: (i,j)\in E} q_{j}^{(c)}$ (breaking ties arbitrarily)
        \STATE set $y_{ij^*}(k)=1$
    \ENDIF
    \ENDFOR
    \FOR{$j \in {M}$}
    \IF {$a_j^{(c)}(k)=1$ \textbf{and} $\max_{i: (i,j) \in E} q_{i}^{(s)}>0$}
        \STATE let $i^* \in \arg\max_{i: (i,j) \in E} q_{i}^{(s)}$ (breaking ties arbitrarily)
        \STATE set $y_{i^*j}(k)=1$
    \ENDIF
    \ENDFOR
\STATE \textbf{output:} matching decision $\BFy(k)$
\end{algorithmic}
}
\end{algorithm}

For simplicity, we allow the max-weight matching policy (Algorithm~\ref{alg: maxweight}) to match an incoming arrival straightaway, which is slightly incompatible with the MDP definition \eqref{MDP} as the Bellman equation does not allow the incoming customer to be matched. One can always consider a system that schedules the matches exactly as governed by Algorithm~\ref{alg: maxweight}, but executes them with a one-step delay. In this new system, the reward $\mathcal{R}(\BFq, \BFz)$ is simply off by an additive factor of $2 \max\{\BFs\}$.

\section{Asymptotically Optimal Policies in the Large-Scale Regime}
\label{sec:large-theta}
\subsection{Fluid Model and Large-Scale Regime} 
\label{sec:fluid}
We consider a fluid approximation of the queueing system where random arrivals are replaced by deterministic arrival processes.
The fluid model is a deterministic optimization problem maximizing the long-run average profit. Suppose customers arrive with constant rates $\BFlambda$ and servers arrive with constant rates $\BFmu$. 
Let $\chi_{ij}$ be the average rate of type $i$ server matched to the type $j$ customer for all $(i,j) \in E$. 
The fluid model is defined as
\begin{subequations}\label{eq:fluid_opt}
\begin{align}
   {\gamma}^{*} = \max_{\BFlambda, \BFmu,\BFchi}\ & \ \inner{F(\BFlambda)}{\BFlambda}-\inner{G(\BFmu)}{\BFmu} \label{fluidobj}\\
    \text{subject to}\ &\ {\lambda}_j = \sum_{i=1}^n \chi_{ij}, \quad \forall j \in {M}, \label{ball} \\
    &\  {\mu}_i = \sum_{j=1}^m \chi_{ij}, \quad \forall i \in {N}, \label{balm} \\
    &\ \chi_{ij}  = 0, \ \  \forall (i,j) \notin E,\ \  \chi_{ij} \geq 0, \ \  \forall (i,j) \in E. \label{consx}
\end{align}
\end{subequations}
We denote an optimal solution to the above fluid problem by $(\BFlambda^*, \BFmu^*, \BFchi^*)$. We assume $\BFlambda^*, \BFmu^* > 0$. For example, if one of the arrival rates is zero, we could simply remove that vertex from the graph $G$. To interpret the fluid model above, note that Eqs \eqref{ball} and \eqref{balm} are the balance equations for the number of customers and servers matched. Eq \eqref{consx} specifies that matching is only allowed among compatible customer-server pairs.
Intuitively, it is easy to see that these constraints are necessary, because if the balance equations do not hold, then some customer or server types will keep accumulating over time. 
Thus, the optimization program Eq~\eqref{eq:fluid_opt} serves as an \emph{upper bound} on the achievable profit under any pricing and matching policy that makes the system stable.
This is formally shown in the following proposition.
The proof can be found in Appendix \ref{app: fluid}.
\begin{proposition}\label{theo: fluid}
The optimal value of the fluid problem Eq~\eqref{eq:fluid_opt} is an upper bound on the long run expected profit rate under any stationary, ergodic policy that makes the system stable.
\end{proposition} 
In the remainder of this section, we analyze the two-sided queueing system in a large-scale regime where the arrival rates of all customer and server types are simultaneously scaled by a factor of $\eta \in \mathbb{N}$.
\begin{definition}[Large-Scale Regime] \label{def: asymptoticregime}
Consider a family of two-sided queueing systems associated with the same bipartite graph $G(N \cup M, E)$ parametrized by $\eta \in \mathbb{N}$. For the $\eta^{\text{th}}$ system, the demand and supply curves satisfy $F^{\eta}(\eta\BFlambda)=F(\BFlambda)$ for all $\BFzero_m \leq \BFlambda \leq \BFlambda_{\max}$ and $G^{\eta}(\eta\BFmu)=G(\BFmu)$ for all $\BFzero_n \leq \BFmu \leq \BFmu_{\max}$. 
\end{definition}
The above large-scale regime is commonly assumed in the dynamic pricing and matching literature \citep{matchingqueues,amywardpricingmatching}, which ensures that supply and demand are balanced as the system scales up.
According to Definition~\ref{def: asymptoticregime}, it is easily verified that the fluid solution to the $\eta^{th}$ scaled system is given by $\eta\olam$ and $\eta\omu$, where $\olam$ and $\omu$ is the optimal solution of the unscaled fluid model Eq~\eqref{eq:fluid_opt}.
\begin{definition}[Profit Loss]\label{def:profit_loss}
The profit loss (denoted by $\lr$) of a policy is the difference between the optimal value of the (scaled) fluid model, denoted by ${\gamma}^{\eta}_{*}$, and the long run average profit (including the penalty incurred due to waiting) under that policy. 
\end{definition}
 The optimal value of the $\eta^{th}$ fluid model is ${\gamma}^{\eta}_{*}=\eta {\gamma}^{*}$. Therefore, if the profit loss of a policy is sublinear in $\eta$, namely $\lr=o(\eta)$, we say the policy is asymptotically optimal in the large-scale regime.
\subsection{Fluid Pricing Policy} \label{subsec: fluidsolution}
Based on the fluid model, we propose a static pricing policy defined as follows:
\begin{align}
    \lambda_j(\BFq)&=\begin{cases}
    \eta\lambda^{*}_j \ &\text{ if } q_j^{(c)} < q_{\max}^{\eta} \\
    0\; &\text{  otherwise }
    \end{cases} \quad \forall j \in {M}, \label{eq:fluid_pricing} \\
    \mu_i(\BFq)&=\begin{cases}
    \eta\mu^{*}_i \ &\text{ if } q_i^{(s)} < q_{\max}^{\eta} \\
    0\; &\text{   otherwise }
    \end{cases} \quad \forall i \in {N}. \nonumber
\end{align}
Here, $q_{\max}^{\eta}$ denotes the maximum queue buffer size; it is a parameter that depends on $\eta$, which will be specified later.
The main intuition of the fluid pricing policy is the following.
When all queues are below their maximum buffer capacity $\mathbf{q}^{\eta}$,
the profit rate of the fluid pricing policy is exactly equal to $\eta {\gamma}^{*}$.
If any customer queue is full, say, $q_j^{(c)}=q_{\max}^{\eta}$, then all future arrivals to queue $j$ will be rejected until at least one customer waiting in queue $j$ is matched. Thus, a fraction of revenue is lost due to customer rejections.
More specially, let ${\gamma}^{\eta}$ be the long run average profit of the fluid pricing policy (excluding waiting costs). Let $\mathbf{I}^{(s)}(q_{\max}^{\eta})$ be a (vector) indicator function representing whether server queues are at the maximum capacity, and let $\mathbf{I}^{(c)}(q_{\max}^{\eta})$ be a (vector) indicator function representing whether customer queues are at the maximum capacity. 
Then, we have
\begin{align}
    \lr={}&\fr-({\gamma}^{\eta}-\inner{\BFs}{\E[\BFq]}) \nonumber \\
    ={}& \eta\bigl(\inner{F(\olam)}{\olam}-\inner{G(\omu)}{\omu}\bigr)-\inner{F(\olam)}{\eta\olam\circ(\BFone-\E[\mathbf{I}^{(c)}
    (q_{\max}^{\eta})])}\nonumber \\
    &\quad -\eta\inner{G(\omu)}{\omu\circ(\BFone-\E[\mathbf{I}^{(s)}(q_{\max}^{\eta})])}+\inner{\BFs}{\E[\BFq]} \nonumber \\
    ={}& \eta\left(\inner{F(\olam)}{(\olam\circ\E[\mathbf{I}^{(c)}(q_{\max}^{\eta})])}-\inner{G(\omu)}{(\omu\circ\E[\mathbf{I}^{(s)}(q_{\max}^{\eta})])}\right)+\inner{\BFs}{\E[\BFq]}, \label{eq: rl}
\end{align}
where the first equality follows from Definition~\ref{def:profit_loss}, and the second equality uses the definition of the fluid pricing policy.
As a result, Eq~\eqref{eq: rl} shows that the profit loss of the fluid pricing policy depends on the parameter $q_{\max}^{\eta}$.
If we increase the buffer capacity $q_{\max}^{\eta}$, then the probability of dropping customers/servers will reduce, i.e., $\E[\mathbf{I}(\BFq_{\max}^{\eta})]$ will decrease. However, increasing the buffer capacity will lead to increasing in the expected queue lengths, which will increase the penalty incurred due to waiting. Thus, we choose buffer capacity to balance the trade-off in order to minimize the overall profit loss. {Precisely, we will see that choosing $q_{\max}^{\eta} \sim \sqrt{\eta}$ will result in $\E[\mathbf{I}(\BFq_{\max}^{\eta})] \sim \eta^{-1/2}$ and $\E[\inner{\BFone_{m+n}}{\BFq}] \sim \sqrt{\eta}$, which attains the optimal profit loss. 
}
\begin{theorem} \label{theo: fluidpricingpolicy}
Suppose a family of two-sided queues is given by the bipartite graph $G(N \cup M,E)$ parameterized by $\eta$.
The profit loss $\lr$ under the fluid pricing policy Eq~\eqref{eq:fluid_pricing} and max-weight matching (Algorithm~\ref{alg: maxweight})
is $O(\sqrt{\eta})$, where $q_{\max}^{\eta}=\con\sqrt{\eta}$ for any positive constant $\con$.
\end{theorem}
The proof of Theorem~\ref{theo: fluidpricingpolicy} can be found in Appendix~\ref{appendix: fluidpricingpolicy}.
In addition, it can be shown that the $O(\sqrt{\eta})$ profit loss rate cannot be improved using any fluid pricing policy. 
The proof of the proposition below is presented in Appendix \ref{appendix: singlelinktwosidedqueue}. 
\begin{proposition} \label{prop: singlelinktwosidedqueue}
For a family of single-link two-sided queues parametrized by $\eta$, any fluid pricing policy will have a profit loss $\lr$ that is at least  $\Omega(\sqrt{\eta})$.
The choice of $q_{\max}^{\eta}=\gamma \sqrt{\eta}$ for any positive constant $\gamma$ provides the optimal profit loss rate $\Theta(\sqrt{\eta})$.
\end{proposition} 
\subsection{Two-Price Policy}
\label{sec:two-price}
A main drawback of the fluid pricing policy is that the prices are not adaptive to changes in the system state.
In this section, we consider another policy that uses \emph{two} different prices for each customer/server type. The proposed two-price policy is built on the two-price policy in \cite{kim2017value} for single server queues. Our contribution lies in a joint analysis of two-price and dynamic matching policies in a multi-type queueing network.
The two-price policy can be viewed as a generalization of the fluid pricing policy.
We introduce additional parameters $\BFtheta \in \mathbb{R}^m_+$, $\BFphi \in \mathbb{R}^n_+$ and $\sigma^{\eta}>0$, which governs the arrival rates of the customers and servers respectively when the queue length is greater than a certain threshold $\tau^{\eta}_{\max}$.
The two-price policy is defined as
\begin{align}
    \lambda_j(\BFq)&=
\begin{cases}
\eta \lambda_j^{*} &\textit{ if } q_j^{(c)} \leq \tau^{\eta}_{\max}\\
\eta \lambda_j^{*}-\theta_j\sigma^{\eta} &\textit{ otherwise }
\end{cases} \quad \forall j \in {M}, \label{eq:two-price_policy-cust}\\
\mu_i(\BFq)&=
\begin{cases}
\eta \mu_i^{*} &\textit{ if } q_i^{(s)} \leq \tau_{\max}^{\eta}\\
\eta \mu_i^{*}-\phi_i\sigma^{\eta} &\textit{ otherwise }
\end{cases} \quad \forall i \in {N}. \nonumber
\end{align}
The policy sets a threshold $\tau_{\max}^{\eta}$ for all customer and server types. It uses the fluid arrival rates when queue lengths are below this threshold, and then reduces the arrival rates by $\theta_j\sigma^{\eta}$ outside this threshold for type $j$ customer. Similarly, the policy reduces the server arrival rates outside the threshold by $\phi_i\sigma^{\eta}$ for type $i$ server. Here, $\tau_{\max}^{\eta}$, $\sigma^{\eta}$, $\BFtheta$ and $\BFphi$ are parameters that will be specified later. (Our convention is to use superscript $\eta$ to denote  any parameter or quantity that is associated with the $\eta^{th}$ scaled system.)
Intuitively, for any type of customer/server, if we increase $\sigma^{\eta}$, the queue length will have a larger negative drift when it exceeds the threshold $\tau_{\max}^{\eta}$, so 
the expected queue length $\E[\inner{\BFone_{m+n}}{\BFq}]$ will be smaller. However, if $\sigma^{\eta}$ are too large, 
the arrival rates outside the threshold $\tau_{\max}^{\eta}$ will be far from the optimal fluid arrival rates, which will result in a larger profit loss. Thus, there is a trade-off between the expected queue length and profit loss. 
For the matching algorithm associated with the two-price policy, here we use a modified version of max-weight matching algorithm (Eq~\eqref{eq:max-weight}). First, we define redundant edges $(i, j) \in E$ as follows:
\begin{definition} \label{def:redundant_edges}
    An edge $(i, j) \in E$ is redundant if and only if $\chi^*_{ij} = 0$ for all optimal solutions of \eqref{eq:fluid_opt}. Denote the set of redundant edges by $E_r \subseteq E$.
\end{definition}
Now, we consider a modified max-weight matching policy that mimics the max-weight matching policy but matches only using the non-redundant edges, i.e., replacing $E$ by $E \backslash E_r$ in Algorithm~\ref{alg: maxweight} (Other matching algorithms will be considered in Section~\ref{sec:MW vs Rand}).
The following theorem provides a bound on the asymptotic performance of the two-price policy as
$\eta$ tends to infinity.
\begin{theorem} \label{theorem: twoprice}
Consider a family of two-sided queues parametrized by $\eta$ represented by the bipartite graph $G(N \cup M,E)$. The profit loss $\lr$ under the two-price policy Eq~\eqref{eq:two-price_policy-cust} and the modified max-weight matching (Algorithm~\ref{alg: maxweight}) is $O(\eta^{1/3})$ for any $\tau_{\max}^{\eta} \leq \eta^{1/3}$, $\sigma^{\eta}=\eta^{2/3}$ and constants $\BFtheta>\BFzero_{m}$, $\BFphi>\BFzero_{n}$.
\end{theorem}
The above theorem shows that the profit loss of the two-price policy is $O(\eta^{1/3})$, which is better than the $O(\sqrt{\eta})$ loss in the fluid pricing policy. 
The proof of the theorem contains two main steps. 
The first step is to show that the system is stable under the two-price policy and the expected queue lengths are bounded.  We also give an upper bound of the expected queue lengths (Lemma~\ref{lemma: posrec}). 
The second step in the proof is to estimate the profit loss $\lr$ (Lemma~\ref{lemma: firstorderterms}) by applying the KKT conditions of the fluid problem.
\begin{lemma} \label{lemma: posrec}
For a system of two-sided queues operating under the two-price policy and the modified max-weight matching algorithm parameterized by $\eta$, the system is positive recurrent for any $\BFtheta>\BFzero_m$, $\BFphi>\BFzero_n$, $\sigma^{\eta}>0$ and $\tau_{\max}^{\eta}>0$. The expected queue lengths are bounded by
\begin{align*}
    \E \left[\inner{\BFtheta}{\BFq^{(c)}}\right]+\E \left[\inner{\BFphi}{\BFq^{(s)}}\right]  \leq &
    \tau_{\max}^{\eta}\left(\sum_{j = 1}^{m} \theta_j  \Pr[q_j^{(c)}\leq \tau_{\max}^{\eta}]+ \sum_{i =1}^{n} \phi_i \Pr[q^{(s)}_i \leq \tau_{\max}^{\eta}]\right) \nonumber\\
    &+ \frac{\eta}{\sigma^{\eta}}\left(\inner{\BFone_n}{\omu} + \inner{\BFone_m}{\olam}\right). 
\end{align*}
\end{lemma}
\begin{lemma} \label{lemma: firstorderterms}
For a system of two-sided queues operating under two-price policy and modified max-weight matching policy, for any $\BFtheta>\BFzero_m$, $\BFphi>\BFzero_n$ and $\tau_{\max}^{\eta}>0$, we have
\begin{align*}
&\sum_{j \in {M}} \left(F_j'(\lambda_j^{*})\lambda_j^{*}+F_j(\lambda_j^{*})\right)\theta_j \Pr[q^{(c)}_j>\tau_{\max}^{\eta}]-
\sum_{i \in {N}} \left(G_i'(\mu_i^{*})\mu_i^{*}+G_i(\mu_i^{*})\right)\phi_i \Pr[q^{(s)}_i>\tau_{\max}^{\eta}] = 0.
\end{align*}
\end{lemma}
\subsection{Lower Bound}
\label{sec:lower_bound}
In this section, we will obtain a lower bound on the profit loss under a broad family of policies, and thus establish that the $O(\eta^{1/3})$ rate obtained by the two-price policy in Theorem~\ref{theorem: twoprice} is optimal. 
In particular, we consider a family  of pricing policies that have the following form:
\begin{align}
    \lambda_j&=\eta\lambda_j^{*}+f_j\left(\frac{\BFq}{\eta^{\alpha}}\right)\eta^{\beta} \quad \forall j \in {M}, \label{eq: pricinggraphgeneralcust}\\
     \mu_i&=\eta\mu_i^{*}+g_i\left(\frac{\BFq}{\eta^{\alpha}}\right)\eta^{\beta} \quad \forall i \in {N}. \label{eq: pricinggraphgeneralser}
\end{align}
The motivation for this policy is as follows.
The first terms in Eqs~\eqref{eq: pricinggraphgeneralcust} and \eqref{eq: pricinggraphgeneralser} (i.e., $\eta\lambda_j^{*}$ and $\eta\mu_i^{*}$) are static and result from the solution of the fluid model; the second terms account for dynamic adjustments as the queue length changes.  We assume the adjustment terms can be further decomposed into two terms: a function that rescales the queue length, $f_j(\cdot)$ or $g_i(\cdot)$, and a term that determines the scaling of price adjustments, $\eta^{\beta}$, for some $1>\beta>0$.
As the arrival rates are scaled up, the average queue length will also increase. Thus, we rescale the queue length in the functions $f_j(\cdot)$ and $g_i(\cdot)$ for all $i \in {N}$ and $j \in {M}$ by $\eta^{\alpha}$ for some $1\geq \alpha \geq 0$. 
For our analysis, we assume the pricing policy to 
 satisfy the following conditions.
\begin{condition}  \label{ass: generalpricinggraph} 
\begin{enumerate}[label=(\alph*)]
    \item  There exist constants ${\BFGamma} \in \mathbb{R}_{+}^{m}$ and ${\BFPsi} \in \mathbb{R}_{+}^{n}$ such that $|f_j\left({\BFq}/{\eta^{\alpha}}\right)| \leq {\Gamma}_j$ for all $j \in {M}$ and $|g_i\left({\BFq}/{\eta^{\alpha}}\right)| \leq {\Psi}_i$ for all $i \in {N}$ for all $\BFq \in S$ and for all $\eta \geq 1$. \label{condition: boundedfg}
    \item  $ 0 < \alpha+\beta \leq 1$. \label{condition: aplusb}
    \item There exist constants $\kappa>0$ and $\delta>0$ such that for all $j \in {M}$, if ${q_j^{(c)}}/{\eta^{\alpha}}>\kappa$, then either $f_j\left({\BFq}/{\eta^{\alpha}}\right)<-\delta$ or there exists $i: (i,j) \in E$ such that $g_i\left({\BFq}/{\eta^{\alpha}}\right)>\delta$ for all $\eta$. Similarly for all $i \in {N}$, if ${q_i^{(s)}}/{\eta^{\alpha}}>\kappa$, then either $g_i\left({\BFq}/{\eta^{\alpha}}\right)<-\delta$ or there exists $j: (i,j) \in E$ such that $f_j\left({\BFq}/{\eta^{\alpha}}\right)>\delta$ for all $\eta$. \label{condition: posrecsystem}
\end{enumerate}
\end{condition}
We now interpret the conditions above.
Condition \ref{ass: generalpricinggraph}\ref{condition: boundedfg} requires the functions $f$ and $g$ to be bounded, given appropriately scaling of the queue lengths $\BFq$ as $\eta$ increases.
Condition \ref{ass: generalpricinggraph}\ref{condition: aplusb} states that the rate of queue length re-scaling ($\alpha$) should not exceed the rate of re-scaling pricing adjustment terms ($1-\beta$). This condition is needed so that the price adjustment terms are sufficiently large to make the system stable. (In the special case of a single-link system, this assumption is not needed; the extension is presented later in Proposition~\ref{prop: twopricesinglelinktwosidedqueue}.)  Condition \ref{ass: generalpricinggraph}\ref{condition: posrecsystem} states that
if a queue is too long, we should either decrease the arrival rate of this queue or increase the arrival rates of those matched to this queue.
Aside from the above conditions, the pricing forms in Eqs (\ref{eq: pricinggraphgeneralcust}-\ref{eq: pricinggraphgeneralser})
are fairly general, because the 
pricing function of any queue can depend on the entire system state vector $(\BFq)$, and we do not make any strong assumptions such as monotonicity, continuity or differentiablity on functions $f$ and $g$.
Finally, we emphasize that our analysis does not require any assumption on the form of matching policies.
The two-price policy in Section~\ref{sec:two-price} satisfies the above condition with
\begin{align}
    f_j(\BFq)=-\theta_j\mathbbm{1}_{q_j^{(c)}>\tau_{\max}} \ (\forall j \in {M}), \quad
    g_i(\BFq)=-\phi_i\mathbbm{1}_{q_i^{(s)}>\tau_{\max}} \ (\forall i \in {N}), \quad \beta=2/3. \label{eq: twopricecondition12}
\end{align}
Now we present the result on the lower bound. 
\begin{theorem} \label{theorem: lowerboundmultiplelink}
For a two-sided queue defined by a graph $G(N \cup M,E)$ operating under any pricing policy of the form Eq~\eqref{eq: pricinggraphgeneralcust} and \eqref{eq: pricinggraphgeneralser} that satisfies Condition \ref{ass: generalpricinggraph}, if the resulting system is stable, there exists a constant $K(F,G,f,g)$ such that
\begin{align*}
    \lr \  \geq\  K\eta^{1/3}.
\end{align*}
\end{theorem}
The details of the proof are deferred to Appendix \ref{appendix: lowerboundmultiplelink}. We present below an intuitive explanation of the rate in the lower bound.
\begin{remark}[Intuitive explanation of $\eta^{1/3}$]
\label{rmk:eta-lb}
The main reason why the profit loss lower bound is of order $\Omega(\eta^{1/3})$ is due to the trade-off between the expected queue length and the loss in revenue. Consider a pricing policy that deviates from the fluid optimal pricing policy by  $\epsilon>0$, that is, for all $\BFq \in S$, we have $|\lambda_j(\BFq)-\lambda_j^*|<\epsilon$ for all $j \in {M}$ and $|\mu_i(\BFq)-\mu_i^*|<\epsilon$ for all $i \in {N}$. 
One can show that under such a policy, the expected queue length is of the order ${1}/{\epsilon}$ and revenue loss is of the order $\eta\epsilon^2$.
Specifically, the queue length can be coupled to that of an $M/M/1$ queue in heavy-traffic with parameter $\epsilon$, whose mean queue length is known to be of  the order ${1}/{\epsilon}$ by the Kingman's bound. 
The loss in revenue can be estimated by the Taylor series expansion of the revenue function. Since we are operating closed to the optimal price of the fluid model, the first order term vanished, and the dominant term of the second order, viz., $\eta\epsilon^2$.
The co-efficient of this term is shown to be strictly positive by analyzing the tail probabilities.
Therefore, we have
\begin{align*}
    \E\left[\inner{\BFone_{n+m}}{\BFq}\right] \sim \frac{1}{\epsilon}, \quad
   \lr\sim\eta\epsilon^2.
\end{align*}
To achieve the optimal trade-off between expected queue length and profit loss, we choose $\epsilon \sim \eta^{-1/3}$, which results in the $\eta^{1/3}$ profit loss in Theorem~\ref{theorem: lowerboundmultiplelink}. 
\end{remark}
We can further relax Condition \ref{ass: generalpricinggraph}\ref{condition: aplusb} in the special case of a single-link system ($m=n=1$) operating under any two-price policy. The result is stated below and the proof can be found in Appendix~\ref{appendix: twopricesinglelinktwosidedqueue}.
\begin{proposition} \label{prop: twopricesinglelinktwosidedqueue}
For a family of single-link two-sided queue parametrized by $\eta$, any two pricing policy given by Eq~\eqref{eq:two-price_policy-cust} with $\sigma^{\eta}=\eta^{\beta}$ for some $\beta<1$ and $\tau_{\max}^{\eta}=\eta^{\alpha}$ for some $\alpha \in \mathbb{R}_+$, will have a profit loss $\lr$ at least  $\Omega({\eta}^{1/3})$.
The choice of $\tau_{\max}^{\eta}=\eta^{1/3}$ and $\sigma^{\eta}=\eta^{2/3}$ and any positive constants $\BFtheta$ and $\BFphi$ provides the optimal profit loss $\Theta({\eta}^{1/3})$.
\end{proposition} 
Before concluding this section, note that without any further assumptions, the dependence of profit loss $\lr$ on the number of customer types $n$ and the number of server types $m$ can be linear in the worst case. To see this, if the system is consisted of $n$
independent single-link two-sided queues, i.e.,
the bipartite graph is such that $|M|=|N|=n$ and $E=\{(i,i): i \in M\}$, then by Theorem~\ref{theorem: lowerboundmultiplelink},
the profit loss $\lr$ is trivially lower bounded by $nK\eta^{1/3}$. In the next section, we impose additional conditions on the bipartite graph and show tight lower and upper bound on the profit-loss in terms of the system size.
\section{Asymptotically Optimal Policies in the Large-System Regime: the Superiority of Max-Weight Matching}
\label{sec:large-n}
In this section, we present further insights into the max-weight matching algorithm in a \emph{large-system} regime, in which both the arrival rates and the number of customer and server types increase. 
First, we will show that max-weight is delay optimal. Under the fluid pricing policy, max-weight matching minimizes the probability of hitting the queue length threshold among all matching policies. Under the two-price policy, max-weight matching minimizes the expected sum of queue lengths among all possible matching policies.
Second, we compare max-weight matching with a randomized matching policy with probabilities specified by the fluid model and show that max-weight has smaller loss in terms of the number of customer/server types.
Third, we prove that the profit loss of max-weight matching achieves a tight lower bound in the large-system regime.
Together, these results show the superiority of the max-weight policy.
We start by establishing the \emph{state space collapse} under max-weight matching.
State space collapse means that all the customer queues are almost equal in length and all the server queues are almost equal in length; hence, with high probability, only customers or only servers are waiting in the system. This implies that max-weight ends up matching the maximum possible number of customer-server pairs, as only the excess customers/servers are waiting in the system. To achieve the state space collapse, we propose a complete resource pooling condition on the compatibility graph. Similar conditions have been proposed for single-sided queues \cite[Assumption 2.4]{CRP_queueing} \cite[Definition 1]{CRP_manufacturing}.
\begin{condition}[Complete Resource Pooling (CRP)] \label{condition: crp}
There exists an optimal solution $(\BFlambda^*, \BFmu^*)$ to the fluid problem Eq~\eqref{eq:fluid_opt} such that
for all non-empty $J \subsetneq {M}$ and for all non-empty $I \subsetneq {N}$, it holds that
\begin{align*}
    \sum_{j \in J} \lambda_j^* < \sum_{i : \exists j \in J, (i,j) \in E} \mu_i^*, \quad \sum_{i \in I} \mu_i^* < \sum_{j: \exists i \in I, (i,j) \in E} \lambda_j^*.
\end{align*}
\end{condition}
It is straightforward to verify that the CRP condition implies the connectedness of the graph $G(N \cup M,E)$. The CRP condition also implies that the optimal solution of the fluid problem is in the interior of the feasible region. The following lemma formalizes this observation.
\begin{lemma} \label{lemma: positive_chi}
If Condition~\ref{condition: crp} is satisfied,  there exists $\BFchi^* \geq \BFzero$ such that $\chi^*_{ij}>0$ for all $(i,j) \in E$ and $(\BFlambda^*,\BFmu^*,\BFchi^*)$ is an optimal solution to the fluid problem Eq~\eqref{eq:fluid_opt}.
\end{lemma}
(All the proofs in this section can be found in Appendix~\ref{appendix: further-analysis-max-weight}.)
The above result is not surprising as it is known in the heavy traffic literature \citep{atilla_srikant,lange2019heavy} that if the arrival rate is approaching a point on the boundary of the capacity region in the interior of a facet, then the system exhibits complete resource pooling. The above lemma implies that $E_r = \emptyset$, i.e., no edges are redundant as defined in Definition~\ref{def:redundant_edges}. Thus, the max-weight policy and the modified max-weight policy are identical, and so, we do not make this distinction for the rest of the paper.
However, the analysis of state space collapse for two-sided queues does not follow immediately from the literature of single-sided queues and is more involved. We propose a Lyapunov function approach and use the drift method to show state space collapse. To simplify the analysis, in this section we restrict to a setting where $m=n$ and there exists a perfect matching in the graph $G(N \cup M,E)$.  
\begin{condition} \label{condition: perfect_matching}
The graph $G(N \cup M, E)$ has a perfect matching. Without loss of generality, we assume that server type $i$ is connected to customer type $i$ for all $i\in[n]$.
\end{condition}
In general, if $m\neq n$, and if the above condition is not satisfied, we show in Appendix~\ref{subsec: relax-perfect-matching} that the pricing and matching problem under a given general graph can be reformulated as a problem under a new graph where the above condition is satisfied. Thus, the results in the following propositions and the theorem can be applied (with minor modifications as shown in Appendix~\ref{subsec: relax-perfect-matching}) even when the above condition does not hold. 
\subsection{Delay Optimality of Max-Weight Matching}
\label{sec:delay-optimal}
We first show the delay optimality of the max-weight matching algorithm among all possible matching algorithms under the fluid pricing policy. 
\begin{proposition} \label{prop: max-weight-optimal-fluid-SSC}
Under the fluid pricing policy with any matching algorithm, we have
\begin{align*}
    q_{\max}^\eta\left(\sum_{j=1}^n \lambda_j^* \Pr\left[q_j^{(c)}=q_{\max}^\eta\right]+\sum_{i=1}^n \mu_i^*\Pr\left[q_i^{(s)}=q_{\max}^\eta\right]\right) \geq \frac{\inner{\BFone_n}{\BFlambda^*}+\inner{\BFone_n}{\BFmu^*}}{2n+1/q_{\max}^\eta}.
\end{align*}
Furthermore, under the fluid pricing policy with the max-weight matching algorithm, if $q_{\max}^\eta\to\infty$ as $\eta\to\infty$, we have
\begin{align*}
    \lim_{\eta \to \infty} q_{\max}^\eta\left( \sum_{j=1}^n \lambda_j^* \Pr\left[q_j^{(c)}=q_{\max}^\eta\right]+\sum_{i=1}^n \mu_i^*\Pr\left[q_i^{(s)}=q_{\max}^\eta\right]\right) = \frac{\inner{\BFone_n}{\BFlambda^*}+\inner{\BFone_n}{\BFmu^*}}{2n}.
\end{align*}
\end{proposition}
The above proposition states that the max-weight algorithm (asymptotically) minimizes the proportion of time spent in the threshold state among all possible matching algorithms, hence minimizing the revenue loss caused by hitting the queue length thresholds. 
Similarly, the max-weight matching algorithm is delay optimal under the two-price policy. The following proposition states that the max-weight algorithm (asymptotically) minimizes the expected total queue length under the two-price algorithm among all possible matching algorithms. 
\begin{proposition} \label{prop: max-weight-optimal-two-price-SSC}
Under the two pricing policy with $\BFtheta=\BFphi=\BFone_n$ and any matching policy, the expected total queue length satisfies
\begin{align*}
    \frac{\sigma^\eta}{\eta}\E[\inner{\BFone_{2n}}{\BFq}] \geq \frac{\inner{\BFone_n}{\BFlambda^*}+\inner{\BFone_n}{\BFmu^*}}{2n} - \frac{\sigma^\eta}{\eta}.
\end{align*}
Furthermore, under the two-price policy with $\BFtheta=\BFphi=\BFone_n$ and the max-weight matching policy, if 
$\lim_{\eta\to\infty}\sigma^\eta/\eta=0$
and $\lim_{\eta\to\infty}\sigma^\eta\tau_{\max}^\eta/\eta=0$, we have
\begin{align*}
    \lim_{\eta\to\infty}\frac{\sigma^\eta}{\eta}\E[\inner{\BFone_{2n}}{\BFq}] =\frac{\inner{\BFone_n}{\BFlambda^*}+\inner{\BFone_n}{\BFmu^*}}{2n}.
\end{align*}
\end{proposition}
Notice that the queue length bound in Proposition~\ref{prop: max-weight-optimal-two-price-SSC} is tighter than the bound in Lemma~\ref{lemma: posrec}, because former requires the CRP condition (Condition~\ref{condition: crp}) whereas the latter does not require such condition.
Together, Propositions \ref{prop: max-weight-optimal-fluid-SSC} and \ref{prop: max-weight-optimal-two-price-SSC} establish the asymptotic delay optimality of the max-weight algorithm.
\subsection{Max-Weight versus Randomized Matching}
\label{sec:MW vs Rand}
In this section, 
we compare the max-weight policy with a
randomized matching policy (defined in Algorithm \ref{alg: random_non_empty}) resulting from the fluid model. 
The randomized matching algorithm matches an incoming arrival
to compatible types at fixed probabilities, which are determined by the fluid solution $\BFchi^*$ (see Eq~\eqref{eq:fluid_opt}). 
If some queues are empty, the probabilities are rescaled proportionally to match only nonempty queues.
Unlike the max-weight algorithm, the randomized matching algorithm does not use information about the queue lengths (except for the emptiness of the queues). 
\begin{algorithm}[!htbp]
\caption{Randomized Matching (Nonempty Queues First)} \label{alg: random_non_empty}
{\fontsize{10}{11}\selectfont
\begin{algorithmic}
\STATE \textbf{input:} new arrival $\BFa(k)$, queue length $\BFq(k)$, the fluid solution $\BFchi^*$ \COMMENT{$k$ is a decision epoch}
\STATE \textbf{initialization:} $\BFy(k)=\BFzero$
    \FOR{$i \in {N}$}
    \IF {$a_i^{(s)}(k)=1$}
        \STATE set $y_{ij}^{}(k)=1$ with probability $\frac{\chi^*_{ij}\mathbbm{1}_{\{q_j^{(c)}>0\}}}{\sum_{j'=1}^m \chi_{ij'}^*\mathbbm{1}_{\{q_{j'}^{(c)}>0\}}}$ for all $j \in {M}$.
    \ENDIF
    \ENDFOR
    \FOR{$j \in {M}$}
    \IF {$a_j^{(c)}(k)=1$ }
        \STATE set $y_{ij}^{}(k)=1$ with probability $\frac{\chi^*_{ij}\mathbbm{1}_{\{q_i^{(s)}>0\}}}{\sum_{i'=1}^n \chi_{i'j}^*\mathbbm{1}_{\{q_{i'}^{(s)}>0\}}}$ for all $i \in {N}$.
    \ENDIF
    \ENDFOR
\STATE \textbf{output:} matching decision $\BFy(k)$
\end{algorithmic}
}
\end{algorithm}
We analyze the profit losses of these two matching algorithms and its dependence
on the number of customer/server types $n$ when $\eta \rightarrow \infty$.
First, we consider the fluid pricing policy.
The theorem below shows that even though both max-weight and randomized matching have $O(\eta^{1/2})$ profit loss,
max-weight matching is order $n^{1/2}$ better than randomized matching policy.
\begin{theorem} \label{theo: fluid_max_weight_SSC}
Suppose a family of two-sided queues is given by the bipartite graph $G(N \cup M, E)$ parametrized by $\eta$. Under the fluid price policy Eq~\eqref{eq:fluid_pricing} and randomized matching policy (Algorithm \ref{alg: random_non_empty}), for $q_{\max}^\eta=\gamma\eta^{1/2}$, we have $\lr=O(\eta^{1/2})$. For any $\gamma>0$, there exists $(\BFlambda^*,\BFmu^*,\BFchi^*)$ satisfying Condition \ref{condition: crp} and \ref{condition: perfect_matching} such that
\begin{align*}
   \liminf_{\eta \rightarrow \infty} \frac{\lr}{\eta^{1/2}} = \Omega(n).
\end{align*}
In addition, under the fluid price policy \eqref{eq:fluid_pricing} and max-weight matching \eqref{eq:max-weight} for $q_{\max}^{\eta}=\sqrt{\eta/n}$, we have
\begin{align*}
   \limsup_{\eta \rightarrow \infty} \frac{\lr}{\eta^{1/2}}\leq n^{1/2}\left(\frac{\inner{\BFone_n}{\BFlambda^*}+\inner{\BFone_n}{\BFmu^*}}{2n} \max_{j \in {N}} F_j(\lambda_j^*)+2\max_{i \in N, j\in M}\{s_i^{(s)},s_j^{(c)}\}\right)=O(n^{1/2}).
\end{align*}
\end{theorem}
Next, we compare the max-weight and randomized matching algorithms for the two-price pricing policy. The theorem below shows that
both algorithms achieve $O(\eta^{1/3})$ profit loss, whereas max-weight is order $n^{2/3}$ better than randomized matching.
\begin{theorem} \label{theo: two_price_maxweight_SSC}
Suppose a family of two-sided queues is given by the bipartite graph $G(N \cup M, E)$ parametrized by $\eta$. Under the two-price policy  Eq~\eqref{eq:two-price_policy-cust} and randomized matching policy (Algorithm \ref{alg: random_non_empty}), for $\sigma^\eta=\eta^{2/3}$ and $\tau_{\max}^\eta=\gamma\eta^{1/3}$, we have $\lr=O(\eta^{1/3})$. For any choice of $\BFtheta>\BFzero$, $\BFphi>\BFzero$ and $\gamma>0$, there exists $(\BFlambda^*,\BFmu^*,\BFchi^*)$ satisfying Condition \ref{condition: crp} and \ref{condition: perfect_matching} such that
\begin{align*}
    \liminf_{\eta \rightarrow \infty} \frac{\lr}{\eta^{1/3}}=\Omega(n).
\end{align*}
In addition, under the two-price policy Eq~\eqref{eq:two-price_policy-cust} and max-weight matching Eq~\eqref{eq:max-weight} with $\BFtheta=\BFone_n$ and $\BFphi=\BFone_n$, $\sigma^\eta=n^{-1/3}\eta^{2/3}$, if $\lim_{\eta\to\infty}\tau_{\max}^\eta/\eta^{1/3}=0$, we have
\begin{align*}
\limsup_{\eta \rightarrow \infty} \frac{\lr}{\eta^{1/3}}\leq{}& \left(\sum_{i \in {N}} \left(\frac{\mu_i^*G_i''(\mu_i^*)}{2}+G_i'(\mu_i^*)\right)-\sum_{j \in {N}} \left(\frac{\lambda_j^*F_j''(\lambda_j^*)}{2}+F_j'(\lambda_j^*)\right)\right. \\
    &\left.+\frac{\max_{i \in N,j\in M}\{s_i^{(s)},s_j^{(c)}\}}{2}\left(\inner{\BFone_n}{\BFlambda^*}+\inner{\BFone_n}{\BFmu^*}\right)\right)n^{-2/3}=O(n^{1/3}).
\end{align*}
\end{theorem}
\subsection{Lower Bound}
\label{sec:lower-bound-n}
In this section, we prove a lower bound on the profit-loss for the large-system regime. We show that under mild assumptions, any pricing and matching policy has a profit loss of $\Omega((n\eta)^{1/3})$. In light of this lower bound result and the upper bound from Theorem~\ref{theo: two_price_maxweight_SSC}, the two-price max-weight policy achieves the optimal rate for the large-system regime.
To prove the lower bound, we consider the following problem instance. We assume the number of customer and server types are equal (i.e. $n=m$) and the bipartite matching graph is a complete graph.
Intuitively, the complete graph presents a best case scenario, as it provides with the maximum flexibility to match any  customer-server pairs.  To define the supply and demand rate, we assume that there exist functions $F,G : \mathbb{R}_+ \rightarrow \mathbb{R}_+$ such that $F_j=F$ for all $j \in M$ and $G_i=G$ for all $i \in N$. It is easily verified that this problem instance satisfies the complete resource pooling condition (Condition~\ref{condition: crp}).
We also assume that the unit waiting cost is $\BFs=\BFone_{2n}$. By symmetry, we can conclude that the fluid solution $(\BFlambda^\star,\BFmu^\star)=(\lambda^\star\BFone_{n},\mu^\star\BFone_{n})$ for some $\lambda^\star=\mu^\star>0$. 
We show the following properties of the optimal pricing and matching policy for this instance.
\begin{proposition} \label{prop: optimal_policy_complete_graph}
Consider a two-sided queueing system defined by a complete bipartite graph $G(M \cup N, M \times N)$. Assume that $n=m$, $F_j=F$ for all $j \in M$, $G_i=G$ for all $i \in N$, and $\BFs=\BFone_{2n}$. Then, there exists an optimal policy $(\BFlambda^\star(\cdot),\BFmu^\star(\cdot),\BFx^\star(\cdot))$ for the DTMDP \eqref{MDP} that satisfies the following:
\begin{enumerate}[label=(\alph*)]
    \item \label{prop: equal_arrival_rate_queues} $\lambda_{j_1}^\star(\BFq)=\lambda_{j_2}^\star(\BFq)$, and $\mu_{i_1}^\star(\BFq)=\mu_{i_2}^\star(\BFq)$ for all $j_1,j_2 \in M$, $i_1,i_2 \in N$, and $\BFq \in \mathbb{Z}_+^{2n}$.
    \item \label{prop: equal_imbalance} $\left(\BFlambda^\star(\BFq_1),\BFmu^\star(\BFq_1)\right)=\left(\BFlambda^\star(\BFq_2),\BFmu^\star(\BFq_2)\right)$ if $\inner{\BFone_n}{\BFq_1^{(c)}}=\inner{\BFone_n}{\BFq_2^{(c)}}$ and $\inner{\BFone_n}{\BFq_1^{(s)}}=\inner{\BFone_n}{\BFq_2^{(s)}}$.
    \item \label{prop: non_idling_matching} $\inner{\BFone_{2n}}{\BFx^\star(\BFq)}=2\min\left\{\inner{\BFone_n}{\BFq^{(c)}},\inner{\BFone_n}{\BFq^{(s)}}\right\}$ for all $\BFq \in \mathbb{Z}_+^{2n}$.
\end{enumerate}
\end{proposition}
Part~\ref{prop: equal_arrival_rate_queues} implies that the optimal prices are equal for all the customer queues and server queues, respectively. Part~\ref{prop: equal_imbalance} implies that the arrival rates depends on the state of the system $\BFq$ only through the total number of customers and total number of servers in the system. Both of these conclusions are intuitive as the problem instance is defined symmetrically for all types. Lastly,  \ref{prop: non_idling_matching} implies that the optimal policy will always match the maximum possible number of customer-server pairs. In particular, there is no incentive to hold a compatible customer-server pair due to the complete graph structure. 
Motivated by the above proposition, we restrict ourselves to a family of \emph{symmetric pricing policies} defined as follows:
\begin{align}
    \lambda_j^\eta(\BFq)&=\eta \lambda^\star + f\left(\frac{\inner{\BFone_n}{\BFq^{(c)}}-\inner{\BFone_n}{\BFq^{(s)}}}{\eta^{\alpha} }\right)\eta^{\beta}  \quad \forall j \in M, \\
    \mu_i^\eta(\BFq)&=\eta \mu^\star + g\left(\frac{\inner{\BFone_n}{\BFq^{(c)}}-\inner{\BFone_n}{\BFq^{(s)}}}{\eta^{\alpha} }\right)\eta^{\beta} \quad \forall i \in N.
\end{align}
The above policy family is based on the one considered in Section~\ref{sec:lower_bound}. However, we make two changes on Eqs \eqref{eq: pricinggraphgeneralcust} and \eqref{eq: pricinggraphgeneralser}: we assume $f_j=f$ ($\forall j \in M$) and $g_i=g$ ($\forall i \in N$), since the prices are symmetric for all customer and server types
by part~\ref{prop: equal_arrival_rate_queues}; the policy depends on the state $\BFq$ through the difference between total customers and total servers, according to \ref{prop: equal_imbalance} and \ref{prop: non_idling_matching}.
We impose some technical conditions on $f$ and $g$.
\begin{condition} \label{cond: n_dep}
\begin{enumerate}[label=(\alph*)]
    \item There exists a constant $\Gamma$ which may depend on $n$ but not on $\eta$, such that $|f(z)| \leq \Gamma$ and $|g(z)| \leq \Gamma$ for all $z \in \mathbb{R}$. \label{cond: bounded_n}
    \item $0<\alpha+\beta < 1$. \label{cond: scaling_n}
    \item There exist constants $\kappa,\delta>0$, such that if $z>\kappa$, then either $f(z)<-\delta \Gamma$ or $g(z)>\delta \Gamma$. In addition, if $z<-\kappa$, then either $f(z)>\delta \Gamma$ or $g(z)<-\delta \Gamma$. \label{cond: stability_n}
\end{enumerate}
\end{condition}
Condition \ref{cond: n_dep}\ref{cond: bounded_n} is analogous to Condition \ref{ass: generalpricinggraph}\ref{condition: boundedfg} implying that $f$ and $g$ are bounded. 
Condition \ref{cond: n_dep}\ref{cond: scaling_n} is similar to Condition \ref{ass: generalpricinggraph}\ref{condition: aplusb}. Lastly, Condition \ref{cond: n_dep}\ref{cond: stability_n} is adapted from Condition \ref{ass: generalpricinggraph}\ref{condition: posrecsystem} using the properties from Proposition~\ref{prop: optimal_policy_complete_graph}.
Now we present the result on the lower bound below.
\begin{theorem} \label{theo: lower_bound_n}
Consider a family of two-sided queues parametrized by $\eta$ and $n$ in the large-system regime. Under any symmetric pricing policy satisfying Condition \ref{cond: n_dep} and any matching policy, we have 
\begin{align*}
    \liminf_{\eta \rightarrow \infty} \frac{\lr}{\eta^{1/3}}=\Omega(n^{1/3}).
\end{align*}
\end{theorem}
\begin{remark}[Intuitive Explanation of $n^{1/3}$] 
We prove the lower bound by considering the complete graph instance with symmetric demand and supply functions. By Proposition~\ref{prop: optimal_policy_complete_graph}, the optimal policy will match the maximum possible number of customer-server pairs, so only the excess customers/servers are waiting in the system. Thus, the system mimics a single link two-sided queue ($m=n=1$) with customer and server arrival rates given by $\inner{\boldsymbol{1}}{\boldsymbol{\lambda}}$ and $\inner{\boldsymbol{1}}{\boldsymbol{\mu}}$ respectively. Under the large-system regime where the number of types is scaled by $n$ and the arrival rate of each type is scaled by $\eta$, the total arrival rates for customers and servers are scaled by $n \eta$.
The profit loss of any pricing and matching policy is lower bounded by the cube root of the total arrival rate (see Remark~\ref{rmk:eta-lb}), so the profit loss in the large-system regime is of order
$\Omega((n\eta)^{1/3})$. Meanwhile, since the max-weight matching policy leads to state space collapse under the CRP condition, the two-sided queueing system behaves like a single link system where all customer and server types are pooled together, so the two-price max-weight policy able to achieve a tight $O((n\eta)^{1/3})$ profit loss rate.
\end{remark}
\section{Numerical Experiments}
\label{sec:simulation}
\subsection{Single-Link Systems}
\begin{figure}[!tbp]
  \centering
  \begin{minipage}[b]{0.48\textwidth}
  \FIGURE
    {\includegraphics[width=\textwidth]{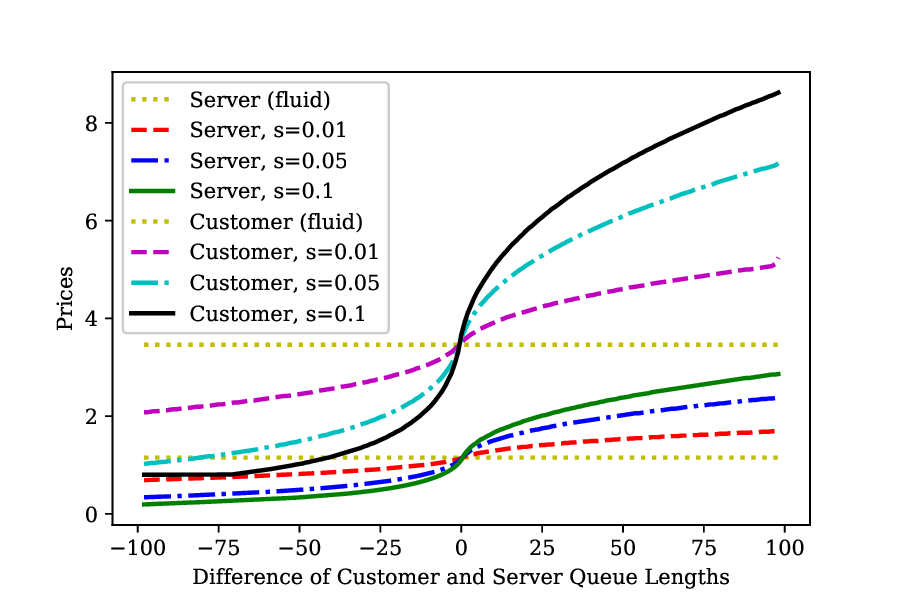}}
    {Optimal pricing policies under different values of penalty coefficients.
    \label{fig: pricing-policies}}
    {}
  \end{minipage}
  \hfill
\begin{minipage}[b]{0.48\textwidth}
  \FIGURE
    {\includegraphics[width=\textwidth]{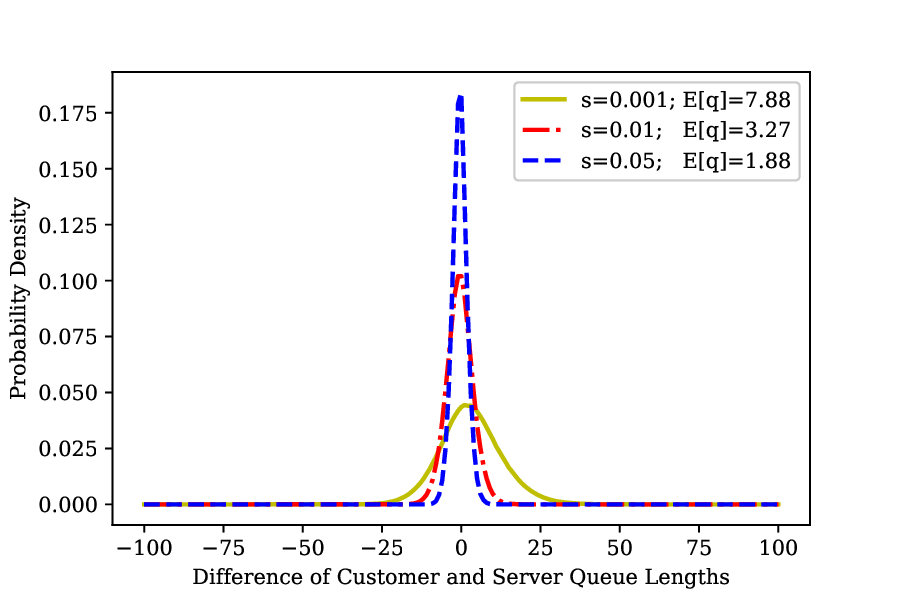}}
    {Stationary distribution of queue length under different penalty coefficients.
    \label{fig: stationary-distribution}}
    {}
  \end{minipage}
\end{figure}

Our first experiment analyzes a single-link system with one server type and one customer type. 
In this case, the system state of the MDP is represented by a single variable, namely, the difference between the customer queue length and the server queue length (a detailed discussion of this system is included in Appendix~\ref{sec:single-link}). We will solve the optimal policy of the MDP and compare it with the fluid pricing policy and two-price policy. 
We assume a supply curve given by $p_1=\lambda^{0.5}$ and a demand curve given by $p_2=4\mu^{-0.5}$.
With these supply and demand curves, the optimal profit of the fluid model is $3.08$ when $\lambda=\mu=4/3$,  $p_1=1.15$ and $p_2=3.46$. 
We then calculate the optimal pricing policy of the long-run average cost MDP using relative value iteration. 
Figure~\ref{fig: pricing-policies} shows the optimal pricing policy under three different values of the penalty coefficient ($s$), as well as the optimal price of the fluid model.
The result shows that the optimal customer price is always above the server price, and both prices are increasing with the queue length difference.
Intuitively, if the system has more customers, 
the customer price should be increased to reduce the customer arrival rate and server price should be increased to increase the server arrival rate.
This observation verifies Proposition~\ref{theo: monotonicity} in Appendix~\ref{sec:single-link}.
As $s$ increases, more weight is given to the waiting cost (or equivalently, customers and servers become more sensitive to delays), so the price increases more steeply as the number of customers and servers waiting in the system increases. 
Figure~\ref{fig: stationary-distribution} show the stationary distribution and the mean of queue length for different values of the penalty coefficient ($s$).
As expected, when $s$ increases, the queue length is more concentrated around 0.
\begin{figure}[!htbp]
  \FIGURE
    {\includegraphics[width=0.5\textwidth]{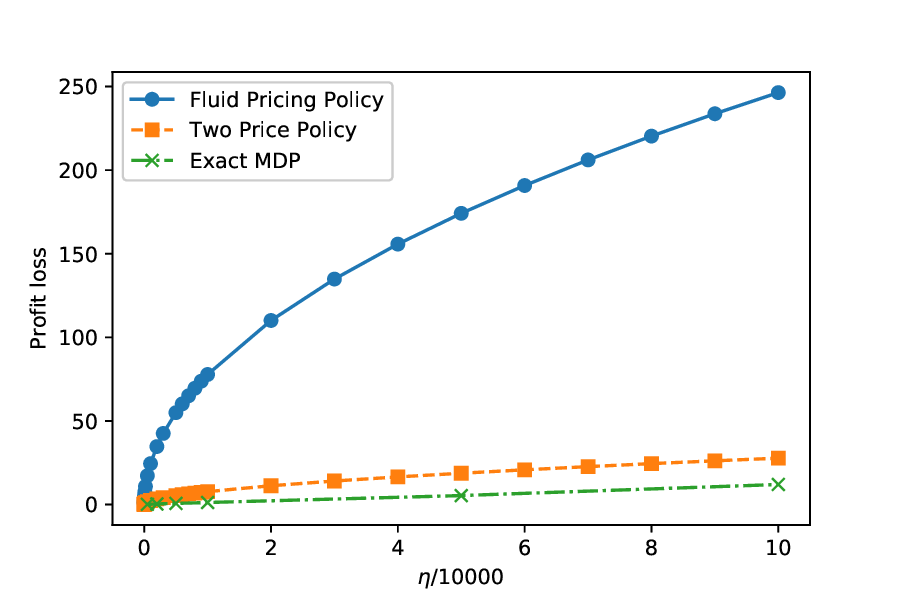}}
    {Performance of two-price and fluid pricing policy compared to the exact solution.
    \label{fig: tpandfluid}}
    {}
\end{figure}
Furthermore, we simulate the profit loss under the fluid pricing policy and two-price policy and compare it with the theoretical result presented before and also with the exact solution obtained by solving the MDP. The result is presented in Figure~\ref{fig: tpandfluid}. The profit loss under the fluid pricing policy has an order of $\sqrt{\eta}$ and that under the two-price policy has an order of $\eta^{1/3}$, verifying Theorem \ref{theo: fluidpricingpolicy} and Theorem \ref{theorem: twoprice}. Also observe that the profit loss under the two-price policy is not much different from that of the optimal profit loss, demonstrating the effectiveness of a two-price policy.
\subsection{Max-Weight v/s Modified Max-Weight} \label{sec: sim_erratum}
We set $n=m=2$ with the compatibility graph given by $E = \{(1, 1), (1, 2), (2, 2)\}$. We assume the unit holding cost is $s=1$. The parameters of the two-price policy are chosen to be $\tau_{\max}^\eta = 0$ and $\sigma^\eta = \eta^{2/3} n^{-1/3}$. The demand curves are given by $F_1(\lambda_1) = 5 - \lambda_1$ and $F_2(\lambda_2) = 4 - \lambda_2$ and the supply curves are given by $G_1(\mu_1) = 1.5 \mu_1$ and $G_2(\mu_2) = \mu_2$. The optimal solution of the fluid problem \eqref{eq:fluid_opt} for these choices of demand and supply curves is $\BFlambda^* = \BFmu^* = \BFone_2$ and $\chi^*_{ij} = 1$ if $i = j$ and $0$ otherwise. As $(\BFlambda^*, \BFmu^*, \BFchi^*)$ is the unique optimal solution of \eqref{eq:fluid_opt} and $\chi_{12}^* = 0$, we conclude that $E_r = \{(1, 2)\}$. As $E_r \neq \emptyset$, the CRP condition is not satisfied by \cite[Lemma 3]{varma2023dynamic}.

For the above-defined parameters, we implement the max-weight matching policy \cite[Algorithm 1]{varma2023dynamic} and the modified max-weight matching policy (Algorithm~\ref{alg: maxweight}) combined with the two-price policy. The max-weight matching policy makes use of all edges in $E$ and we break the ties in favor of the edges $(1, 1)$ and $(2, 2)$ over the edge $(1, 2)$. The modified max-weight matching policy only uses the edges $\{(1, 1), (2, 2)\}$ as the edge $(1, 2)$ is redundant.

\begin{figure}[!tbh]
  \centering
\FIGURE
{
\begin{tikzpicture}
\node at (0, 0) {\includegraphics[scale=0.55]{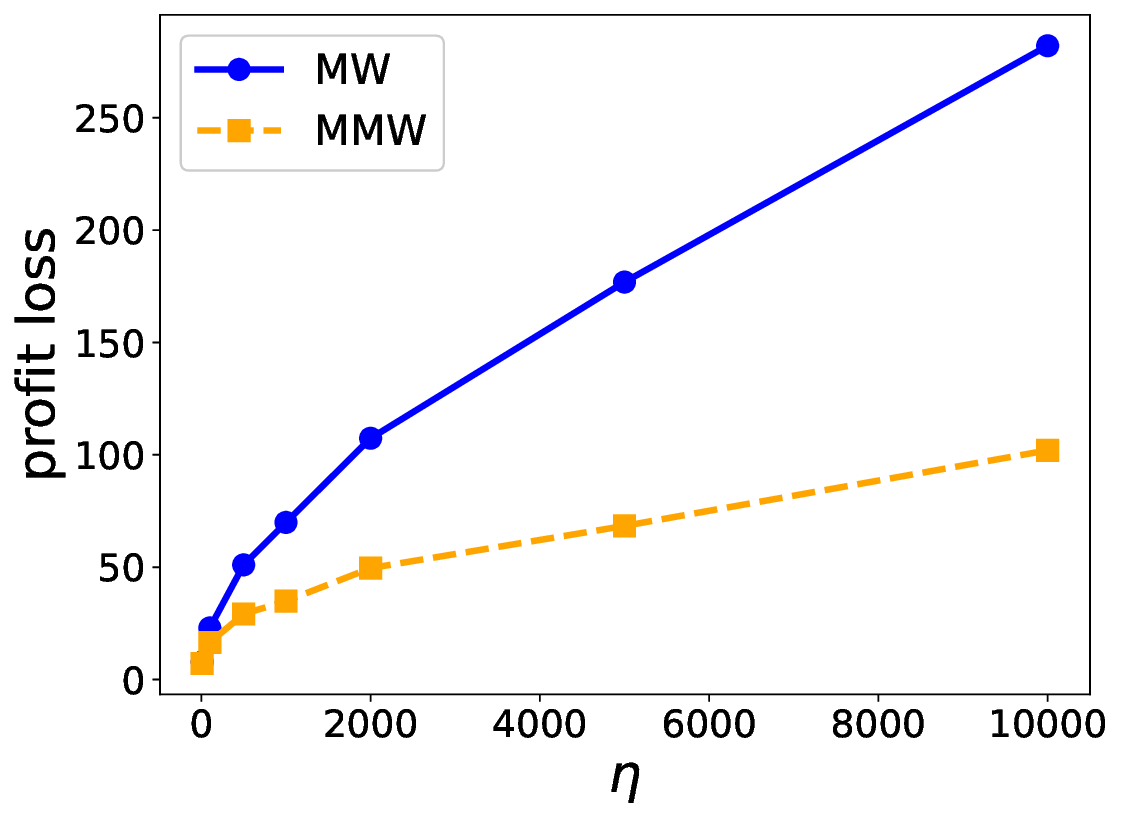}};
\node at (8, 0) {\includegraphics[scale=0.55]{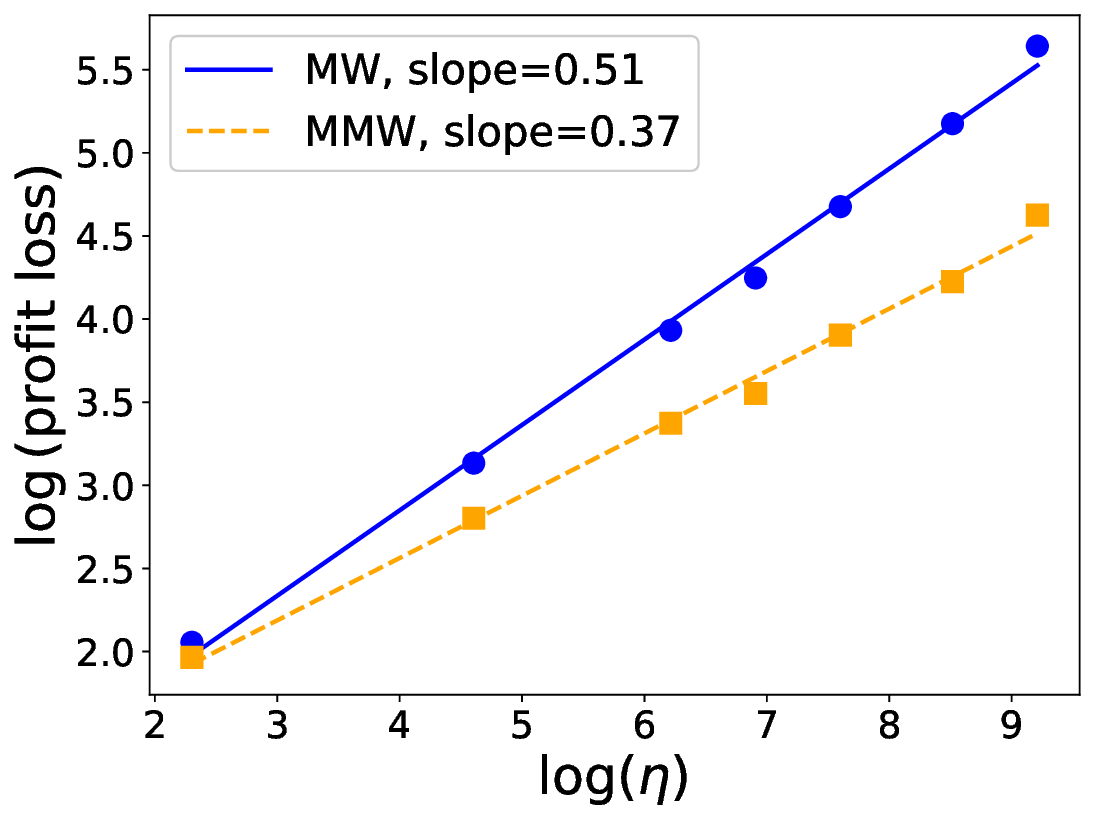}};
\end{tikzpicture}
}{Profit-loss for the max-weight and modified max-weight matching policy combined with two-price policy. For reproducibility, the code is open source and can be found here \cite{varma_erratum_code}.
    \label{fig: mw_vs_mmw}}
    {}
\end{figure}

We report the profit loss for $\eta \in \{10, 100, 500, 1000, 2000, 5000, 10000\}$ in Figure~\ref{fig: mw_vs_mmw}. As seen in the left frame, the modified max-weight matching policy has a smaller profit-loss compared to the max-weight matching policy. In particular, as seen in the right frame, the profit-loss scales approximately as $\Theta(\sqrt{\eta})$ and $\Theta(\eta^{1/3})$ for the max-weight and modified max-weight matching policy respectively. Thus, Theorem~\ref{theorem: twoprice} may not hold under the max-weight matching policy but continues to hold for the modified max-weight policy.

\begin{figure}[!tbp]
  \centering
  \begin{minipage}[b]{0.49\textwidth}
  \FIGURE
    {\includegraphics[width=\textwidth]{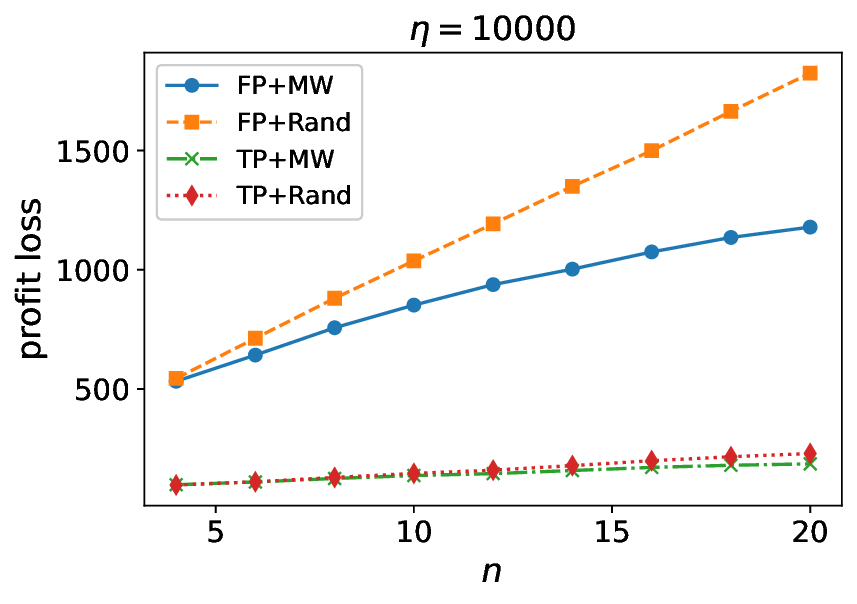}}
    {{Profit losses under FP+MW, FP+Rand, TP+MW and TP+Rand for different $n$ when $\eta=10000$.}
    \label{fig: profitloss eta10000 mnnotequal}}
    {}
  \end{minipage}
  \begin{minipage}[b]{0.49\textwidth}
  \FIGURE
    {\includegraphics[width=\textwidth]{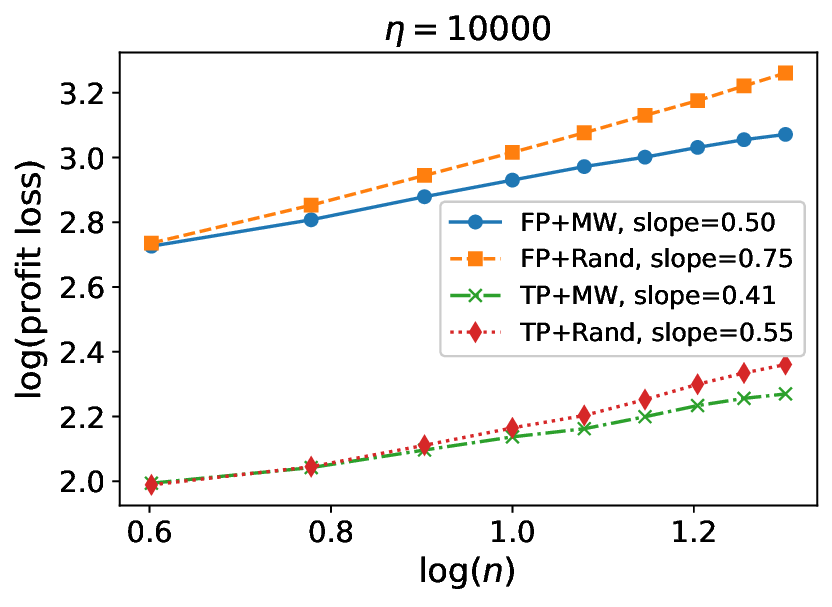}}
    {Log-log plot of profit losses under FP+MW, FP+Rand, TP+MW and TP+Rand for different $n$ when $\eta=10000$.
    \label{fig: logprofitloss eta10000 mnnotequal}}
    {}
  \end{minipage}
\end{figure}
\vspace{-30pt}
\subsection{Systems with Multiple Types}
\begin{figure}[!tbp]
  \centering
  \begin{minipage}[b]{0.49\textwidth}
  \FIGURE
    {\includegraphics[width=\textwidth]{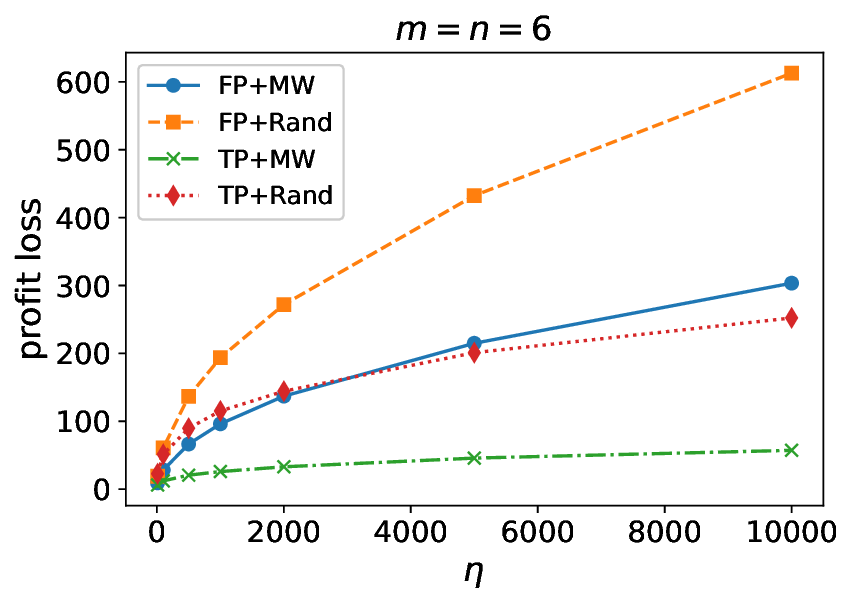}}
    {{Profit losses under FP+MW, FP+Rand, TP+MW and TP+Rand for different $\eta$ when $m=n=6$}.
    \label{fig: profitloss mnequal}}
    {}
  \end{minipage}
  \begin{minipage}[b]{0.49\textwidth}
  \FIGURE
    {\includegraphics[width=\textwidth]{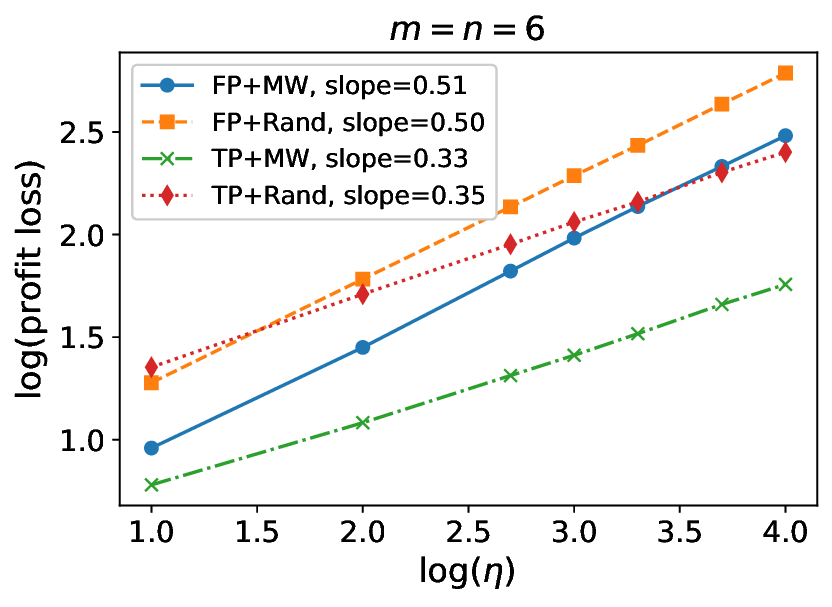}}
    {Log-log plot of profit losses under FP+MW, FP+Rand, TP+MW and TP+Rand for different $\eta$ when $m=n=6$.
    \label{fig: logprofitloss mnequal}}
    {}
  \end{minipage}
\end{figure}

\begin{figure}[!tbp]
  \centering
  \begin{minipage}[b]{0.49\textwidth}
  \FIGURE
    {\includegraphics[width=\textwidth]{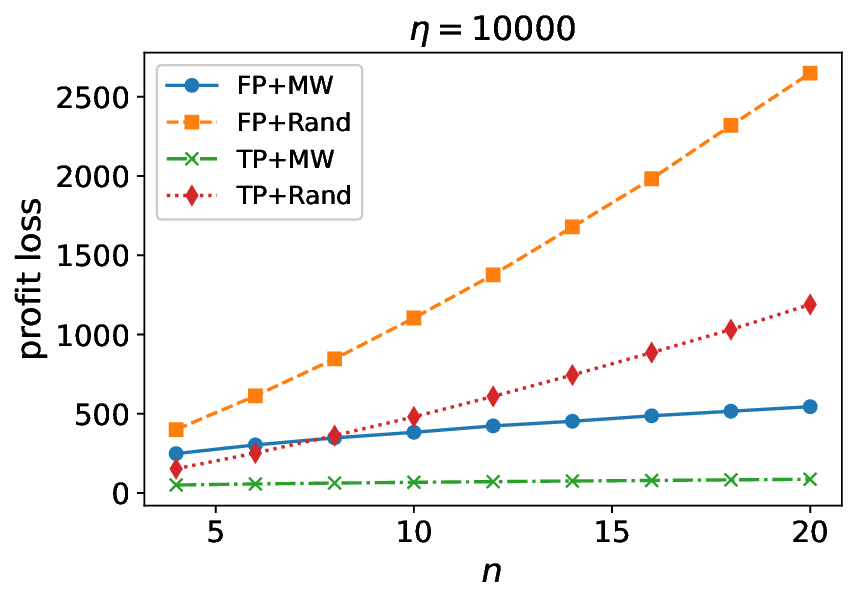}}
    {{Profit losses under FP+MW, FP+Rand, TP+MW and TP+Rand for different $n$ when $\eta=10000$.}
    \label{fig: profitloss eta10000 mnequal}}
    {}
  \end{minipage}
  \begin{minipage}[b]{0.49\textwidth}
  \FIGURE
    {\includegraphics[width=\textwidth]{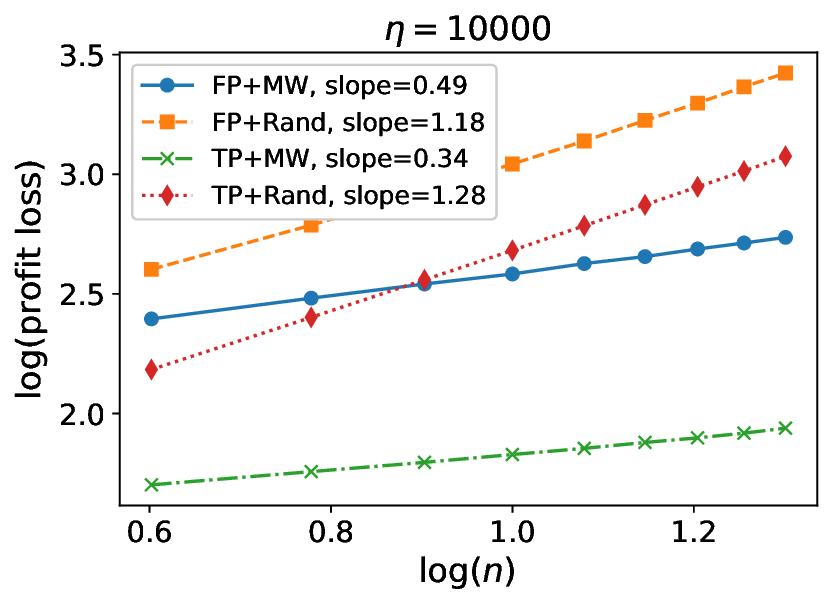}}
    {Log-log plot of profit losses under FP+MW, FP+Rand, TP+MW and TP+Rand for different $n$ when $\eta=10000$.
    \label{fig: logprofitloss eta10000 mnequal}}
    {}
  \end{minipage}
\end{figure}

\begin{figure}[!tbp]
  \centering
  \begin{minipage}[b]{0.49\textwidth}
  \FIGURE
    {\includegraphics[width=\textwidth]{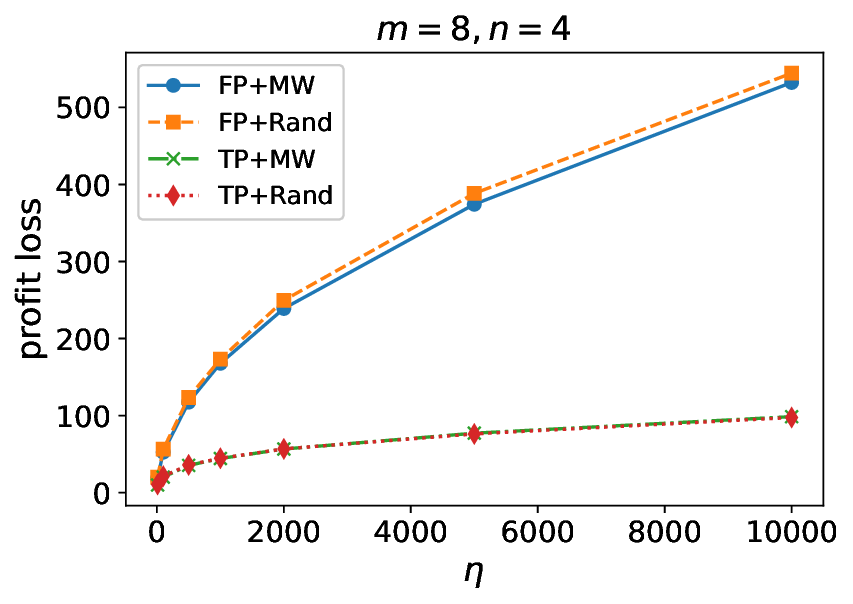}}
    {{Profit losses under FP+MW, FP+Rand, TP+MW and TP+Rand for different $\eta$ when $m=4$ and $n=2$.}
    \label{fig: profitloss mnnotequal}}
    {}
  \end{minipage}
  \begin{minipage}[b]{0.49\textwidth}
  \FIGURE
    {\includegraphics[width=\textwidth]{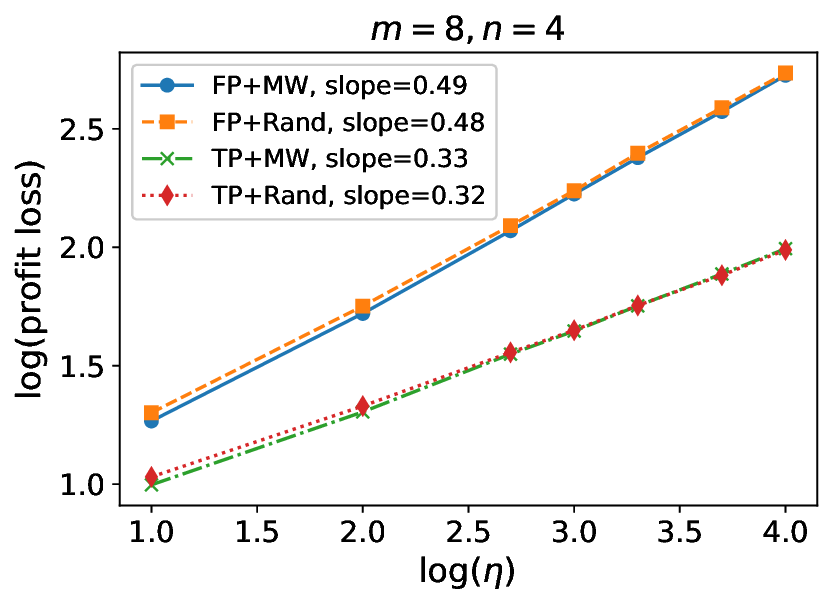}}
    {Log-log plot of profit losses under FP+MW, FP+Rand, TP+MW and TP+Rand for different $\eta$ when $m=4$ and $n=2$
    \label{fig: logprofitloss mnnotequal}}
    {}
  \end{minipage}
\end{figure}

Next, we analyze the general two-sided queues with multiple customer and server types.
We examine the profit losses under the following four different algorithms:
\begin{enumerate}
    \item FP+MW represents the fluid pricing (Eq \eqref{eq:fluid_pricing}) and max-weight matching (Eq~\eqref{eq:max-weight}) policy;
    \item FP+Rand represents the fluid pricing (Eq~\eqref{eq:fluid_pricing}) and randomized matching (Algorithm \ref{alg: random_non_empty}) policy;
    \item TP+MW represents the two-price policy (Eq \eqref{eq:two-price_policy-cust}) with max-weight matching (Eq \eqref{eq:max-weight});
    \item TP+Rand represents the two-price policy  (Eq \eqref{eq:two-price_policy-cust}) with randomized matching (Algorithm \ref{alg: random_non_empty}).
\end{enumerate}
In this numerical experiment, we first consider a setting where the number of servers and the number of customers are equal ($m=n$) and CRP condition (Condition \ref{condition: crp}) is satisfied. We assume the compatibility graph is given by
\begin{align*}
    E=\{(i,j)\in[n]\times[n]: j\in\{i+k\}\cup\{(i+k-n)^+\}, k=0,1,2,3 \}.
\end{align*}
The demand and supply curves are given by
\[
F_j(\lambda_j)=2-{\lambda_j}/{2},\ \forall j\in[m],
\quad\text{and}\quad
G_i(\mu_i)={\mu_i}/{2}\ \forall j\in[n],
\]
respectively.
We assume the unit holding cost is $s=1$.
The parameter of the fluid pricing policy is set to 
$    q_{\max}^{\eta}=2\sqrt{{\eta}/{n}}.$
The parameters of the two-price policy are chosen to be
$ \tau_{\max}^{\eta}=0$ and $\sigma^{\eta}=\eta^{2/3}n^{-1/3}$.
We report the profit loss for $\eta \in \{10, 100, 500, 1000, 2000, 5000, 10000\}$ when $m=n=6$ in Figure~\ref{fig: profitloss mnequal}. We find that when $\eta$ is larger, the profit loss of TP+MW grows the slowest, followed by the profit loss of TP+Rand, FP+MW, and FP+Rand. This result confirms the advantage of the two-price policy over the fluid pricing policy, as well as the advantage of the max-weight matching policy over the randomized matching policy. Figure~\ref{fig: logprofitloss mnequal} shows the same plot in logarithmic scale. Note that the slope of log-log plot in Figure~\ref{fig: logprofitloss mnequal} can be interpreted as the order of profit loss with respect to $\eta$. The fitted slopes of FP+MW and TP+MW are 0.51 and 0.33 respectively. This is consistent with Theorem~\ref{theo: fluidpricingpolicy} and Theorem~\ref{theorem: twoprice}, which state that FP+MW and TP+MW have the orders of profit loss with respect to $\eta$ of $O(\eta^{1/2})$ and $O(\eta^{1/3})$ respectively. Figure~\ref{fig: logprofitloss mnequal} shows that FP+MW and FP+Rand yield the same order of profit loss with respect to $\eta$ of approximately 1/2. Moreover, the two-price policy combined with either max-weight matching or randomized matching yields the same order of profit loss with respect to $\eta$ of approximately 1/3. That is, choosing max-weight or randomized matching does not affect order of profit loss with respect to $\eta$.
Our next experiment investigates how the profit loss changes with the number of customer and server types. We first consider a setting when the types of two sides are balanced ($m=n$).
Figure \ref{fig: profitloss eta10000 mnequal} shows the profit loss for $n \in \{4,6,8,\ldots,20\}$ when $\eta=10000$ (a large $\eta$ is chosen so that the asymptotic trend becomes clear). Figure \ref{fig: logprofitloss eta10000 mnequal} shows the same plot in logarithmic scale. It can be observed that the profit losses when a pricing policy is combined with the randomized matching policy grow faster than those when a pricing policy is combined with the max-weight matching policy as $n$ increases. In other words, the max-weight matching policy performs better than the randomized matching policy. Figure \ref{fig: logprofitloss eta10000 mnequal} suggests that the orders of profit loss with respect to $n$ of FP+MW and FP+Rand are 0.49 and 1.18 respectively,
which are close to those predicted by Theorem \ref{theo: fluid_max_weight_SSC}. Moreover, 
the orders of profit loss with respect to $n$ of TP+MW and TP+Rand are 0.34 and 1.28,
which are close to those predicted by Theorem \ref{theo: two_price_maxweight_SSC}.
Clearly, in both cases, the max-weight algorithm performances much better than the randomized matching algorithm for large $n$.
We also consider the setting where the number of server queues and the number of customer queues are not equal. Specifically, we assume that the number of server queues is twice as many as the number of customer queues, i.e., $m=2n$. The compatibility graph is given by
\begin{align*}
    E=\{(i,j)\in[2n]\times[n]: j\in\{i+k\}\cup\{(i+k-n)^+\}, k=0,1\}.
\end{align*}
The demand and supply curve are assumed to be
\begin{align*}
    F_j(\lambda_j)=6-\lambda_j, \,\forall j \in [m] \qquad\text{and}\qquad G_i(\mu_i)=\mu_i,\, \forall i\in[n],
\end{align*}
respectively. The parameters of pricing policy are similar to the previous case when $m=n$.
We report the profit loss for $\eta \in \{10, 100, 500, 1000, 2000, 5000, 10000\}$ when $m=4$ and $n=2$ in Figure \ref{fig: profitloss mnnotequal}. The result shows that when $\eta$ is larger, the profit loss when using the two-price policy grows significantly slower, compared to when using the fluid pricing policy. Moreover, we can observe that in this case the benefit of the max-weight matching policy over the randomized matching policy when combined with any pricing policy is negligible. The same plot in logarithmic scale (shown in Figure \ref{fig: logprofitloss mnnotequal}) shows that the fitted orders of profit loss with respect to $\eta$ of FP+MW and FP+Rand are 0.49 and 0.48 and those of TP+MW and TP+Rand are 0.33 and 0.32. This observation confirms the results from Theorem 
\ref{theo: fluidpricingpolicy} and Theorem \ref{theo: fluid_max_weight_SSC} as well as Theorem \ref{theorem: twoprice} and Theorem \ref{theo: two_price_maxweight_SSC}, which state that the orders of profit loss with respect to $\eta$ of FP+MW and FP+Rand are 1/2 and those of TP+MW and TP+Rand are 1/3 respectively.
Figure \ref{fig: profitloss eta10000 mnnotequal} shows the profit loss for $n \in \{4,6,8,\ldots,20\}$ and $m=2n$ when $\eta=10000$. Figure \ref{fig: logprofitloss eta10000 mnnotequal} shows the same plot in logarithmic scale. These figures shows the superiority of the max-weight policy over the randomized matching policy discussed in Section \ref{sec:MW vs Rand}. 
\section{Conclusion}
In this paper, we present a model of dynamic pricing and matching for two-sided queueing systems. The system is formulated as a Markov decision process, and a fluid approximation model is considered.  We presented a fluid pricing and max-weight matching policy and showed that it achieves $O(\sqrt{\eta})$ optimality rate. Furthermore, we proposed a dynamic pricing and max-weight policy, which achieves $O(\eta^{1/3})$ optimality rate. We also show that this scaling of $O(\eta^{1/3})$ matches the lower bound for a broad family of policies.
We also demonstrate the advantage of max-weight matching over randomized matching. Under the complete resource pooling condition, we show that max-weight matching achieves $O(\sqrt{n})$ and $O(n^{1/3})$ optimality rates for static and two-price policies, respectively,
where $n$ is the number of customer and server types. In comparison, the randomized matching policy may have an $\Omega(n)$ optimality rate.
\bibliographystyle{informs2014}
\bibliography{references}
\ECSwitch
\ECHead{E-Companion}
\renewcommand{\theHsection}{A\arabic{section}}
\begin{APPENDICES}{}
\section{MDP Analysis}
\label{sec:mdp-appendix}
\subsection{Proof of Proposition~\ref{prop: model-equivalence}}
\label{sec:proof_of_equivalence}
\proof{Proof of Proposition~\ref{prop: model-equivalence}.}
Since the delay-sensitive model and the holding cost model have the same states, actions, and transition rates, a control policy $\pi$ induces the same stationary distribution under the two models. Let 
$(\BFlambda^{\pi}(\BFq),\BFmu^{\pi}(\BFq))$ be the arrival rates under this policy when the system state is $\BFq$, and let $\E^{\pi}[\cdot]$ be the expectation operator under the stationary distribution induced by this policy.
The long-run average profit of the delay sensitive model is given by
\begin{align*}
     &\quad \E^{\pi}\left[\sum_{j=1}^m \lambda^{\pi}_j(\BFq)\left(F_j(\lambda^{\pi}_j(\BFq))-s^{(c)}_j \bar{w}_j^{(c)}\right)-\sum_{i=1}^n \mu_i^{\pi}(\BFq)\left(G_i(\mu^{\pi}_i(\BFq))+s^{(s)}_i \bar{w}_i^{(s)}\right)\right] \\
     &= \E^{\pi}\left[\sum_{j=1}^m \lambda^{\pi}_j(\BFq)F_j\left(\lambda^{\pi}_j(\BFq)\right)-\sum_{i=1}^n \mu^{\pi}_i(\BFq)G_i\left(\mu^{\pi}_i(\BFq)\right)\right]-\sum_{j=1}^m s^{(c)}_j\E[\lambda^{\pi}_j(\BFq)]\bar{w}_j^{(c)} \\
     & \quad \quad -\sum_{i=1}^n s^{(s)}_i
     \E^{\pi}[\mu^{\pi}_i(\BFq)]\bar{w}_i^{(s)}  \\
     &=\E^{\pi}\left[\sum_{j=1}^m \lambda^{\pi}_j(\BFq)F_j\left(\lambda^{\pi}_j(\BFq)\right)-\sum_{i=1}^n \mu_i^{\pi}(\BFq)G_i\left(\mu_i^{\pi}(\BFq)\right)\right]-\sum_{j=1}^m s^{(c)}_j\E^{\pi} [ q_j^{(c)}]-\sum_{i=1}^n s^{(s)}_i  \E^{\pi}[q_i^{(s)}].
\end{align*}
The last step applies Little's Law:
$\E[\lambda^{\pi}_j(\BFq)]\bar{w}_j^{(c)} = \E^{\pi} [ q_j^{(c)}]$,
$\E^{\pi}[\mu^{\pi}_i(\BFq)]\bar{w}_i^{(s)}=\E^{\pi}[q_i^{(s)}]$.
Note that the last line of the above equation is the long-run average profit under the holding cost model, so the proof is complete.
$\Halmos$
\endproof
\subsection{Monotonicity of the Optimal Prices (Single-Link Two-Sided Queues)}
\label{sec:single-link}

In this section, we consider a special case with $n=1$ and $m=1$, i.e., a single-link two-sided queue given in Figure~\ref{fig:ex_single_link}. The goal of this section is to analyze the optimal pricing policy for this special case, which will motivate our pricing policies for more complex systems.

\begin{figure}[!htbp]
\TABLE{A single-link two-sided queue. \label{fig:ex_single_link}}
    {\begin{tabular}{c}
         \begin{tikzpicture}[scale=0.7]
\draw[black, very thick] (0,0) -- (2,0) -- (2,1) -- (0,1);
\draw[black, very thick] (8,0) -- (6,0) -- (6,1) -- (8,1);
\node[black] at (1,1.5) {\footnotesize Customer};
\node[black] at (7,1.5) {\footnotesize Server};
\draw[black,thick]  (2.25, 0.5) edge[<->]  (5.75, 0.5);
\end{tikzpicture}
    \end{tabular}}
    {}
\end{figure}

In single-link systems, there is no incentive for the system operator to hold customers or servers. Whenever possible, we should match the incoming arrival immediately. 
Thus, at any point of time, there can only be either customers or servers waiting in the system (Here, we allow to match an incoming arrival as well for simplicity, which is slightly incompatible with the definition $\BFx \in X(\BFq)$). This enables up to reduce the state space by letting $q=q^{(c)}-q^{(s)}$, the difference between the number of customers and servers waiting in the system. Note that, $q$ can be either positive or negative. Using $q$ as the system state, the Bellman equation \eqref{MDP} becomes
\begin{align}
    h(q)=&\max_{\substack{\lambda \in [0, \lambda_{\max}] \\ \mu \in [0, \mu_{\max}]}} \left[\frac{F(\lambda)\lambda-G(\mu)\mu}{c} -\frac{s|q|}{c} -\frac{\gamma}{c} + \frac{\lambda}{c}h(q + 1) +\left(1-\frac{\mu+\lambda}{c}\right)h(q)+\frac{\mu}{c}h(q-1)\right], \forall q \in S \label{eq: singlelinkbellman}
\end{align}
where $c$ is a uniformization parameter (see Definition~\ref{ass: upperboundedrates}) and $S$ denotes the state space. 
(We omit the subscripts for customer and server types, since $n=m=1$). We now present the monotonicity result below.
\begin{proposition}\label{theo: monotonicity}
For a single-link two-sided queue, there exists an optimal pricing policy $\BFp(q) = (p^{(s)}(q), p^{(c)}(q))$, where both the server price $p^{(s)}(q)$ and the customer price $p^{(c)}(q)$ increases monotonically with the system state $q$.
\end{proposition} 
This result motivates us to search for the optimal pricing policy in the restricted space of monotonic pricing policies, which will be presented in Section~\ref{sec: LP approx}.

To prove Proposition~\ref{theo: monotonicity}, we first show that the difference of the optimal bias functions, $\Delta h(q)\overset{\Delta}{=}h(q)-h(q-1)$, is monotonically decreasing in $q$. As the relative value iteration may not converge in general for average cost MDP with countable state space, we work with the discounted MDP with the following Bellman equation with discount factor $\alpha \in (0, 1)$:
\begin{align*}
    V_{\alpha}(q) = \max_{\substack{\lambda \in [0, \lambda_{\max}] \\ \mu \in [0, \mu_{\max}]}} \left[\frac{F(\lambda)\lambda-G(\mu)\mu}{c} -\frac{s|q|}{c} + \alpha\frac{\lambda}{c}V_\alpha(q + 1) +\alpha\left(1-\frac{\mu+\lambda}{c}\right)V_\alpha(q)+\alpha\frac{\mu}{c}V_\alpha(q-1)\right]
\end{align*}
We consider the value iteration method to compute $V_\alpha(q)$.
The value iteration starts with an arbitrary initial value $V_{\alpha, 0}(q)$; therefore we choose an initial value function such that $\Delta V_{\alpha, 0}(q)$ is decreasing. In each iteration $k=1,2,\cdots$, for all $q \in S$, the value iteration algorithm updates the value function by
\begin{align}
    V_{\alpha, k+1}(q)=&\max_{\lambda \in [0, \lambda_{\max}], \mu \in [0, \mu_{\max}]} \left[\frac{F(\lambda)\lambda-G(\mu)\mu}{c} -\frac{s|q|}{c} + \alpha\frac{\lambda}{c}V_{\alpha, k}(q + 1)\right. \nonumber\\
    &\left.+\alpha\left(1-\frac{\mu+\lambda}{c}\right)V_{\alpha, k}(q)+\alpha\frac{\mu}{c}V_{\alpha, k}(q-1)\right] \label{eq:discounted_bellman}
\end{align}
As $k\to \infty$, we have $V_{\alpha, k}(q) \to V_\alpha(q)$ (e.g., see Theorem 3.4 of \citet{feinberg2018convergence}). Moreover, we will then take $\alpha \to 1$ to obtain results on the optimal bias function.
These claims are formalized in the following lemma.

\begin{lemma} \label{lemma: monotonic}
The difference of the optimal bias function $\Delta h^{*}(q)$ that solves \eqref{eq: singlelinkbellman} is monotonically decreasing in $q$.
\end{lemma}

\proof{Proof of Lemma~\ref{lemma: monotonic}.}
The proof is by induction using value iteration for the discounted problem. In the base case, as we can choose any initial value function in the value iteration algorithm, we pick a value function such that $\Delta V_{\alpha, 0}(q)$ is monotonically decreasing in $q$, say, $V_{\alpha, 0}(q)=0$. 

Suppose we are at iteration $k$.
Assume that $\Delta V_{\alpha, k}(q)$ is monotonically decreasing for all $q$.
We will now calculate $\Delta V_{\alpha, k+1}(q+1)-\Delta V_{\alpha, k+1}(q)$ and show that it is nonpositive. 
The value iteration step can be rewritten as
\begin{align}
    V_{\alpha, k+1}(q)={}& \frac{\mathcal{R}(\mu^{*}(q),\lambda^{*}(q))-s|q|}{c}+ \alpha\frac{\lambda^\star(q)}{c}\Delta V_{\alpha, k}(q + 1)+\alpha V_{\alpha, k}(q)-\alpha\frac{\mu^\star(q)}{c}\Delta V_{\alpha, k}(q),
\end{align}
where $(\lambda^\star(q), \mu^\star(q))$ is the maximizer of \eqref{eq:discounted_bellman} with the value function $V_{\alpha, k}$ and $\mathcal{R}(\mu,\lambda)$ is equal to $F(\lambda)\lambda-G(\mu(q))\mu$. Thus, we have
\begin{align}
    &\Delta V_{\alpha, k+1}(q+1) - \Delta V_{\alpha, k+1}(q) \nonumber\\
    ={}& \alpha\left(\Delta V_{\alpha, k}(q+1) - \Delta V_{\alpha, k}(q)\right) 
    \nonumber\\
    &+\left[\alpha\lambda^\star(q+1)\Delta V_{\alpha, k}(q + 2)-\alpha\mu^\star(q+1)\Delta V_{\alpha, k}(q+1)+\mathcal{R}(\mu^{*}(q+1),\lambda^{*}(q+1))-s|q+1|\right]/c\nonumber \\
    &-2[\alpha\lambda^\star(q)\Delta V_{\alpha, k}(q + 1)-\alpha\mu^\star(q)\Delta V_{\alpha, k}(q)+\mathcal{R}(\mu^{*}(q),\lambda^{*}(q))-s|q|]/c \nonumber\\
    &+[\alpha\lambda^\star(q-1)\Delta V_{\alpha, k}(q)-\alpha\mu^\star(q-1)\Delta V_{\alpha, k}(q-1)+\mathcal{R}(\mu^{*}(q-1),\lambda^{*}(q-1))-s|q-1|]/c  \label{big}
\end{align}
Because  $(\lambda^{*}(q),\mu^{*}(q))$ is the maximizer for state $q$, we have
\begin{align}
    \mathcal{R}(\mu^{*}(q),&\lambda^{*}(q))+\alpha\lambda^{*}(q)\Delta V_{\alpha, k}(q+1)-\alpha\mu^{*}(q)\Delta V_{\alpha, k}(q)\geq \mathcal{R}(\mu^{*}(q+i),\lambda^{*}(q+i))\nonumber\\
    &+\alpha\lambda^{*}(q+i)\Delta V_{\alpha, k}(q+1)-\alpha\mu^{*}(q+i)\Delta V_{\alpha, k}(q) \quad \text{for } i \in \{-1,1\}. \label{geq}
\end{align}
Plugging \eqref{geq} to \eqref{big}, and using the induction hypothesis and noting that $2|q| - |q+1|-|q-1| \leq 0$, we get 
\begin{align}
     &\Delta V_{\alpha, k+1}(q+1)-\Delta V_{\alpha, k+1}(q) \leq \alpha \lambda^{*}(q+1)\big(\Delta V_{\alpha, k}(q+2)
    -\Delta V_{\alpha, k}(q+1)\big)/c 
     \nonumber \\
     & + \alpha (1 -\mu^{*}(q+1)/c-\lambda^{*}(q-1)/c)\big(\Delta V_{\alpha, k}(q+1)-\Delta V_{\alpha, k}(q)\big)  \nonumber \\
    &\; +\alpha\mu^{*}(q-1)\big(\Delta V_{\alpha, k}(q)-\Delta V_{\alpha, k}(q-1)\big)/c. \nonumber 
\end{align}
Since the uniformization constant is chosen as $c \geq \lambda_{\max}+\mu_{\max}$, we have $1-\mu^{*}(q+1)/c -\lambda^{*}(q-1)/c \geq 0$ for all $q$. By the induction hypothesis in iteration $k$, $\Delta V_{\alpha, k}(q+1)-\Delta V_{\alpha, k}(q) \leq 0$ for all $q$. Thus, the last line of the above inequality is nonpositive.

This proves that the relative value iteration step preserves the monotonicity of $\Delta V_{\alpha, k}(q)$ . As $k\to\infty$, $\Delta V_{\alpha, k}(q)$  converges to $\Delta V_{\alpha}^{*}(q)$ (e.g., see Theorem 3.4 of \citet{feinberg2018convergence}).
Since the limit of a decreasing function is decreasing, $\Delta V_{\alpha}^{*}(q)$ is monotonically decreasing. Lastly, by \citet[Theorem 3.2]{cavazos1989weak}, there exists a subsequence of discount factors $\{\alpha_n\}$ and reference states $\{q_n\}$ such that 
$$
V_{\alpha_n}^\star(q) - V_{\alpha_n}^\star(q_n) \to h^\star(q) \quad \text{as} \quad \alpha_n \to 1,
$$
where, the convergence above is  pointwise. Additive normalization and pointwise convergence preserve the monotonicity property; thus, the resulting bias is such that $\Delta h^\star(q)$ is monotonically decreasing. This completes the proof.
$\Halmos$
\endproof

\proof{Proof of Proposition~\ref{theo: monotonicity}.}
The equation to be maximize in \eqref{eq: singlelinkbellman} is separable in $\lambda$ and $\mu$.
The domain of $\lambda$ is $[0,\lambda_{\max}]$ and the domain of $\mu$ is $[0,\mu_{\max}]$.
By Assumption 2, it is  also concave in $\lambda$ and $\mu$. Thus, $(\lambda^*, \mu^*)$ is an optimal solution to \eqref{eq: singlelinkbellman} either if it is on the boundary of the domain or 
if it satisfies the first order necessary condition.

First, we consider the boundary case. We will show that if $\lambda^{*}(q_0)=0$ for some $q_0$,
then $\lambda^{*}(q)=0$ for all $q > q_0$. Suppose $\lambda^{*}(q_0)=0$. Then we have
\[
F(\lambda)\lambda+\Delta h^{*}(q_0)\lambda \ \leq\  0, \quad \forall \lambda \in [0,\lambda_{\max}].
\]
By Lemma~\ref{lemma: monotonic},
$\Delta h^{*}(q)$ is decreasing in $q$, so
\begin{align*}
 F(\lambda)\lambda+\Delta h^{*}(q_0+k)\lambda
 \ &=\ F(\lambda)\lambda+\Delta h^{*}(q_0)\lambda+\left(\Delta h^{*}(q_0 + k)-\Delta h^{*}(q_0)\right)\lambda \\
 \ &\leq\ \left(\Delta h^{*}(q_0 + k) -\Delta h^{*}(q_0)\right)\lambda \\
 \ &\leq\ 0, \quad \forall k \geq 1.
\end{align*}
The above inequality implies that $\lambda^{*}(q_0+k)=0$. Similarly, if $\lambda^{*}(q_0)=\lambda_{\max}$ then $\lambda^{*}(q)=\lambda_{\max}$ for all $q<q_0$. To see this, suppose
\[
    F(\lambda)\lambda+\Delta h^{*}(q_0)\lambda \leq F(\lambda_{\max})\lambda_{\max}+\Delta h^{*}(q_0)\lambda_{\max}, \quad \forall \lambda \in [0,\lambda_{\max}]. 
\]
For any $k\geq 1$,
adding $\left(\Delta h^{*}(q_0-k) - \Delta h^{*}(q_0)\right)\lambda$ on both sides of the above inequality, we get
\begin{align*}
    F(\lambda)\lambda+\Delta h(q_0-k)\lambda &\leq F(\lambda_{\max})\lambda_{\max}    + \Delta h^{*}(q_0)\lambda_{\max} + \left(\Delta h^{*}(q_0-k) - \Delta h^{*}(q_0)\right)\lambda \\
   &\leq F(\lambda_{\max})\lambda_{\max}    + \Delta h^{*}(q_0)\lambda_{\max} + \left(\Delta h^{*}(q_0-k) - \Delta h^{*}(q_0)\right)\lambda_{\max}  \\
   &= F(\lambda_{\max})\lambda_{\max} +\Delta h^{*}(q_0-k)\lambda_{\max},
\end{align*}
where the second inequality follows as $\lambda_{\max}\geq\lambda$ and $\Delta h^{*}(q_0-k)-\Delta h^{*}(q_0) \geq 0$. Thus, $\lambda^{*}(q)=\lambda_{\max}$ is the maximizer for all $q<q_0$. 
Using the same argument, we can show that if $\mu^{*}(q_0)=\mu_{\max}$ then $\mu(q)^{*}=\mu_{\max}$ for all $q \geq q_0$; if $\mu^{*}(q_0)=0$ then $\mu^{*}(q)=0$ for all $q \leq q_0$. 

We now consider the case when the optimal solution in \eqref{eq: singlelinkbellman} is in the interior of the domain, i.e., $\lambda_{\max}>\lambda^{*}(q)>0$ and $\mu_{\max}>\mu^{*}(q)>0$. 
By the first-order condition, we have
\begin{align*}
    [F(\lambda^{*}(q))\lambda^{*}(q)]'+\Delta h^{*}(q+1)&=0, \\
    [G(\mu^{*}(q))\mu^{*}(q)]'+\Delta h^{*}(q)&=0. 
\end{align*}
Note that the derivatives are well defined because $F$ and $G$ are differentiable by Assumption \ref{ass: monotonic}.
Since $\Delta h^{*}(q)$ is monotonically decreasing in $q$ by Lemma \ref{lemma: monotonic}, it implies that  $[F(\lambda^{*}(q))\lambda^{*}(q)]'$ and $[G(\mu^{*}(q))\mu^{*}(q)]'$ are monotonically increasing in $q$. 
By Assumption \ref{ass: concavity}, $[F(\lambda)\lambda]'$ is decreasing in $\lambda$ and $[G(\mu)\mu]'$ is increasing in $\mu$. Thus, by the chain rule, $\lambda^{*}(q)$ is monotonically decreasing in $q$ and $\mu^{*}(q)$ is monotonically increasing in $q$ when $\lambda_{\max}>\lambda^{*}>0$ and $\mu_{\max}>\mu^{*}>0$. 

Combining the boundary case and the interior case, we prove that $\lambda^{*}(q)$ is monotonically decreasing in $q$. As the demand curve is monotonically decreasing by Assumption \ref{ass: monotonic}, the optimal customer price $(p^{(c)})^{*}(q)$ is monotonically increasing. Also, $\mu^{*}(q)$ is monotonically increasing in $q$. As the supply curve is monotonically increasing by Assumption \ref{ass: monotonic}, the optimal server price $(p^{(s)})^{*}(q)$ is monotonically increasing. $\Halmos$
\endproof

\subsection{LP-based Approximation Algorithm} \label{sec: LP approx}
Throughout this section, we assume the queue length is bounded by some fixed constant. Specifically, we assume that the state space is given by $S = \{\BFq \in \mathbb{Z}^{m+n}: \BFzero \leq \BFq \leq \BFone_{m+n} q_{\max}$\} for some $q_{\max} <\infty$.
The Bellman equation \eqref{MDP} for the uniformized DTMDP defined in the Section~\ref{sec: model} can be rewritten as a semi-infinite linear program \citep[see e.g.][]{bertsekas2012dynamic}:
\begin{align}
    \min_{(\gamma,\BFh)}\ &\  \gamma    \label{eq: obj-LP}\\
\text{subject to }\ &
    \ \gamma \geq  \mathcal{R}(\BFq,\BFz)+ c\E[h(\BFq') \mid \BFq,\BFz]-c h(\BFq) \quad
    \forall \BFq,\BFz \in S \times Z(\BFq), \label{eq: cons-LP}
\end{align}
where $\mathcal{R}(\BFq,\BFz)$ is the expected profit (see Eq \eqref{revenue}) and $\E[h(\BFq') \mid \BFq,\BFz]$ is the expectation of the bias function after one transition.
However, it is difficult to solve this linear program directly, since it has a separate decision variable $h(\BFq)$ for each state $\BFq$, and
the number of constraints in  \eqref{eq: cons-LP} is equal to the number of state-action pairs.
Due to the curse of dimensionality, the state space will increase exponentially with the customer and server types. 
\subsubsection{Polynomial Bias Function Approximation.}
We consider an approximation of the MDP in the value space. In particular, we approximate the bias function $h(\BFq)$ by a polynomial function:
\begin{align}
    h(\BFq) \approx \sum_{l=1}^r \left(\sum_{i=1}^n b_{l_1}^{(s)}(q_i^{(s)})^l+\sum_{j=1}^m b_{l_1}^{(c)}(q_j^{(c)})^l\right)\overset{\Delta}{=}\sum_{l=1}^r \inner{\BFb_l}{\BFq^l}, \label{eq: upperboundvaluefunction}
\end{align}
for some degree $r \in \mathbb{Z}_+$. Here, $\BFb_l$ is a vector $(b_{l_1}^{(s)},\hdots,b_{l_n}^{(s)},b_{l_1}^{(c)},\hdots,b_{l_m}^{(c)})$ for all $l \in [r]$.
If we apply this approximation to the linear programming formulation \eqref{eq: obj-LP}, it in fact leads to an upper bound approximation of the original MDP (see e.g.\ \cite{adelman2007dynamic}). 
\begin{proposition} \label{theo: valuefunctionapprox}
Suppose $h(\BFq)$ is replaced with $\sum_{l=1}^r \inner{\BFb_l}{\BFq^l}$ in the optimization problem \eqref{eq: obj-LP}.
The optimal objective value is an upper bound of the optimal average cost of the original MDP.
\end{proposition}
\proof{Proof of Proposition~\ref{theo: valuefunctionapprox}.}
Rewriting the Bellman equation using the approximation of the bias function gives us the following semi-infinite linear program:
\begin{align}
    \min \ &\ \gamma   \nonumber \\
    \text{subject to}\ & \ \gamma \geq  \mathcal{R}(\BFq,\BFz)+ c\E\left[\sum_{l=1}^r \inner{\BFb_l}{\BFq'}\mid \BFq,\BFz\right]-c \sum_{l=1}^r \inner{\BFb_l}{\BFq^l} \quad
    \forall \BFq,\BFz \in S \times Z(\BFq). \nonumber
\end{align}
The decision variables in the above optimization problem are $\gamma$ and $\BFb_l \ \forall l \in [r]$. 
Let the optimal solution to the above optimization problem be $\gamma^*$ and $\BFb_l^{*}$ $(\forall l \in [r])$. Now define
\begin{align*}
    h(\BFq)=\sum_{l=1}^{r}\inner{\BFb_l^{*}}{\BFq^l}\ \forall \BFq \in S.
\end{align*}
Since $h(\BFq)$ and $\gamma^{*}$ are a feasible solution to the optimization problem \eqref{eq: obj-LP}, the optimal value of \eqref{eq: obj-LP} is less than or equal to $\gamma^{*}$. $\Halmos$
\endproof
By approximating the bias function by a polynomial function of degree $r$, the number of variables is reduced from $q_{\max}^{m+n}$ to  $(m+n) \times r$. We will later see that this approximation reduces the computational time of the semi-infinite linear program drastically.
As the degree of the polynomial $r$ increases, the upper bound becomes tighter. Therefore, $r$ can be selected to balance the trade-off between approximation accuracy and the computational time.
\subsubsection{Matching Policy under Value Function Approximation.} 
We now focus on matching policies when $h(\BFq)$ is approximated by polynomial functions, in particular, linear and quadratic functions.
We denote by $(\BFlambda^{*}(\BFq), \BFmu^{*}(\BFq))$ the optimal pricing decision for state $\BFq$ under bias function approximation.
We can rewrite the Bellman equation as follows:
\begin{align*}
     \frac{\gamma}{c}\ = \frac{\mathcal{R}(\BFq,\BFlambda^{*}(\BFq),\BFmu^{*}(\BFq))}{c} +\max_{\BFx \in X(\BFq)}\left\{\E\left[h(\BFq')\mid\BFq,\BFz\right]-h(\BFq) \right\}, \quad \forall \BFq \in S.
\end{align*}
The optimal matching decision is given by
\begin{equation}\label{eq:Bellman-matching-desision}
    \BFx^{*}(\BFq)=\argmax_{\BFx \in X(\BFq)}\left\{\E\left[h(\BFq')\mid\BFq,\BFz\right]-h(\BFq) \right\} \quad \forall \BFq \in S.
\end{equation}
Suppose we approximate the bias functions by linear functions, say, $h(\BFq)=\inner{-\BFb}{\BFq}$ for some $\BFb \in \mathbb{R}^{m+n}$, then we have
\begin{align*}
    \BFx^{*}(\BFq)&=\argmax_{\BFx \in X(\BFq)}\left\{\E\left[h(\BFq')\mid\BFq,\BFz\right]-h(\BFq) \right\}  \\
    &=\argmax_{\BFx \in X(\BFq)}\left\{\sum_{j=1}^m \lambda_j^{*}(\BFq) \left(h(\BFq+\BFe_j^{(c)}-\BFx)-h(\BFq)\right)+\sum_{i=1}^n \mu_i^{*}(\BFq) \left(h(\BFq+\BFe_i^{(s)}-\BFx)-h(\BFq)\right) \right\} \\
    &=\argmax_{\BFx \in X(\BFq)}\left\{\left(\sum_{j=1}^m \lambda_j^{*}(\BFq) + \sum_{i=1}^n \mu_i^{*}(\BFq)\right)
    \inner{\BFb}{\BFx} \right\} \\
    &=\argmax_{\BFx \in X(\BFq)}\left\{
    \inner{\BFb}{\BFx} \right\}.
\end{align*}
Note that we ignore terms that are independent of $\BFx$ in the second to last equality.
Thus, the optimal matching decision is given by a maximum-weight bipartite matching problem
with fixed (vertex) weights $\BFb$ and feasibility set $X(\BFq)$.
Alternatively, suppose we approximate the bias function by a \emph{quadratic} function of queue lengths given by
$
    h(\BFq)=\inner{-\BFone}{\BFq^2}. 
$
By Eq~\eqref{eq:Bellman-matching-desision}, the optimal matching decision is given by
\begin{align}
     \BFx^{*}(\BFq) =&\argmax_{\BFx \in X(\BFq)}\left\{\E\left[h(\BFq')\mid\BFq,\BFz\right]-h(\BFq) \right\}  \nonumber\\
     ={}&\argmax_{\BFx \in X(\BFq)}\left\{\sum_{j=1}^m \lambda_j^{*}(\BFq) \left(h(\BFq+\BFe_j^{(c)}-\BFx)-h(\BFq)\right)+\sum_{i=1}^n \mu_i^{*}(\BFq) \left(h(\BFq+\BFe_i^{(s)}-\BFx)-h(\BFq)\right) \right\}\nonumber \\
     ={}&\argmax_{\BFx \in X(\BFq)}\left\{\left(\sum_{j=1}^m \lambda_j^{*}(\BFq)+ \sum_{i=1}^n \mu_i^{*}(\BFq) \right)
     \left(
     2\inner{\BFq}{\BFx}-\inner{\BFone}{\BFx^2} 
     \right) + 
     \sum_{j=1}^m 2\lambda_j^{*}(\BFq)x^{(c)}_j +
     \sum_{i=1}^n 2\mu_i^{*}(\BFq)x^{(s)}_i
     \right\}.\nonumber
\end{align}   
We consider two cases. If $X(\BFq)=\{\BFzero\}$, namely, there is no feasible matching given state $\BFq$, then we trivially have $\BFx^*(\BFq)=\BFzero$. 
Otherwise, if $X(\BFq)\neq \{\BFzero\}$, then there exists a matching decision $\BFx \in X(\BFq)$ such that $2\inner{\BFq}{\BFx}-\inner{\BFone}{\BFx^2} >0$, which means that $\BFx=\BFzero$ cannot be a maximizer of the above equation. This means that under quadratic value function approximation, customers and servers should be matched instantly when they arrive.  Recall that there is at most one new arrival per period in the unformized DTMDP. So, we either have $\inner{\BFone}{\left(\BFx^*(\BFq)\right)^2}=0$ when $X(\BFq)=\{\BFzero\}$ or  $\inner{\BFone}{\left(\BFx^*(\BFq)\right)^2}=1+1=2$ when $X(\BFq)\neq\{\BFzero\}$. 
In either case, we can remove the term $-\inner{\BFone}{\BFx^2}$ in the above equation and get
\begin{align}
      \BFx^{*}(\BFq){=}{}& \argmax_{\BFx \in X(\BFq)}\left\{\left(\sum_{j=1}^m \lambda_j^{*}(\BFq)+ \sum_{i=1}^n \mu_i^{*}(\BFq) \right)
     \left(
     2\inner{\BFq}{\BFx} 
     \right) + 
     \sum_{j=1}^m 2\lambda_j^{*}(\BFq)x^{(c)}_j +
     \sum_{i=1}^n 2\mu_i^{*}(\BFq)x^{(s)}_i
     \right\} \label{eq: quadraticapproxmaxweight}\\
     {\approx}{}& \argmax_{\BFx \in X(\BFq)}\left\{ 
     \left(\sum_{j=1}^m \lambda_j^{*}(\BFq)+ \sum_{i=1}^n \mu_i^{*}(\BFq) \right)
     \left(
     2\inner{\BFq}{\BFx} 
     \right) \right\} \nonumber \\
     {=}{}& \argmax_{\BFx \in X(\BFq)}\left\{ 
     \inner{\BFq}{\BFx} 
      \right\}, \nonumber
\end{align}
where the approximation step is motivated by the fact that the value of $(\sum_{j=1}^m \lambda_j^{*}(\BFq)+ \sum_{i=1}^n \mu_i^{*}(\BFq))$ dominates either $\lambda^*_j(\BFq)$ or $\mu^*_{i}(\BFq)$, especially when $n$ and $m$ are large. This approximation leads to the max-weight matching policy defined in Section~\ref{subsec:max-weight matching}.
\subsubsection{Constraint Generation Algorithm.}
In this section, we use the constraint generation technique to solve the optimization problem \eqref{eq: obj-LP} given approximated bias functions. The constraint generation steps are described in Algorithm~\ref{alg1}.
\begin{algorithm}
\caption{Constraint Generation with Bias Function Approximation}\label{alg1}
{\fontsize{10}{13}\selectfont
\begin{algorithmic}[2]
\STATE \textbf{Initialization:} $\BFb^0_l=0,\ \forall l \in [r], \gamma^0=-\infty, k=0, \epsilon = +\infty$
\STATE \textbf{Initialize Master LP:} $\{ \min \gamma, \text{subject to } \varnothing\}$ \hfill{Master Problem: $(LP^0)$}
\WHILE{$\epsilon > \text{tolerance}$}
    \STATE $T^k(\BFq,\BFz) \overset{\Delta}{=} \mathcal{R}(\BFq,\BFz)+c\E\left[\sum_{l=1}^r \inner{\BFb_l^k}{(\BFq^l)'}\mid \BFq,\BFz\right]$   $-c\sum_{l=1}^r\inner{\BFb_l^k}{\BFq^l}$
        \STATE $\delta^k(\BFq) \leftarrow \max_{\BFz\in Z(\BFq)}T^k(\BFq,\BFz)$ \COMMENT{Sub-Problem}
        \IF{$\delta^k(\BFq) > \gamma^k$}
            \STATE $\BFq^k,\BFz^k=\arg\max_{\BFq,\BFz}T^k(\BFq,\BFz)$
            \STATE add the constraint $\gamma \geq T^k(\BFq^k,\BFz^k)$ to the master LP \hfill{ $(LP^{k})$}
        \ENDIF
    \STATE let $(\gamma^{k+1},\BFb^{k+1}_l \ \forall l \in [r])$ be the solution to $(LP^k)$ \COMMENT{Master-Problem}
    \STATE $\epsilon 
    \leftarrow \gamma^{k+1}-\gamma^k$
    \STATE $k \leftarrow k+1$
\ENDWHILE
\STATE \textbf{Output:} $\gamma \leftarrow \gamma^k,\BFb_l \leftarrow \BFb^k_l, \ \forall l \in [r]$
\end{algorithmic}
}
\end{algorithm}
The algorithm starts with some initial values of the weights $\BFb$, say $\BFb=\BFzero$. The algorithm also maintains a master LP that approximates \eqref{eq: obj-LP}.
In each iteration, the algorithm finds a violating constraint by solving the following sub-problem
\begin{align}
    \max_{\BFq \leq q_{\max}\BFone_{m+n}, \BFz \in Z(\BFq)}\mathcal{R}(\BFq,\BFz)+c\E\left[\sum_{l=1}^r \inner{\BFb_l^k}{({\BFq}^l)'}\mid \BFq,\BFz \right]-c\sum_{l=1}^r\inner{\BFb_l^k}{\BFq^l}. \label{eq: subproblem}
\end{align}
If the optimal value of the above subproblem is larger than the optimal value of the master-problem, then a violating constraint is found and added to the master-problem. We then solve the master-problem to get the updated values of $\BFb_l \ \forall l \in [r]$ and $\gamma$. This process is repeated until either no violating constraint is found or the improvement is less than some toleration.
\subsubsection{Simulation Results.}
We present some simulation results obtained by approximating the bias function by polynomial functions. 
Figures~\ref{fig: polyapproxvsMDP0.01} and \ref{fig: polyapproxvsMDP0.05} show the pricing policy obtained by linear and quadratic approximation of the bias function, as well as the optimal pricing policy. We observe that both linear and quadratic approximation result in pricing policies that are monotonic, which is consistent with the monotonic structure of the optimal pricing policy. As the penalty coefficient $s$ increases, the approximation error becomes larger. 
\begin{figure}[b]
  \centering
  \begin{minipage}[b]{0.48\textwidth}
  \FIGURE
    {\includegraphics[width=\textwidth]{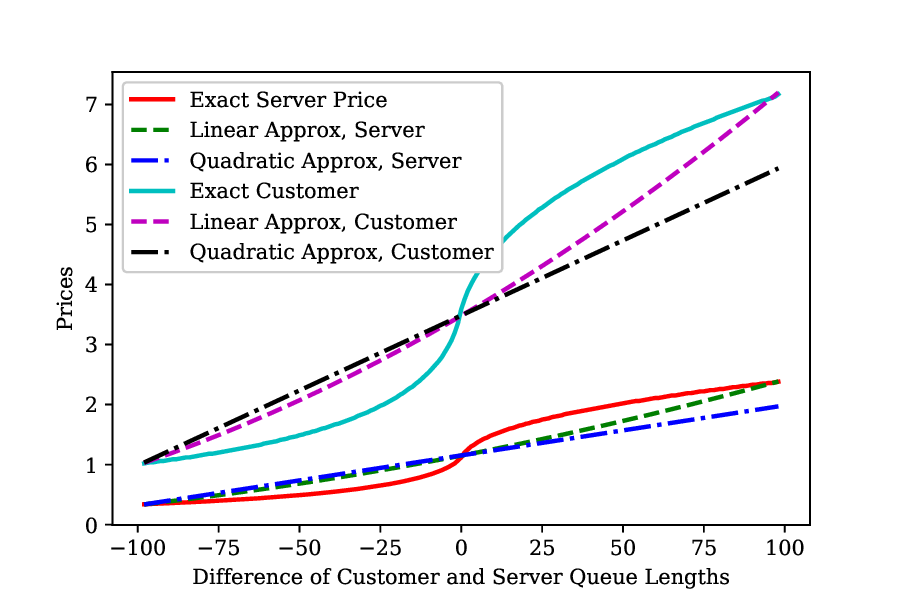}}
    {Comparison of pricing policies with bias approximation ($s=0.05$).
    \label{fig: polyapproxvsMDP0.05}}
    {}
  \end{minipage}
  \hfill
    \begin{minipage}[b]{0.48\textwidth}
  \FIGURE
      {\includegraphics[width=\textwidth]{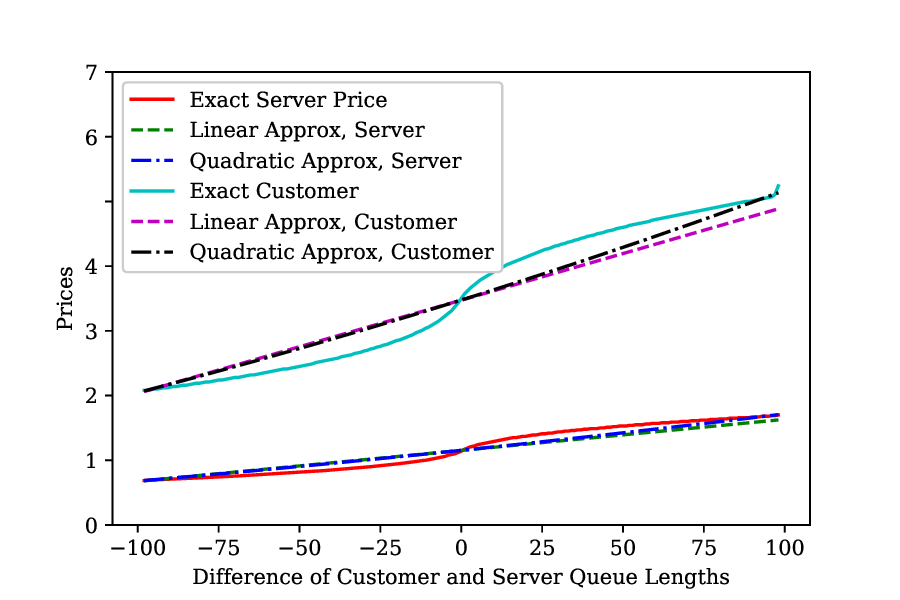}}
    {Comparison of pricing policies with bias approximation ($s=0.01$).
    \label{fig: polyapproxvsMDP0.01}}
    {}
  \end{minipage}
\end{figure}
In the linear approximation case, we also compare the solution from the constraint generation algorithm to the least squares fit of the exact bias function. 
Our experiment uses both constant elasticity supply/demand functions and linear supply/demand function.
The result is presented in Tables~\ref{tab: loglog} and
\ref{tab: linear}. 
We observe again that the approximation error increases with $s$.  
We can also see that the optimal objective value $\gamma$ obtained by constraint generation is an upper bound on the optimal value of the MDP. This verifies Proposition~\ref{theo: valuefunctionapprox}.
\begin{table}
\TABLE
{Comparison of constraint generation solution with the optimal solution with constant elasticity supply and demand curves.
\label{tab: loglog}}
{\begin{tabular}{c|ccc|ccc|cc}
\hline
\multirow{2}{*}{$s$} & \multicolumn{3}{c|}{Cons. Generation} & \multicolumn{3}{c|}{Least Squares Fit} & \multicolumn{2}{c}{\% Error} \\
                     & $b_0$       & $b_1$       & $\gamma$       & $b_0$        & $b_1$        & $\gamma$        & $b_0$         & $b_1$         \\ \hline 
0.01                 & -1.73       & -0.007      & 3.07      & -1.69        & -0.008       & 3.03       & 2\%           & 10\%          \\
0.02                 & -1.73       & -0.010      & 3.07      & -1.80        & -0.013       & 3.03       & 4\%           & 25\%          \\
0.05                 & -1.73       & -0.013      & 3.06      & -1.89        & -0.019       & 3.03       & 9\%           & 34\%          \\
0.1                  & -1.73       & -0.014      & 3.06      & -2.05        & -0.025       & 2.95       & 16\%          & 43\%          \\
0.2                  & -1.73       & -0.016      & 3.06      & -2.59        & -0.034       & 2.8        & 24\%          & 53\%          \\
0.5                  & -1.73       & -0.017      & 3.06      & -2.59        & -0.051       & 2.48       & 33\%          & 67\%          \\ \hline
\end{tabular}}
{}
\end{table}
\begin{table}
\TABLE
{Comparison of constraint generation solution with the solution with linear supply and demand curves.
\label{tab: linear}}
{\begin{tabular}{c|ccc|ccc|cc}
\hline
\multirow{2}{*}{$s$} & \multicolumn{3}{c|}{Cons. Generation} & \multicolumn{3}{c|}{Least Squares Fit} & \multicolumn{2}{c}{\% Error} \\
                     & $b_0$        & $b_1$      & $\gamma$       & $b_0$         & $b_1$       & $\gamma$        & $b_0$         & $b_1$         \\ \hline
0.01                 & -0.014       & -2.49      & 3.11      & -0.016        & -2.51       & 3.06       & 12\%          & 1\%           \\
0.02                 & -0.020       & -2.49      & 3.10      & -0.023        & -2.51       & 3.02       & 13\%          & 1\%           \\
0.05                 & -0.032       & -2.48      & 3.09      & -0.037        & -2.52       & 2.93       & 14\%          & 1\%           \\
0.1                  & -0.048       & -2.48      & 3.07      & -0.054        & -2.53       & 2.81       & 12\%          & 2\%           \\
0.2                  & -0.071       & -2.46      & 3.04      & -0.081        & -2.55       & 2.63       & 12\%          & 3\%           \\
0.5                  & -0.121       & -2.44      & 2.98      & -0.144        & -2.58       & 2.24       & 16\%          & 5\%           \\ \hline
\end{tabular}}
{}
\end{table}
\section{Asymptotic Optimality of the Fluid Pricing Policy}
\subsection{Proof of Proposition \ref{theo: fluid}} \label{app: fluid}
First note that, under a given pricing and matching policy, if $\E[q^{(c)}_j]=+\infty$ for some $j \in {M}$ or $\E[q^{(s)}_i]=+\infty$ for some $i \in {N}$, then $\mathcal{R}(\BFq,\BFz)=-\infty$ and the theorem is trivially true, as the optimal objective function value \eqref{fluidobj} is greater than or equal to 0 because $\tilde{\BFlambda}=\BFzero_m$ and $\tilde{\BFmu}=\BFzero_n$ is a feasible solution of the fluid optimization problem. So without loss of generality, we assume that all queue lengths have finite expectations.
We will now show the following lemma.
\begin{lemma} \label{claim-necessaryfluidlp}
    For any stationary pricing and matching policy under which the system is stable and $\E[q^{(s)}_i]<\infty, \E[q^{(c)}_j]<\infty$ for all $i \in {N}, j \in {M}$, the expectation of actions $\E[\BFlambda(\BFq)], \E[\BFmu(\BFq)], c\cdot \E[\BFy(\BFq)]$ satisfy
    the constraints in the fluid LP \eqref{ball}--\eqref{consx}.
\end{lemma}
\proof{Proof of Lemma~\ref{claim-necessaryfluidlp}.}
We consider the uniformized DTMC induced by a pricing and matching policy such that the system is stable. By the stability condition, for any period $k$, we have
\begin{align}
    \E[\BFq(k+1)]=\E[\BFq(k)], \nonumber
\end{align}
where the expectation is 
with respect to the stationary distribution of the uniformized DTMC. Let $\BFa(k)$ be the new arrival at period $k$ and let $\BFx(k)$ be the matching decision at period $k$. (Recall that $\BFy(k)\in \mathbb{Z}^{nm}_+$ represents the matching decision for each customer-server pair and $\BFx(k) \in \mathbb{Z}^{n+m}_+$ represents the number of matches for each customer and server type.) We can rewrite the stability condition as
\[
    \E[\BFa(k)]=\E[{\BFx}(k)]. 
\]
Thus, we have $c\E[{\BFx}(k)] = c\E[\BFa(k)]=c\E\left[\E[\BFa(k)\mid \BFq(k)]\right]=\left(\E[\BFlambda(\BFq(k))],\E[\BFmu(\BFq(k))]\right)$, where $\BFlambda(\BFq(k))$ and $\BFmu(\BFq(k))$ are the arrival rates under the given pricing policy in state $\BFq(k)$. By Eq~\eqref{eq:xq}, there exists $y_{ij}(k) \geq 0$ for all $i \in {N}$ and $j \in {M}$ such that,
\begin{align}
    x_j^{(c)}(k)&=\sum_{i=1}^n y_{ij}(k) \leq q_j^{(c)}(k)  \quad \forall j \in {M}, \nonumber\\
    x_i^{(s)}(k)&=\sum_{j=1}^m y_{ij}(k) \leq q_i^{(s)}(k)  \quad \forall i \in {N}, \nonumber\\
    y_{ij}(k)&= 0  \quad \forall (i,j) \notin E. \nonumber
\end{align}
Since the matching policy is stationary, the expectation of the matching decision will not depend on $k$. Taking expectation on both sides with respect to the stationary distribution and defining $\chi_{ij}\overset{\Delta}{=}c\E[y_{ij}(k)]$, we have 
\begin{align}
    \E\left[{\lambda}_j(k)\right] = 
    c\E\left[{x}_j^{(c)}(k)\right] &=\sum_{i=1}^n \chi_{ij}(k) \leq c\E\left[q_j^{(c)}(k)\right] < \infty  \quad \forall j \in {M}, \nonumber\\
    \E\left[{\mu}_i(k)\right] = c\E\left[{x}_i^{(s)}(k)\right] &=\sum_{j=1}^m \chi_{ij}(k) \leq c\E\left[q_i^{(s)}(k)\right] < \infty  \quad \forall i \in {N}, \nonumber\\
    c\E\left[y_{ij}(k)\right] &= \chi_{ij}(k)= 0  \quad \forall (i,j) \notin E. \nonumber
\end{align}
Thus, for any pricing and matching policy under which the system is stable, the expectation of actions satisfy the constraints \eqref{ball}--\eqref{consx}. $\Halmos$ 
\endproof
\proof{Proof of Proposition \ref{theo: fluid}.}
Note that the Uniformized DTMC is aperiodic as we will always have transition from a state back to itself.
By the ergodic theorem for Markov chains, the long run average profit for a given policy is $\E[\mathcal{R}(\BFq,\BFz)]$. 
Also, we have 
\[
\E[\mathcal{R}(\BFq,\BFz)] \leq \E[\inner{F(\BFlambda)}{\BFlambda}-\inner{G(\BFmu)}{\BFmu}] \leq \inner{F(\E[\BFlambda])}{\E[\BFlambda]}-\inner{G(\E[\BFmu])}{\E[\BFmu]},
\]
where the first inequality holds by excluding waiting costs, and the second inequality follows from Jensen's Inequality and Assumption~\ref{ass: concavity}. Thus, the optimal value of the fluid problem \eqref{fluidobj} provides an upper bound for the average profit under any stationary pricing and matching policy. By Lemma~\ref{claim-necessaryfluidlp}, a given policy should also satisfy the constraints \eqref{ball}--\eqref{consx}.
$\Halmos$
\endproof
\subsection{Proof of Theorem \ref{theo: fluidpricingpolicy}} \label{appendix: fluidpricingpolicy}
\proof{Proof of Theorem \ref{theo: fluidpricingpolicy}.}
By Eq \eqref{eq: rl}, the profit-loss $\lr$ is bounded by
\begin{align}
    \lr \leq & \eta\left(\inner{F(\olam)}{(\olam\circ\E[\mathbf{I}^{(c)}(q_{\max}^{\eta})])}-\inner{G(\omu)}{(\omu\circ\E[\mathbf{I}^{(s)}(q_{\max}^{\eta})]}\right)+\inner{\BFs}{\E[\BFq]} \nonumber \\
    \leq & \eta \inner{F(\olam)}{(\olam\circ\E[\mathbf{I}^{(c)}(q_{\max}^{\eta})])} + \inner{\BFs}{\E[\BFq]}. \label{eq:loss-terms}
\end{align}
To bound the first term, we define a function of customer queue lengths $V^{(c)}(\BFq)=\sum_{j=1}^m (q_j^{(c)})^2$. 
Consider the unformized DTMC under the fluid pricing and max-weight matching policy.
Suppose $\BFq$ is the system state in any given period \emph{after} the matching decision $\BFx$ (or equivalently, $\BFy$) has been taken. Similarly, let $\BFq'$ be the state in the next period after the matching decision has been taken. We have
\begin{align}
    &\eta c \left(\E[V^{(c)}(\BFq')\mid \BFq] - V^{(c)}(\BFq)\right) \nonumber \\
    ={}&\sum_{j=1}^m \etal_j \left[ \sum_{j' \in {M}/\{j\}} \left(q_{j'}^{(c)}\right)^2+\left(q_j+\left(1-\sum_{i=1}^m y_{ij}\right)\mathbbm{1}_{\left\{q_j^{(c)}< q_{\max}^{\eta}\right\}}\right)^2\right]\nonumber \\
    &+\sum_{i=1}^m \etam_i \left[ \sum_{j=1}^m \left(q_j^{(c)}-y_{ij}\mathbbm{1}_{\{q_i^{(s)}<q_{\max}^\eta\}}\right)^2\right]-\left[\sum_{i=1}^m \etam_i+\sum_{j=1}^m \etal_j\right] \sum_{j=1}^m \left(q_j^{(c)}\right)^2 \nonumber\\
    {=}{}& \sum_{j=1}^m \etal_j \left(1-\sum_{i=1}^n y_{ij}\right)\left(1+2q_j^{(c)}\right)\mathbbm{1}_{\left\{q_j^{(c)}< q_{\max}^{\eta}\right\}}+\sum_{i=1}^n \etam_i \sum_{j=1}^m y_{ij}\left(1-2q_j^{(c)}\right)\mathbbm{1}_{\{q_i^{(s)}<q_{\max}^\eta\}} \nonumber \\
    {\leq} {}& \sum_{i=1}^n \etam_i + \sum_{j=1}^m \etal_j+2\left[\sum_{j=1}^m\etal_j q_j^{(c)}\left(1-\sum_{i=1}^n y_{ij}\right)\mathbbm{1}_{\left\{q_j^{(c)}< q_{\max}^{\eta}\right\}}-\sum_{i=1}^n\etam_i\sum_{j=1}^m q_j^{(c)} y_{ij}\right] \label{eq: drift_q2_fluid} \\
    {=} {} & \sum_{i=1}^n \etam_i + \sum_{j=1}^m \etal_j+2\left[\sum_{j=1}^m \etal_j q_j^{(c)}\mathbbm{1}_{\left\{q_j^{(c)}< q_{\max}^{\eta}\right\}}\mathbbm{1}_{\left\{\max_{i': (i',j)\in E}q_{i'}^{(s)}=0\right\}}-\sum_{i=1}^n \etam_i \left(\max_{j': (i,j')\in E}q_{j'}^{(c)}\right)\right] \nonumber \\
    {=}{}& 2\sum_{(i,j)\in E}\eta \chi_{ij}^{*}\left[1+q_j^{(c)}\mathbbm{1}_{\left\{q_j^{(c)}< q_{\max}^{\eta}\right\}}\mathbbm{1}_{\left\{\max_{i':(i',j)\in E}q_{i'}^{(s)}=0\right\}}-\max_{j':(i,j') \in E}q_{j'}^{(c)}\right]. \label{eq: drift_const_price}
\end{align}
The first equality holds because the next arrival is a customer of type $j$
with probability ${\etal_j}/{(\eta c)}$, 
and is a server of type $i$ with probability ${\etam_i}/{(\eta c)}$.
The second equality follows as $y_{ij} \in \{0,1\}$. The first inequality follows from the fact that $(1-\sum_{i=1}^n y_{ij}) \leq 1$ for all $j \in {M}$ and $\sum_{j=1}^m y_{ij} \leq 1$ for all $i \in {N}$, because there can be at most one matching in each period under the max-weight policy. The third equality follows from the definition of the max-weight matching policy \eqref{eq:max-weight}. The last equality follows because $(\olam,\omu,\BFchi^{*})$ is the optimal solution to the fluid model and satisfies \eqref{ball}--\eqref{consx}. 
We take expectation with respect to the steady state distribution of $\BFq$.
As $\BFq \leq \BFone q_{\max}^{\eta}$, $V^{(c)}(\BFq)$ is bounded, so $\E[V^{(c)}(\BFq)]$ is finite.
By the assumption of that $\BFq$ follows the steady state distribution, we have $\mathbf{E}[V^{(c)}(\BFq')-V^{(c)}(\BFq)]=0$, so the expectation of the left hand side is 0. 
Therefore,
\begin{align}
    0 \leq {}& \sum_{(i,j) \in E}\chi_{ij}^{*}\left[1+\E\left[q_j^{(c)} \mathbbm{1}_{\left\{q_j^{(c)} < q_{\max}^{\eta}\right\}}\mathbbm{1}_{\left\{\max_{i':(i',j)\in E}q_{i'}^{(s)} =0\right\}}\bigg] -\E\bigg[\max_{j':(i,j') \in E}q_{j'}^{(c)} \right]\right] \nonumber \\
    {\leq}{}& \sum_{(i,j) \in E}\chi_{ij}^{*}\left[1+\E\left[q_j^{(c)} \right]-\E\left[I^{(c)}_j(q_{\max}^{\eta})\right]q_{\max}^{\eta}-\E\left[\max_{j':(i,j') \in E}q_{j'}^{(c)} \right]\right] \nonumber \\
   {\leq} {}& \sum_{(i,j) \in E}\chi_{ij}^{*}\left[1-\left[I^{(c)}_j(q_{\max}^{\eta})\right]q_{\max}^{\eta}\right] \nonumber \\
    = {}& \sum_{j=1}^m \lambda_j^{*} - q_{\max}^{\eta}\sum_{j=1}^m \lambda_j^{*} \left[I^{(c)}_j(q_{\max}^{\eta})\right]. \label{proof1}
\end{align}
The second inequality above holds because $\mathbbm{1}_{\max_{i'\in N(j)}q_{i'}^{(s)} =0} \leq 1$ and $q^{(c)}_j \leq q^{\eta}_{\max}$. (Recall that $I^{(c)}_j(q_{\max}^{\eta}) \overset{\Delta}{=} \mathbbm{1}_{\left\{q_j^{(c)} = q_{\max}^{\eta}\right\}}$.)
The third inequality holds because $\max_{j' \in N(i)}q_{j'}^{(c)}  \geq q_j^{(c)} $ for all $j$ such that $(i,j) \in E$. By substituting $q_{\max}^{\eta}=\gamma \sqrt{\eta}$ in \eqref{proof1} for an arbitrary positive constant $\gamma$, we get
\[
    \gamma\sqrt{\eta}\sum_{j=1}^m \lambda_j^{*} \E\left[I^{(c)}_j(q_{\max}^{\eta})\right] \leq \sum_{j=1}^m \lambda_j^{*},
\]
which implies
\[
\sum_{j=1}^m F_j(\lambda_j^{*})\lambda_j^{*} \E\left[I^{(c)}_j(q_{\max}^{\eta})\right] \leq \frac{1}{\gamma\sqrt{\eta}}\max_{j \in {M}}F_j(\lambda_j^{*})\sum_{j=1}^m \lambda_j^{*}.
\]
Thus, the first term $\inner{F(\olam)}{\olam\circ\E\left[\BFI^{(c)}(q_{\max}^{\eta})\right]}$ of \eqref{eq:loss-terms} is $O(1/\sqrt{\eta})$. 
Now we will bound the second term in \eqref{eq:loss-terms}.
The queue length $\BFq$ under the fluid pricing policy always satisfies $\BFq \leq \BFone_{n+m} q_{\max}^{\eta}$. Thus, it is trivially true that $\inner{\BFs}{\E[\BFq]} \leq q_{\max}^{\eta}\inner{\BFone_{n+m}}{\BFs}=\gamma\sqrt{\eta}\inner{\BFone_{n+m}}{\BFs}$. Thus we can upper bound the profit loss $\lr$ by using \eqref{eq: rl} as follows:
\begin{align}
    \lr &\leq \sqrt{\eta}\left(\max_{j \in {M}}\frac{F_j(\lambda_j^{*})}{\gamma}\sum_{j=1}^m \lambda_j^{*}+ \gamma\inner{\BFone_{n+m}}{\BFs}\right)=O(\sqrt{\eta}). \quad \Halmos \label{eqrl}
\end{align}
\endproof
\subsection{Proof of Proposition \ref{prop: singlelinktwosidedqueue}} \label{appendix: singlelinktwosidedqueue}
\proof{Proof of Proposition \ref{prop: singlelinktwosidedqueue}.}
To prove the lower bound of profit loss, we only need to consider 
single-link queues, that is, $n=m=1$.
In the proof, we omitted the subscript for the type of customers and servers. Under the fluid pricing policy, the steady state distribution of $q=q^{(c)}-q^{(s)}$ is uniform in $[-q^{\eta}_{\max}, q^{\eta}_{\max}]$,
as $q$ behaves like a symmetric simple random walk. Thus, the expected value of the sum of queue length ($q^{(s)}+q^{(c)}$) can be computed in terms of the buffer capacity $q_{\max}^{\eta}$ as follows: 
\begin{align*}
    \E[q^{(s)}(k)+q^{(c)}(k)]&=\E[|q^{(s)}(k)-q^{(c)}(k)|] \\
    &=\frac{q_{\max}^{\eta}(q_{\max}^{\eta}+1)}{2q_{\max}^{\eta}+1}.
\end{align*}
The probabilities of $q^{(s)}=q_{\max}^{\eta}$ and $q^{(c)}=q_{\max}^{\eta}$ are
\begin{align*}
    \Pr[q^{(s)}=q_{\max}^{\eta}]=\frac{1}{2q_{\max}^{\eta}+1}, \\
    \Pr[q^{(c)}=q_{\max}^{\eta}]=\frac{1}{2q_{\max}^{\eta}+1}.
\end{align*}
By Eq \eqref{eq: rl}, the expected profit loss is bounded by
\begin{align}
    \lr &\leq \frac{F(\lambda^{*})\lambda^{*}-G(\mu^{*})\mu^{*}}{2q_{\max}^{\eta}+1}\eta+(s^{(s)}+s^{(c)})\frac{q_{\max}^{\eta}(q_{\max}^{\eta}+1)}{2q_{\max}^{\eta}+1}, \nonumber \\
    \lr &\geq \frac{F(\lambda^{*})\lambda^{*}-G(\mu^{*})\mu^{*}}{2q_{\max}^{\eta}+1}\eta+\min\{s^{(s)},s^{(c)}\}\frac{q_{\max}^{\eta}(q_{\max}^{\eta}+1)}{2q_{\max}^{\eta}+1}. \label{eq: sl-rl}
\end{align}
The inequalities above show that the profit loss is minimized with respect to $\eta$ if $q_{\max}^{\eta}=\gamma \sqrt{\eta}$ for some positive constant $\gamma$. To see this, if $q_{\max}^{\eta}=\gamma \eta^{0.5+\epsilon}$ for some $\epsilon>0$, then due to the second term in \eqref{eq: sl-rl}, $\lr=\Theta(\eta^{0.5+\epsilon})$. On the other hand if $q_{\max}^{\eta}=\gamma \eta^{0.5-\epsilon}$ for some $\epsilon>0$, then due to the first term in \eqref{eq: sl-rl}, $\lr=\Theta(\eta^{0.5+\epsilon})$.
Therefore, by taking $q_{\max}^{\eta}=\gamma \sqrt{\eta}$ in \eqref{eq: sl-rl}, the optimal profit loss is $\lr =\Theta(\sqrt{\eta})$.  $\Halmos$
\endproof
\section{Asymptotic Optimality of the Two-Price Policy}
Throughout the proof, take $\eta$ sufficiently large that
$\eta\lambda_j^*-\theta_j\sigma^\eta\geq0$ and
$\eta\mu_i^*-\phi_i\sigma^\eta\geq0$ for all $i \in N,j \in M$.
\subsection{Proof of Lemma \ref{lemma: posrec}} \label{appendix: posrec}
\proof{Proof of Lemma \ref{lemma: posrec}.}
We start by defining two functions below:
\begin{align}
    V^{(s)}(\BFq)\overset{\Delta}{=}\inner{\BFone_n}{(\BFq^{(s)})^2}, \quad V^{(c)}(\BFq)\overset{\Delta}{=}\inner{\BFone_m}{(\BFq^{(c)})^2}. \nonumber
\end{align}
Consider the uniformized DTMC for the $\eta^{th}$ scaled system under the two-price pricing policy. 
Suppose $\BFq$ is the system state \emph{after} the matching decision $\BFx$ (or equivalently, $\BFy$) has been taken.
We use $\BFq'$ to denote the state in the next period.
The drift of $V^{(c)}(\BFq)$ after one transition  is given by
\begin{align}
    &\eta c\left(\E\left[V^{(c)}(\BFq')\mid \BFq\right] - V^{(c)}(\BFq)\right)\nonumber \\
    ={}&\sum_{j=1}^m \left(\etal_j-\theta_j\sigma^{\eta}\mathbbm{1}_{\left\{q_j^{(c)}>\tau_{\max}^{\eta}\right\}}\right)\left(\left(q_j^{(c)}+1-\sum_{i=1}^n y_{ij}^{}\right)^2-\left(q_j^{(c)}\right)^2\right) \nonumber \\
    &+\sum_{i=1}^n\left(\etam_i-\phi_i\sigma^{\eta}\mathbbm{1}_{\left\{q_i^{(s)}>\tau_{\max}^{\eta}\right\}}\right)\sum_{j=1}^m \left(\left(q_j^{(c)}-y_{ij}^{}\right)^2-\left(q_j^{(c)}\right)^2\right) \nonumber \\
    \overset{(a)}{=}{}& \sum_{j=1}^m \left(\etal_j-\theta_j\sigma^{\eta}\mathbbm{1}_{\left\{q_j^{(c)}>\tau_{\max}^{\eta}\right\}}\right)\left(1-\sum_{i=1}^n y_{ij}^{}\right)\left(1+2q_j^{(c)}\right) \nonumber \\
    & \; +\sum_{i=1}^n \left(\etam_i-\phi_i\sigma^{\eta}\mathbbm{1}_{\left\{q_i^{(s)}>\tau_{\max}^{\eta}\right\}}\right)\sum_{j=1}^m y_{ij}^{}\left(1-2q_j^{(c)}\right) \nonumber \\
    \overset{(b)}{\leq}{}& \eta\inner{\BFone_m}{\olam}+\eta\inner{\BFone_n}{\omu}+2\sum_{j=1}^m \etal_j\left(1-\sum_{i=1}^n y_{ij}^{}\right)q_j^{(c)}-2\sum_{i=1}^n\etam_i\sum_{j=1}^m y_{ij}^{}q_j^{(c)} \nonumber \\
    & -2\sigma^{\eta}\sum_{j=1}^m \theta_j\left(1-\sum_{i=1}^n y_{ij}^{}\right)q_j^{(c)}\mathbbm{1}_{\left\{q_j^{(c)}>\tau_{\max}^{\eta}\right\}}+2\sigma^{\eta}\sum_{i=1}^n \phi_i\sum_{j=1}^m y_{ij}^{}q_j^{(c)}\mathbbm{1}_{\left\{q_i^{(s)}>\tau_{\max}^{\eta}\right\}} \nonumber \\
    \overset{(c)}{=}{}& \eta\inner{\BFone_m}{\olam}+\eta\inner{\BFone_n}{\omu}+2\sum_{j=1}^m \etal_j q_j^{(c)}\mathbbm{1}_{\left\{\max_{i':(i',j)\in E \backslash E_r}q_{i'}^{(s)}=0\right\}}-2\sum_{i=1}^n\etam_i\max_{j':(i,j')\in E \backslash E_r}q_{j'}^{(c)}\nonumber \\
    & -2\sigma^{\eta}\sum_{j=1}^m \theta_j q_j^{(c)}\mathbbm{1}_{\left\{q_j^{(c)} >\tau_{\max}^{\eta}\right\}}\mathbbm{1}_{\left\{\max_{i':(i',j)\in E \backslash E_r}q_{i'}^{(s)}=0\right\}}+2\sigma^{\eta}\sum_{i=1}^n \phi_i\max_{j':(i,j' )\in E \backslash E_r}q_{j'}^{(c)}\mathbbm{1}_{\left\{q_i^{(s)}>\tau_{\max}^{\eta}\right\}} \label{eq: lemma1proofgamma} \\
    \overset{(d)}{=}{}& \eta\inner{\BFone_m}{\olam}+\eta\inner{\BFone_n}{\omu}+2\eta\sum_{(i,j) \in E \backslash E_r}\chi^{*}_{ij} \left( q_j^{(c)}\mathbbm{1}_{\left\{\max_{i':(i',j)\in E \backslash E_r}q_{i'}^{(s)}=0\right\}}-\max_{j':(i,j')\in E \backslash E_r}q_{j'}^{(c)} \right)\nonumber \\ &-2\sigma^{\eta}\sum_{j=1}^m \theta_j q_j^{(c)}\mathbbm{1}_{\left\{q_j^{(c)} >\tau_{\max}^{\eta}\right\}} \label{eq: SSC_MW} \\
    \overset{(e)}{\leq}{}& \eta\inner{\BFone_m}{\olam}+\eta\inner{\BFone_n}{\omu}-2\sigma^{\eta}\sum_{j=1}^m \theta_j q_j^{(c)}\mathbbm{1}_{\left\{q_j^{(c)} >\tau_{\max}^{\eta}\right\}}. \label{eq: driftv2}
\end{align}
Under the modified max-weight matching policy, 
new arrivals are immediately matched to 
compatible counterparts if their queues are nonempty. Thus, $\sum_{i=1}^n y_{ij}^{}$ and $y_{ij}^{}$ (for all $j$) are either 1 or 0. 
Step $(a)$ then follows from the fact that
$(1-\sum_{i=1}^n y_{ij}^{})^2=1-\sum_{i=1}^n y_{ij}^{}$ and $(y_{ij}^{})^2=y_{ij}^{}$.
Step (b) holds because
$\etal_j-\theta_j\sigma^{\eta}\mathbbm{1}_{\left\{q_j^{(c)}>\tau_{\max}^{\eta}\right\}}<\etal_j$ as $\theta_j\sigma^{\eta}>0$ and $1-\sum_{i=1}^n y_{ij}^{} \leq 1$, and because $\etam_i-\phi_i\sigma^{\eta}\mathbbm{1}_{\left\{q_i^{(s)}>\tau_{\max}^{\eta}\right\}}<\etam_i$ as $\phi_i\sigma^{\eta}>0$ and $\sum_{j=1}^m y_{ij}^{}\leq 1$. 
Step $(c)$ follows from the definition of the modified max-weight policy in \eqref{eq:max-weight}.
In Eq~\eqref{eq: lemma1proofgamma}, the event
$q_i^{(s)}>\tau_{\max}^{\eta}$ implies that all compatible types of $i$ in $E \backslash E_r$ have empty queues. Similarly, the event
$q^{(c)}_j > \tau_{\max}^{\eta}$ implies that all compatible types of $j$ in $E \backslash E_r$ have empty queues.
Then, step $(d)$ follows from the definition of the optimal fluid solution $(\BFlambda^*, \BFmu^*, \BFxi^*)$ given by \eqref{eq:fluid_opt}. 
Lastly, step $(e)$ follows because $q_j^{(c)} \leq \max_{j':(i,j') \in E \backslash E_r} q_{j'}^{(c)}$. 
By the same argument, we can bound the drift of $V^{(s)}(\BFq)$ by
\begin{align}
   \eta c\left( \E[V^{(s)}(\BFq')\mid \BFq] - V^{(s)}(\BFq)\right)
    \leq \eta\inner{\BFone_m}{\olam}+\eta\inner{\BFone_n}{\omu}-2\sigma^{\eta}\sum_{i=1}^n \phi_i q_i^{(s)}\mathbbm{1}_{\left\{q_i^{(s)} >\tau_{\max}^{\eta}\right\}}. \label{eq: driftv1}
\end{align}
Define $V(\BFq)\overset{\Delta}{=}V^{(s)}(\BFq)+V^{(c)}(\BFq)$.
Adding \eqref{eq: driftv2} and \eqref{eq: driftv1}, we can bound the drift of $V(\BFq)$ by
\begin{align}
    \eta c\left( \E[V(\BFq')\mid \BFq] - V(\BFq)\right)\leq B-2\sigma^{\eta}\sum_{i=1}^n \phi_i q_i^{(s)}\mathbbm{1}_{\left\{q_i^{(s)} >\tau_{\max}^{\eta}\right\}}-2\sigma^{\eta}\sum_{j=1}^m \theta_j q_j^{(c)}\mathbbm{1}_{\left\{q_j^{(c)} >\tau_{\max}^{\eta}\right\}}, \nonumber
\end{align}
where $B\overset{\Delta}{=}2\eta\inner{\BFone_m}{\olam}+2\eta\inner{\BFone_n}{\omu}>0$ is a positive constant that is independent of $\BFq$. 
Consider the following finite set:
\begin{align}
    \mathcal{B}^{\eta}=\left\{\BFq \in \mathbb{Z}^{m+n}_{+}: q^{(s)}_i \leq \max\left\{\frac{B}{\phi_i\sigma^{\eta}},\tau_{\max}^{\eta}\right\},q^{(c)}_j \leq \max\left\{\frac{B}{\theta_j\sigma^{\eta}},\tau_{\max}^{\eta}\right\} \ \forall i \in {N}, j \in {M}\right\}. \nonumber
\end{align}
Outside the finite set $\mathcal{B}^{\eta}$, the drift of the Lyapunov function $V(\BFq)$ is strictly negative. Specifically, we have
\begin{align}
\eta c\left( \E[V(\BFq')\mid \BFq] - V(\BFq)\right) \leq -B, \quad \forall \BFq \notin \mathcal{B^{\eta}}. \label{eq: foster-lyapunov}
\end{align}
Thus, by the Foster-Lyapunov theorem  \citep{srikantbook}, the uniformized DTMC under the two-price and max-weight policy is positive recurrent for any $\eta$. The first half of the lemma is proved. 
To prove the second half of the lemma, we apply the moment bound theorem \citep{HajekComm} to Eq~\eqref{eq: foster-lyapunov}, which gives that
\begin{align}
    \sigma^{\eta}\E\left[\sum_{i=1}^n \phi_i q_i^{(s)}\mathbbm{1}_{\left\{q_i^{(s)} >\tau_{\max}^{\eta}\right\}}+\sum_{j=1}^m \theta_j q_j^{(c)}\mathbbm{1}_{\left\{q_j^{(c)} >\tau_{\max}^{\eta}\right\}}\right] \leq \frac{B}{2}, \nonumber
\end{align}
where the expectation is taken over the stationary distribution under the two-price and max-weight policy.
By substituting $\mathbbm{1}_{\left\{q^{(s)}_i>\tau_{\max}^{\eta}\right\}}=1-\mathbbm{1}_{\left\{q^{(s)}_i\leq \tau_{\max}^{\eta}\right\}}$ and then using the fact that $\E[q_i^{(s)} \mathbbm{1}_{\left\{q_i^{(s)}\leq \tau_{\max}^{\eta}\right\}}] \leq \tau_{\max}^{\eta}$, we get the desired result.
$\Halmos$
\endproof
\subsection{Proof of Lemma \ref{lemma: firstorderterms}} \label{appendix: firstorderterms}
Before proving Lemma \ref{lemma: firstorderterms}, we first prove the following claim.
\begin{claim} \label{claim}
Consider a pricing and matching policy that has the following form:
  \begin{align}
      {\lambda}_j(\BFq)&=\eta\lambda^{*}_j+\tilde{f}_j(\BFq,\eta) \label{eq: specialformcust} \quad \forall j \in {M}, \\
     {\mu}_i(\BFq)&=\eta\mu^{*}_i+\tilde{g}_i(\BFq,\eta) \quad \forall i \in {N}. \label{eq: specialformserver}
\end{align}    
Suppose the system is positive recurrent under this policy and $\E\left[\langle \mathbf{1}_{m+n}, \mathbf{q}\rangle\right] < \infty$. Then for any optimal solution $(\BFlambda^*, \BFmu^*, \BFchi^*)$ of \eqref{eq:fluid_opt}, we have
    \begin{align*}
        &\sum_{j \in {M}} \left(F_j'(\lambda^{*}_j)\lambda_j^{*}+F_j(\lambda_j^{*})\right)\E[\tilde{f}_j(\BFq,\eta)]-\sum_{i \in {N}} \left(G_i'(\mu_i^{*})\mu_i^{*}+G_i(\mu_i^{*})\right)\E[\tilde{g}_i(\BFq,\eta)]=-\sum_{(i,j) \in E} {\xi}_{ij}{\bar\chi}_{ij}\mathbbm{1}_{{\chi}^{*}_{ij}=0},
    \end{align*}
    where ${\bar\chi}_{ij}$ and ${\xi}_{ij}$ are some nonnegative constants and $(\E[\BFlambda(\BFq)]/\eta, \E[\BFmu(\BFq)]/\eta, {\bar\BFchi}/\eta)$ is a feasible solution of \eqref{eq:fluid_opt}.
\end{claim}
\proof{Proof of Claim~\ref{claim}.}
Denote the objective function \eqref{fluidobj} of the fluid model by $r: \mathbb{R}^{m+n} \to \mathbb{R}$. Denote the left-hand sides of the constraints \eqref{ball}--\eqref{balm} by 
$\BFh : \mathbb{R}^{m+n+mn} \rightarrow \mathbb{R}^{m+n}$. Recall that $(\BFlambda^*, \BFmu^*, \BFchi^*)$ is an optimal solution to the fluid problem. Since the fluid problem is convex and satisfies (relaxed) Slater's condition,
by the KKT conditions, 
there exist Lagrange multipliers $\BFkappa \in  \mathbb{R}^{m+n}$ and $\BFxi \in \mathbb{R}_{+}^{m \times n}$
such that
\begin{align}
    &\nabla r(\BFlambda^{*},\BFmu^{*},{\BFchi}^{*})+\nabla \BFh(\BFlambda^{*},\BFmu^{*},{\BFchi}^{*})^\top \BFkappa+ 
    \sum_{(i,j)\in E}  \xi_{ij}\BFe_{ij}+ \sum_{(i,j) \notin E}  {\xi}_{ij}\BFe_{ij}
    =\BFzero, \label{eq: lemma2optimality}
\end{align}
where 
\begin{align*}
    \nabla r(\olam,\omu,\BFchi^{*})= \left( F'(\olam)\olam+  F(
    \olam),\ - G'(\omu)\omu-G(\omu),\ \BFzero \right) \nonumber
\end{align*}
and $\BFe_{ij} \in \mathbb{R}^{m+n+mn}$ is a vector with the component $m + n + (i,j)$ being 1 and all other components being 0. 
Let $\bar\BFlambda = \E[\BFlambda(\BFq)]$, $\bar\BFmu = \E[\BFmu(\BFq)]$.
Since the system is positive recurrent under the given policy with $\E\left[\langle \mathbf{1}_{m+n}, \mathbf{q}\rangle\right] < \infty$, by Lemma~\ref{claim-necessaryfluidlp}, there exists a vector $\bar\BFchi$ such that
the constraints \eqref{ball}--\eqref{consx} of the fluid model are satisfied by $(\bar\BFlambda, \bar\BFmu, \bar\BFchi)$.
We define a vector $\BFd \in \mathbb{R}^{n+m+mn}$ 
 given by
\begin{align*}
     \BFd&=(\olam,\omu,{\BFchi}^{*})-(\bar\BFlambda,\bar\BFmu,\bar\BFchi)/\eta
     = (- \E[\tilde{\BFf}(\BFq,\eta)]/\eta,\  
     - \E[\tilde{\BFg}(\BFq,\eta)]/\eta,\ {\BFchi}^{*} - \bar\BFchi/\eta).
\end{align*}
Because the function $\BFh$ is linear, it holds that
\begin{equation}
    \inner{\nabla \BFh}{\BFd} = \BFh(\olam,\omu,{\BFchi}^{*}) - \BFh(\bar\BFlambda/\eta,\bar\BFmu/\eta,\bar\BFchi/\eta) = \BFzero.  \label{eq: lemma2equality1}
\end{equation}
By complementary slackness, we have
\begin{equation} \label{eq: lemma2equality2}
    \xi_{ij}({\chi}^{*}_{ij}-\bar{\chi}_{ij}/\eta)
    = - \frac{1}{\eta}\xi_{ij}\bar{\chi}_{ij}\mathbbm{1}_{{\chi}^{*}_{ij}=0},
    \quad \forall (i,j)  \in E.
\end{equation}
Moreover, since $\chi^{*}_{ij}=0$ and $\bar\chi_{ij}=0$ for all $(i,j) \notin E$, we have
\begin{align}
    \xi_{ij}({\chi}^{*}_{ij}-\bar{\chi}_{ij}/\eta)& {=} 0, \quad \forall (i,j)  \notin E. \label{eq: lemma2inequality}
\end{align}
Thus, by taking the inner product on both side by $\BFd$ in \eqref{eq: lemma2optimality} and using \eqref{eq: lemma2equality1}, \eqref{eq: lemma2equality2} and \eqref{eq: lemma2inequality}, we get the desired result. 
$\Halmos$
\endproof
\proof{Proof of Lemma \ref{lemma: firstorderterms}.}
By the definition of redundant edges (Definition~\ref{def:redundant_edges}), if $(i, j) \in E \backslash E_r$, then, there exists an optimal solution $(\BFlambda^{*, ij}, \BFmu^{*, ij}, \BFchi^{*, ij})$ of \eqref{eq:fluid_opt} such that $\chi^{*, ij}_{ij} > 0$. As \eqref{eq:fluid_opt} is a convex optimization problem, $\BFlambda^* = \frac{1}{|E\backslash E_r|}\sum_{ij \in E \backslash E_r}\BFlambda^{*, ij}$, $\BFmu^* = \frac{1}{|E\backslash E_r|}\sum_{ij \in E \backslash E_r}\BFmu^{*, ij}$, and $\BFchi^* = \frac{1}{|E\backslash E_r|}\sum_{ij \in E \backslash E_r}\BFchi^{*, ij}$ is also an optimal solution with $\chi^*_{ij} > 0$ for all $(i, j) \in E \backslash E_r$. We consider this optimal solution and apply Claim~\ref{claim} to prove the lemma.
According to definition of the two-price policy \eqref{eq:two-price_policy-cust}, for all $j\in M$, we have
\begin{align}
    \E[\lambda_j^{(c)}(\BFq)]={}&\E\big[\lambda_j^{(c)}(\BFq)\mid q_j^{(c)}\leq \tau_{\max}^{\eta}\big]\Pr[q_j^{(c)}\leq \tau_{\max}^{\eta}]+\E\big[\lambda_j^{(c)}(\BFq) \mid q_j^{(c)}> \tau_{\max}^{\eta}\big]\Pr[q_j^{(c)}> \tau_{\max}^{\eta}] \nonumber \\
    ={}& \eta\lambda^{*}_j-\theta_j\sigma^{\eta}\Pr[q_j^{(c)}> \tau_{\max}^{\eta}] \quad \forall j \in M. \label{eq: lambda_tp}
\end{align}
Similarly, for all $i\in N$, we have
\begin{align}
    \E[\mu_i^{(s)}(\BFq)]=\eta\mu^{*}_i-\phi_i\sigma^{\eta}\Pr[q_i^{(s)}> \tau_{\max}^{\eta}] \quad \forall i \in N. \label{eq: mu_tp}
\end{align}
By Lemma~\ref{lemma: posrec}, the system is positive recurrent under the two-price policy and modified max-weight matching policy for any $\bth>\BFzero_m$,  $\bph>\BFzero_n$ and $\sigma^{\eta} > 0$ with $\E[\inner{\BFone_{m+n}}{\BFq}] < \infty$. Thus, by Lemma~\ref{claim-necessaryfluidlp}, there exists a vector $\BFchi^{TP}$ such that $\E\left[\BFlambda(\BFq)\right]$ (given by \eqref{eq: lambda_tp}), $\E\left[\BFmu(\BFq)\right]$ (given by \eqref{eq: mu_tp}), and $\BFchi^{TP}$ satisfy the constraints \eqref{ball}--\eqref{consx} of the fluid model. Also, the two-price policy falls under the form given in \eqref{eq: specialformcust} and \eqref{eq: specialformserver} with $\tilde{f}_j(\BFq,\eta)=-\theta_j\sigma^{\eta}\mathbbm{1}_{q_j^{(c)}> \tau_{\max}^{\eta}}$ for all $j \in M$ and $\tilde{g}_i(\BFq,\eta) = -\phi_i\sigma^{\eta}\mathbbm{1}_{q_i^{(s)}> \tau_{\max}^{\eta}}$ for all $i \in N$. Thus, by Claim \ref{claim}, we have 
\begin{align*}
    &\quad \sum_{j \in {M}} \left(F_j'(\lambda_j^{*})\lambda_j^{*}+F_j(\lambda_j^{*})\right)\theta_j \Pr[q^{(c)}_j>\tau_{\max}^{\eta}]-
\sum_{i \in {N}} \left(G_i'(\mu_i^{*})\mu_i^{*}+G_i(\mu_i^{*})\right)\phi_i \Pr[q^{(s)}_i>\tau_{\max}^{\eta}]\\
&= \frac{1}{\sigma^\eta}\sum_{(i,j) \in E} {\xi}_{ij}{\chi}^{TP}_{ij}\mathbbm{1}_{{\chi}^{*}_{ij}=0} \ =0,
\end{align*}
where the last equality follows the definition of modified max-weight policy. In particular, $\chi^{TP}_{ij} = 0$ for all $(i, j) \in E_r$ under the modified max-weight matching policy and $\chi_{ij}^* > 0$ for all $(i, j) \in E \backslash E_r$. This completes the proof.
$\Halmos$
\endproof
\subsection{Proof of Theorem \ref{theorem: twoprice}} \label{appendix: twopriceproof}
\proof{Proof of Theorem \ref{theorem: twoprice}.} We will first calculate the profit loss given by \eqref{eq: rl} as follows:
\begin{align}
    {\lr}    ={}&\fr-({\gamma}^{\eta}-\inner{\BFs}{\E[\BFq]}) \nonumber \\
    ={}& \eta \inner{F(\olam)}{\olam} -\eta \inner{G(\omu)}{\omu}-\eta\left(\sum_{j \in {M}} F_j(\lambda_j^{*})\lambda_j^{*} \Pr[q^{(c)}_j \leq \tau_{\max}^{\eta}]+\right. \nonumber \\
    & \sum_{j \in {M}} F_j\left(\lambda_j^{*}-\frac{\theta_j\sigma^{\eta}}{\eta}\right)\left(\lambda_j^{*}-\frac{\theta_j\sigma^{\eta}}{\eta}\right)\Pr[q^{(c)}_j > \tau_{\max}^{\eta}]-\sum_{i \in {N}} G_i(\mu_i^{*})\mu_i^{*} \Pr[q^{(s)}_i \leq \tau_{\max}^{\eta}] \nonumber \\
    & \left.-\sum_{i \in {N}} G_i\left(\mu_i^{*}-\frac{\phi_i\sigma^{\eta}}{\eta}\right)\left(\mu_i^{*}-\frac{\phi_i\sigma^{\eta}}{\eta}\right)\Pr[q^{(s)}_i > \tau_{\max}^{\eta}]\right)+\inner{\BFs}{\E[\BFq]}. \nonumber \\
    ={}&\eta\sum_{j \in {M}} \left(F_j(\lambda_j^{*})\lambda_j^{*}-F_j\left(\lambda_j^{*}-\frac{\theta_j\sigma^{\eta}}{\eta}\right)\left(\lambda_j^{*}-\frac{\theta_j\sigma^{\eta}}{\eta}\right)\right) \Pr[q^{(c)}_j > \tau_{\max}^{\eta}] \nonumber \\
    &-\eta\sum_{i \in {N}} \left(G_i(\mu_i^{*})\mu_i^{*}-G_i\left(\mu_i^{*}-\frac{\phi_i\sigma^{\eta}}{\eta}\right)\left(\mu_i^{*}-\frac{\phi_i\sigma^{\eta}}{\eta}\right)\right) \Pr[q^{(s)}_i > \tau_{\max}^{\eta}]+\inner{\BFs}{\E[\BFq]}. \label{eq: theo2_profit_loss} 
\end{align}
We apply Taylor's theorem to the terms $F_j(\lambda_j^*-\theta_j\sigma^{\eta}/\eta)$ and $G_i(\mu_i^*-\phi_i\sigma^{\eta}/\eta)$.
For customer type $j$, there exists $\tilde{\lambda}^*_j \in [\lambda^*_j - \theta_j \sigma^\eta/\eta, \lambda^*_j]$ such that, using $\sigma^{\eta}=\eta^{2/3}$, we have
\begin{align*}
    &\quad {\eta F_j(\lambda_j^{*})\lambda_j^{*}-\eta F_j\left(\lambda_j^{*}-\frac{\theta_j\sigma^{\eta}}{\eta}\right)\left(\lambda_j^{*}-\frac{\theta_j\sigma^{\eta}}{\eta}\right)}  \nonumber \\
    &=\eta F_j(\lambda_j^{*})\lambda_j^{*}-\eta\left(F_j(\lambda_j^{*})-\frac{\theta_j\sigma^{\eta}}{\eta}F'_j(\lambda_j^{*})+\frac{1}{2}\left(\frac{\theta_j\sigma^{\eta}}{\eta}\right)^2F''_j(\tilde{\lambda}_j^{*})\right)\left(\lambda_j^{*}-\frac{\theta_j\sigma^{\eta}}{\eta}\right) \nonumber \\
    &= \left(F_j(\lambda_j^{*})+F'_j(\lambda_j^{*})\lambda_j^{*}\right)\theta_j\sigma^{\eta}-\left(\frac{1}{2}F''_j(\tilde{\lambda}_j^{*})\lambda_j^{*}+F'_j(\lambda_j^{*})\right)\frac{(\theta_j\sigma^{\eta})^2}{\eta} + O(\eta^{1/3}) \nonumber \\
    &= \left(F_j(\lambda_j^{*})+F'_j(\lambda_j^{*})\lambda_j^{*}\right)\theta_j\sigma^{\eta}+O(\eta^{1/3}),
\end{align*}
where the last equality holds because $(\sigma^{\eta})^2/\eta = \eta^{1/3}$.
Similar inequality holds for the servers. So, by Eq~\eqref{eq: theo2_profit_loss}, we have
\begin{align*}
    L^\eta{=}{}& \sum_{j \in {M}} \left(F_j'(\lambda_j^{*})\lambda_j^{*}+F_j(\lambda_j^{*})\right)\theta_j\sigma^{\eta} \Pr[q^{(c)}_j>\tau_{\max}^{\eta}] \nonumber \\
&\; -\sum_{i \in {N}} \left(G_i'(\mu_i^{*})\mu_i^{*}+G_i(\mu_i^{*})\right)\phi_i\sigma^{\eta} \Pr[q^{(s)}_i > \tau_{\max}^{\eta}]+O(\eta^{1/3}) + \inner{\BFs}{\E[\BFq]}.
\end{align*}
By Lemma \ref{lemma: firstorderterms}, we can simplify the first two terms in the above equation and get
\begin{align}
    \lr={}& \inner{\BFs}{\E[\BFq]} + O(\eta^{1/3}). \label{eq: rleta} 
\end{align}
Now, using Lemma \ref{lemma: posrec}, we can upper bound the expected queue length as
\begin{align}
    \E \left[\inner{\BFtheta}{\BFq^{(c)}}\right]+\E \left[\inner{\BFphi}{\BFq^{(s)}}\right]  \leq{} &
    \tau_{\max}^{\eta}\left(\sum_{j = 1}^{m} \theta_j  \Pr[q_j^{(c)}\leq \tau_{\max}^{\eta}]+ \sum_{i =1}^{n} \phi_i \Pr[q^{(s)}_i \leq \tau_{\max}^{\eta}]\right) \nonumber\\
    &+ \frac{\eta}{\sigma^{\eta}}\left(\inner{\BFone_n}{\omu} + \inner{\BFone_m}{\olam}\right), \nonumber
\end{align}
which implies that
\begin{align}
    \min_{i \in N, j \in M}\left\{\frac{\phi_i}{s_i^{(s)}},\frac{\theta_j}{s_j^{(c)}}\right\}\E[\inner{\BFs}{\BFq}]\leq{}&
    \tau_{\max}^{\eta}\left(\sum_{j = 1}^{m} \theta_j  + \sum_{i =1}^{n} \phi_i \right)+ \frac{\eta}{\sigma^{\eta}}\left(\inner{\BFone_n}{\omu} + \inner{\BFone_m}{\olam}\right). \nonumber
\end{align}
By substituting $\tau_{\max}^{\eta}\leq\eta^{1/3}$, $\sigma^{\eta}=\eta^{2/3}$ and noting that
$\min_{i \in N, j \in M}\{{\phi_i}/{s_i^{(s)}},{\theta_j}/{s_j^{(c)}}\}>0$, we have $\E[\inner{\BFs}{\BFq}] {=}O(\eta^{1/3})$.
By substituting this term in \eqref{eq: rleta} and substituting $\sigma^\eta=\eta^{2/3}$, we have the desired result.
$\Halmos$
\endproof
\section{Lower Bounds}
\subsection{Proof of Theorem \ref{theorem: lowerboundmultiplelink}} \label{appendix: lowerboundmultiplelink}
First, we will present a lemma that provides a lower bound on the expected value of the sum of the queue length $\E\left[\inner{\BFone_{n+m}}{\BFq}\right]$.
\begin{lemma} \label{lemma: boundedf2g2}
If Condition \ref{ass: generalpricinggraph}\ref{condition: boundedfg} is satisfied, then
\begin{align}
   \E[\inner{\BFone_{n+m}}{\BFq}]&\geq \frac{\eta^{1-\beta}\inner{\BFone_m}{\olam}}{2\left(\inner{\BFone_n}{{\BFPsi}}+\inner{\BFone_m}{{\BFGamma}}\right)}-0.5.
\end{align}
In addition, if Condition \ref{ass: generalpricinggraph}\ref{condition: aplusb} and \ref{ass: generalpricinggraph}\ref{condition: posrecsystem} are satisfied, then there exist $\epsilon>0$ and $\mathcal{M}>0$ such that
\begin{align}
   \sum_{j=1}^m \E\left[f_j^2\left(\frac{\BFq}{\eta^{\alpha}}\right)\right]+ \sum_{i=1}^n \E\left[g_i^2\left(\frac{\BFq}{\eta^{\alpha}}\right)\right] &\geq \epsilon  \quad \forall \eta > \mathcal{M}. \label{eq: lowerboundf2}
\end{align}
\end{lemma}
\proof{Proof of Lemma \ref{lemma: boundedf2g2}.}
 We define the ``imbalance'' process of the two-sided queueing system as
  \begin{align}
     z(k)=\inner{\BFone_m}{\BFq^{(c)}(k)}-\inner{\BFone_n}{\BFq^{(s)}(k)}. \label{eq: imbalance}
\end{align}
 Now, we will define a new DTMC $\{\tilde{z}(k): k \in \mathbb{Z}_{+}\}$ and couple it with the uniformized DTMC $\{\BFq(k),z(k): k \in \mathbb{Z}_{+}\}$ with uniformization constant $c$ such that $|z(k)| \geq \tilde{z}(k)$ for all $k \in \mathbb{Z}_{+}$ for all $\eta \geq 1$ if $z(0)=\tilde{z}(0)$.
\begin{figure}
    \TABLE {Coupled Birth and Death Process  \label{fig:my_label}}
    {\begin{tabular}{c}
        \begin{tikzpicture}[->, >=stealth', auto, semithick, node distance=2.5cm]
\tikzstyle{every state}=[fill=white,draw=black,thick,text=black,scale=1]
\node[state]    (0)                     {$0$};
\node[state]    (1)[right of=0]   {$1$};
\node[state]    (2)[right of=1]   {$2$};
\node[state]    (nminusone)[right of=2]   {\tiny{$n-1$}};
\node[state]    (n)[ right of=nminusone]   {$n$};
\node[state,white]    (nplusone)[ right of=n] {};
\path
(0) edge[loop left]     node{$1-\tilde{P}_{12}^{\eta}$}         (0)
(0) edge[bend right,below]     node{$\tilde{P}_{12}^{\eta}$}         (1)
(1) edge[bend right,below]     node{$\tilde{P}_{12}^{\eta}$}         (2)
(2) edge[bend right,below,dashed]     node{$\tilde{P}_{12}^{\eta}$}         (nminusone)
(nminusone) edge[bend right,below]     node{$\tilde{P}_{12}^{\eta}$}         (n)
(n) edge[bend right,below]     node{$\tilde{P}_{12}^{\eta}$}         (nplusone)
(nplusone) edge[bend right,above]     node{$\tilde{P}_{10}^{\eta}$}         (n)
(n) edge[bend right,above]     node{$\tilde{P}_{10}^{\eta}$}         (nminusone)
(nminusone) edge[bend right,above,dashed]     node{$\tilde{P}_{10}^{\eta}$}         (2)
(2) edge[bend right,above]     node{$\tilde{P}_{10}^{\eta}$}         (1)
(1) edge[bend right,above]     node{$\tilde{P}_{10}^{\eta}$}         (0)
(1) edge[loop,above]     node{$\tilde{P}_{11}^{\eta}$}         (1)
(2) edge[loop,above]     node{$\tilde{P}_{11}^{\eta}$}         (2)
(nminusone) edge[loop,above]     node{$\tilde{P}_{11}^{\eta}$}         (nminusome)
(n) edge[loop,above]     node{$\tilde{P}_{11}^{\eta}$}         (n);
\end{tikzpicture}
    \end{tabular}}
    {}
\end{figure}
The state space of $\{\tilde{z}(k): k \in \mathbb{Z}_{+}\}$ is $\mathbb{Z}_{+}$ and the transition matrix is given by
\begin{align}
    \tilde{P}_{ij}^{\eta}=\begin{cases}
    \tilde{P}_{10}^{\eta}\overset{\Delta}{=}\left(\inner{\BFone_n}{\BFmu^{*}}+\left(\inner{\BFone_n}{{\BFPsi}}+\inner{\BFone_m}{{\BFGamma}}\right)\eta^{\beta-1}\right)/c
     &\textit{if } j=i-1,\ \forall j \in \mathbb{Z}_+,\ \forall i >0, \\
     \tilde{P}_{12}^{\eta}\overset{\Delta}{=}\left(\inner{\BFone_m}{\BFlambda^{*}}-\left(\inner{\BFone_n}{{\BFPsi}}+\inner{\BFone_m}{{\BFGamma}}\right)\eta^{\beta-1}\right)/c &\textit{if } j=i+1,\ \forall j \in \mathbb{Z}_+,\  \forall i >0, \\
    \tilde{P}_{11}^{\eta}\overset{\Delta}{=}1- \left(\inner{\BFone_n}{\BFmu^{*}}+\inner{\BFone_m}{\BFlambda^{*}}\right)/c& \textit{if } j=i, \ \forall j \in \mathbb{Z}_+,\ \forall i >0, \\
     \left(\inner{\BFone_m}{\BFlambda^{*}}-\left(\inner{\BFone_n}{{\BFPsi}}+\inner{\BFone_m}{{\BFGamma}}\right)\eta^{\beta-1}\right)/c&\textit{if } j=1, \ i=0, \\
    1- \left(\inner{\BFone_m}{\BFlambda^{*}}-\left(\inner{\BFone_n}{{\BFPsi}}+\inner{\BFone_m}{{\BFGamma}}\right)\eta^{\beta-1}\right)/c& \textit{if } j=0, \ i=0, \\
    0 &\textit{otherwise},
    \end{cases}  \nonumber
\end{align}
where $c$ is the uniformization constant given by Definition~\ref{ass: upperboundedrates}. 
See Figure~\ref{fig:my_label} for an illustration of the new DTMC.
Note that, under Condition \ref{ass: generalpricinggraph}\ref{condition: boundedfg}, we have $\tilde{P}_{ij}^{\eta} \geq \Pr[|z(k+1)|=j,\ | |z(k)|=i,\ \BFq(k)=\bar{\BFq}]$ for all $\bar{\BFq} \in S$ if $j<i$ and  $\tilde{P}_{ij}^{\eta} \leq \Pr[|z(k+1)|=j,\ | |z(k)|=i,\ \BFq(k)=\bar{\BFq}]$ for all $\bar{\BFq} \in S$ if $j>i$ for all $\eta>0$. Thus, we can couple these system using a common source of randomness, such that $\tilde{z}(k) \leq |z(k)|$ for all $k \in \mathbb{Z}_{+}$ on each sample path. Hence, we have $\Pr[\tilde{z}(k)\leq \mathcal{K}]\geq \Pr[z(k) \leq \mathcal{K}]$ for all $k \geq 1$ and $\mathcal{K} \in \mathbb{R}$ 
Thus, in the limit as $k \rightarrow \infty$, we have $\Pr[\tilde{z}(\infty)\leq \mathcal{K}]\geq \Pr[z(\infty) \leq \mathcal{K}]$ where, $\tilde{z}(\infty)$ and $z(\infty)$ denotes random variables following the stationary distributions of $\{\tilde{z}(k): k \geq 1\}$ and $\{z(k): k \geq 1\}$, respectively. So, we have
\begin{align}
   \E\left[\tilde{z}(\infty)\right] &\leq \E\left[|z(\infty)|\right]. \label{eq: lemmaexpectationbelow}
\end{align}
The stationary distribution of  $\{\tilde{z}(k): k \geq 1\}$ is given by
\begin{align*}
    \pi_i&=
    \pi_0\left(\frac{\tilde{P}_{12}^{\eta}}{\tilde{P}_{10}^{\eta}}\right)^i \ \ \forall i>0, \qquad
    \sum_{i=-\infty}^{\infty} \pi_i =1.
\end{align*}
Solving for $\pi_i$, we get
\begin{align*}
\pi_i=\left(1-\frac{\tilde{P}_{12}^{\eta}}{\tilde{P}_{10}^{\eta}}\right)\left(\frac{\tilde{P}_{12}^{\eta}}{\tilde{P}_{10}^{\eta}}\right)^i \quad \forall i \in \mathbb{Z}_{+}.
\end{align*}
Thus, we have
\begin{align}
    \E[\tilde{z}(\infty)]&=\frac{\tilde{P}_{12}^{\eta}}{\tilde{P}_{10}^{\eta}-\tilde{P}_{12}^{\eta}} {=}\frac{\eta^{1-\beta}\inner{\BFone_m}{\olam}}{2\left(\inner{\BFone_n}{{\BFPsi}}+\inner{\BFone_m}{{\BFGamma}}\right)}-0.5. \label{eq: lemmaexpectationtilde}
\end{align}
where the second equality uses the definition of
$\tilde{P}^{\eta}_{ij}$ above
and the fact that $\inner{\BFone_n}{\BFmu^{*}}=\inner{\BFone_m}{\BFlambda^{*}}$. Finally, note that $|z| \leq \inner{\BFone_{n+m}}{\BFq}$. Thus, we have
\begin{align*}
   \E[\inner{\BFone_{n+m}}{\BFq}] \ \geq\ \E[|z(\infty)|] \ {\geq}\ \E[\tilde{z}(\infty)]&\ {=}\  \frac{\eta^{1-\beta}\inner{\BFone_m}{\olam}}{2\left(\inner{\BFone_n}{{\BFPsi}}+\inner{\BFone_m}{{\BFGamma}}\right)}-0.5,
\end{align*}
where the second inequality uses \eqref{eq: lemmaexpectationbelow} and the last equality uses
\eqref{eq: lemmaexpectationtilde}.
This completes the first part of the proof. 
Next, as $\tilde{z}(\infty)\leq |z(\infty)|$ almost surely, for any constant $\mathcal{K}$, we have
\begin{align}
\Pr[|z(\infty)|>\mathcal{K}] \geq \Pr[\tilde{z}(\infty)>\mathcal{K}]=\sum_{i=\mathcal{K}+1}^{\infty} \pi_i=\left(\frac{\tilde{P}_{12}^{\eta}}{\tilde{P}_{10}^{\eta}}\right)^{\mathcal{K}+1}=\left(\frac{\inner{\BFone_m}{\BFlambda^{*}}-\left(\inner{\BFone_n}{{\BFPsi}}+\inner{\BFone_m}{{\BFGamma}}\right)\eta^{\beta-1}}{\inner{\BFone_n}{\BFmu^{*}}+\left(\inner{\BFone_n}{{\BFPsi}}+\inner{\BFone_m}{{\BFGamma}}\right)\eta^{\beta-1}}\right)^{\mathcal{K}+1}\hspace{-15pt}. \label{eq: lb_couping_imbalance}
\end{align}
Let $\mathcal{K}=(n+m)\kappa\eta^{\alpha}$. When $\alpha=0$, we have
\begin{align*}
    \lim_{\eta \rightarrow \infty}&\left(\frac{\inner{\BFone_m}{\BFlambda^{*}}-\left(\inner{\BFone_n}{{\BFPsi}}+\inner{\BFone_m}{{\BFGamma}}\right)\eta^{\beta-1}}{\inner{\BFone_n}{\BFmu^{*}}+\left(\inner{\BFone_n}{{\BFPsi}}+\inner{\BFone_m}{{\BFGamma}}\right)\eta^{\beta-1}}\right)^{(n+m)\kappa+1}=1^{(n+m)\kappa+1}=1.
\end{align*}
For $\alpha>0$, we have
\begin{align*}
    & \quad \lim_{\eta \rightarrow \infty} \left(\frac{\inner{\BFone_m}{\BFlambda^{*}}-\left(\inner{\BFone_n}{{\BFPsi}}+\inner{\BFone_m}{{\BFGamma}}\right)\eta^{\beta-1}}{\inner{\BFone_n}{\BFmu^{*}}+\left(\inner{\BFone_n}{{\BFPsi}}+\inner{\BFone_m}{{\BFGamma}}\right)\eta^{\beta-1}}\right)^{(n+m)\kappa\eta^{\alpha}}\\
    &=
    \lim_{\eta\rightarrow \infty} \left[\left(1-\frac{2\left(\inner{\BFone_n}{{\BFPsi}}+\inner{\BFone_m}{{\BFGamma}}\right)}{\inner{\BFone_m}{\BFmu^{*}}\eta^{1-\beta}+\left(\inner{\BFone_n}{{\BFPsi}}+\inner{\BFone_m}{{\BFGamma}}\right)}\right)^{\frac{\inner{\BFone_n}{\BFmu^*}\eta^{1-\beta}}{3\left(\inner{\BFone_n}{{\BFPsi}}+\inner{\BFone_m}{{\BFGamma}}\right)}}\right]^{\frac{3K(n+m)\left(\inner{\BFone_n}{{\BFPsi}}+\inner{\BFone_m}{{\BFGamma}}\right)}{\inner{\BFone_n}{\BFmu^*}}\eta^{\alpha+\beta-1}} \\
    & {\geq} \lim_{\eta \rightarrow \infty}\left(1-\frac{2}{3+\frac{3\left(\inner{\BFone_n}{{\BFPsi}}+\inner{\BFone_m}{{\BFGamma}}\right)}{\inner{\BFone_n}{\BFmu^*}}\eta^{\beta-1}}\right)^{\frac{3K(n+m)\left(\inner{\BFone_n}{{\BFPsi}}+\inner{\BFone_m}{{\BFGamma}}\right)}{\inner{\BFone_n}{\BFmu^*}}\eta^{\alpha+\beta-1}}\\
    &=\begin{cases}
    (1/3)^0=1 &\textit{if } \alpha+\beta<1, \\
    b\overset{\Delta}{=}\left(\frac{1}{3}\right)^{\frac{3K(n+m)\left(\inner{\BFone_n}{{\BFPsi}}+\inner{\BFone_m}{{\BFGamma}}\right)}{\inner{\BFone_n}{\BFmu^*}}} &\textit{if } \alpha+\beta=1,
    \end{cases}
\end{align*}
where the first inequality follows from the Bernoulli's inequality, which says that $(1+x)^r \geq 1+rx$ if $x \geq -2$. 
Thus, we have
\begin{align*}
    & \lim_{\eta \rightarrow \infty} \left(\frac{\inner{\BFone_m}{\BFlambda^{*}}-\left(\inner{\BFone_n}{{\BFPsi}}+\inner{\BFone_m}{{\BFGamma}}\right)\eta^{\beta-1}}{\inner{\BFone_n}{\BFmu^{*}}+\left(\inner{\BFone_n}{{\BFPsi}}+\inner{\BFone_m}{{\BFGamma}}\right)\eta^{\beta-1}}\right)^{(n+m)\kappa\eta^{\alpha}+1} \\
    &{=}\lim_{\eta \rightarrow \infty}\left(\frac{\inner{\BFone_m}{\BFlambda^{*}}-\left(\inner{\BFone_n}{{\BFPsi}}+\inner{\BFone_m}{{\BFGamma}}\right)\eta^{\beta-1}}{\inner{\BFone_n}{\BFmu^{*}}+\left(\inner{\BFone_n}{{\BFPsi}}+\inner{\BFone_m}{{\BFGamma}}\right)\eta^{\beta-1}}\right)^{(n+m)\kappa\eta^{\alpha}}\hspace{-15pt}\cdot \lim_{\eta \rightarrow \infty}\left(\frac{\inner{\BFone_m}{\BFlambda^{*}}-\left(\inner{\BFone_n}{{\BFPsi}}+\inner{\BFone_m}{{\BFGamma}}\right)\eta^{\beta-1}}{\inner{\BFone_n}{\BFmu^{*}}+\left(\inner{\BFone_n}{{\BFPsi}}+\inner{\BFone_m}{{\BFGamma}}\right)\eta^{\beta-1}}\right) \\
    &=\begin{cases}
    1 &\textit{if } \alpha+\beta<1, \\
    b &\textit{if } \alpha+\beta=1,
    \end{cases}
\end{align*}
where the first equality holds because the limit of the product of two sequences is equal to the product of the limits of the two sequences. Thus, when Condition \ref{ass: generalpricinggraph}\ref{condition: aplusb} is satisfied, namely $\alpha+\beta \leq 1$, for any given $b\delta^2>\epsilon>0$, there exists a constant $\mathcal{M}$ such that for all $\eta>\mathcal{M}$, it holds that
    \begin{align*}
        \Pr\left[|z(\infty)|>(n+m)\kappa\eta^{\alpha}\right] \geq \frac{\epsilon}{\delta^2}.
    \end{align*}
By the definition of $z$ in Eq~\eqref{eq: lemmaexpectationtilde}, the event
\begin{align*}
    \{|z|>(n+m)\kappa\eta^{\alpha}\} \subset \left\{\|\BFq\|_{\infty}>\kappa\eta^{\alpha}\right\}.
\end{align*}
Also, by Condition \ref{ass: generalpricinggraph} \ref{condition: aplusb}, we have 
\begin{align*}
    \sum_{j=1}^m f_j^2\left(\frac{\BFq}{\eta^{\alpha}}\right)+ \sum_{i=1}^n g_i^2\left(\frac{\BFq}{\eta^{\alpha}}\right) \geq \delta^2 \quad  \forall \BFq : \|\BFq\|_{\infty}>\kappa\eta^{\alpha}.
\end{align*}
Thus for all $\eta>\mathcal{M}$, we have
\begin{align*}
   \sum_{j=1}^m \E\left[f_j^2\left(\frac{\BFq}{\eta^{\alpha}}\right)\right]+ \sum_{i=1}^n \E\left[g_i^2\left(\frac{\BFq}{\eta^{\alpha}}\right)\right] &\geq \delta^2\Pr\left[|z(\infty)|>(n+m)\kappa\eta^{\alpha}\right] \geq\delta^2\frac{\epsilon}{\delta^2}
    =\epsilon. \Halmos
\end{align*}
\endproof
\begin{claim} \label{claim: uniformconvergence}
It holds that
\begin{align*}
    A_j \overset{\Delta}{=} -F'_j(\lambda^{*}_j)-\frac{1}{2}\lambda^{*}_j F_j''(\lambda^{*}_j)&>0 , \quad \forall j \in {M},\nonumber \\
 B_i \overset{\Delta}{=}\ G'_i(\mu^{*}_i)+\frac{1}{2}\mu^{*}_iG''_i(\mu^{*}_i)&>0 \quad \forall i 
 \in {N}.\nonumber
\end{align*}
There exists a constant $\eta'>0$ such that for all $\eta>\eta'$ and $\tilde{c}\overset{\Delta}{=}\min_{j \in {M}, i \in {N}}\left\{{A_j}/{\lambda_j^{*}},\ {B_i}/{\mu_i^*}\right\}$, it holds that
    \begin{align*}
         \sup_{\BFq \in S} |F''_j(\tilde{\lambda}_j(\BFq))-F''_j(\lambda^{*}_j)| < \tilde{c}  \quad \forall j \in {M}, \\
         \sup_{\BFq \in S} |G''_i(\tilde{\mu}_i(\BFq))-G''_i(\mu^{*}_i)| < \tilde{c}  \quad \forall i \in {N}.
    \end{align*}
\end{claim}
\proof{Proof of Claim \ref{claim: uniformconvergence}.}
As $-F_j(\lambda_j)\lambda_j$ for all $j \in {M}$ and $G_i(\mu_i)\mu_i$ for all $i \in {N}$ are strictly convex by Assumption \ref{ass: concavity}, we have
\begin{align}
    \frac{d^2}{d\lambda_j^2}\left(-F_j(\lambda_j)\lambda_j\right)=2 A_j&> 0 \;\; \forall j \in {M}, \quad
     \frac{d^2}{d\mu_i^2}\left(G_i(\mu_i)\mu_i\right)=2 B_i> 0 \;\; \forall i 
     \in {N}.\nonumber
\end{align}
Since $F''_j(.)$ is continuous by assumption, given $\tilde{c}=\min_{i,j}\{\frac{A_j}{\lambda^*_j},\frac{B_i}{\mu^*_i}\}>0$, there exists $\delta'>0$ such that for any $l \in [\lambda^{*}_j-\delta',\lambda^{*}_j+\delta']$ we have
\begin{align*}
    |F''_j(l)-F''_j(\lambda^{*}_j)| < \tilde{c}.
\end{align*}
As $\beta<1$, we have $\eta^{\beta-1}\rightarrow 0$ as $\eta\to\infty$. Consider $\eta'$ such that for all $\eta>\eta'$ we have $|\eta^{\beta-1}|<\delta'/\max_{j \in {M}, i\in {N}}\{{\Gamma}_j,{\Psi}_i\}$, which implies that $\tilde{\lambda}_j(\BFq) \in [\lambda^{*}_j-\delta',\lambda^{*}_j+\delta']$ for all $\BFq$. Thus,
\begin{align*}
    \sup_{\BFq \in S} |F''_j(\tilde{\lambda}_j(\BFq))-F''_j(\lambda^{*}_j)| < \tilde{c}.
\end{align*}
The proof for $G''_i(\mu)$ is similar.
$\Halmos$
\endproof
\proof{Proof of Theorem \ref{theorem: lowerboundmultiplelink}.}
Without loss of generality, we can assume that 
\begin{align}
    \E[\inner{\BFone_{n+m}}{\BFq}]<\infty; \label{eq: boundedqueueexpectedtheorem}
\end{align}
otherwise, the holding cost will be infinity and the lower bound holds trivially. We calculate the profit loss and use Lemma \ref{lemma: boundedf2g2} to lower bound the expected queue length. We have
\begin{align}
\lr={}&\eta\inner{F(\olam)}{\olam}-\eta\inner{G(\omu)}{\omu}-\E[\inner{F^{\eta}(\BFlambda^{\eta})}{\BFlambda^{\eta}}-\inner{G^{\eta}(\BFmu^{\eta})}{\BFmu^{\eta}}-s\inner{\BFone_{n+m}}{\BFq}] \nonumber \\
    {=}{}&\eta\inner{F(\olam)}{\olam}-\eta\inner{G(\omu)}{\omu}-\E\left[\sum_{j=1}^m \left(\lambda^{*}_j+f_j\left(\frac{\BFq}{\eta^{\alpha}}\right)\eta^{\beta-1}\right)F_j\left(\lambda^{*}_j+f_j\left(\frac{\BFq}{\eta^{\alpha}}\right)\eta^{\beta-1}\right)\right.\nonumber\\
    &\left.-\sum_{i=1}^n \left(\mu^{*}_i+g_i\left(\frac{\BFq}{\eta^{\alpha}}\right)\eta^{\beta-1}\right)G_i\left(\mu^{*}_i+g_i\left(\frac{\BFq}{\eta^{\alpha}}\right)\eta^{\beta-1}\right)-\inner{\BFs}{\BFq}\right], \label{eq: lowerboundtheoremtaylor}
    \end{align}
   where the second equality follows from the definition of the asymptotic regime (Definition \ref{def: asymptoticregime}) and the assumptions on the pricing policy Eq~\eqref{eq: pricinggraphgeneralcust}.
   Now, we will use Taylor's theorem to expand each function $F_j$ and $G_i$ individually. For the term involving $F_j$, we have
    \begin{align*}
    &\E\left[\left(\eta\lambda_j^{*}+f_j\left(\frac{\BFq}{\eta^{\alpha}}\right)\eta^{\beta}\right)F_j\left(\lambda_j^{*}+f_j\left(\frac{\BFq}{\eta^{\alpha}}\right)\frac{\eta^{\beta}}{\eta}\right)\right] \\
    {=}{}&\E\left[\left(\eta\lambda_j^{*}+f_j\left(\frac{\BFq}{\eta^{\alpha}}\right)\eta^{\beta}\right)\left(F_j\left(\lambda_j^{*}\right)+f_j\left(\frac{\BFq}{\eta^{\alpha}}\right)\eta^{\beta-1}F_j'(\lambda_j^{*})+f_j^2\left(\frac{\BFq}{\eta^{\alpha}}\right)\frac{\eta^{2\beta-2}}{2}F''_j\left(\tilde{\lambda}_j(\BFq)\right)\right)\right] \\
    ={}&\eta\lambda_j^{*}F_j(\lambda_j^{*})+\left(F'_j(\lambda_j^{*})\lambda_j^{*}+F_j(\lambda_j^{*})\right)\E\left[f_j\left(\frac{\BFq}{\eta^{\alpha}}\right)\right]\eta^{\beta}+F'_j(\lambda_j^{*})\E\left[f^2_j\left(\frac{\BFq}{\eta^{\alpha}}\right)\right]\eta^{2\beta-1}\nonumber\\
    &+\lambda_j^{*}\E\left[f^2_j\left(\frac{\BFq}{\eta^{\alpha}}\right)F''_j\left(\tilde{\lambda}_j(\BFq)\right)\right]\frac{\eta^{2\beta-1}}{2}+\E\left[f^3_j\left(\frac{\BFq}{\eta^{\alpha}}\right)F''_j\left(\tilde{\lambda}_j(\BFq)\right)\right]\frac{\eta^{3\beta-2}}{2},
\end{align*}
where the first term follows from Taylor's theorem and $\tilde{\lambda}_j(\BFq) \in [\lambda^{*}_j,\lambda^{*}_j+f_j\left(\frac{\BFq}{\eta^{\alpha}}\right)\eta^{\beta-1}]$. Similarly, we can use the Taylor's theorem for all $G_i$ in \eqref{eq: lowerboundtheoremtaylor} by noting that 
$\tilde{\mu}_i(\BFq)\in[\mu^{*}_i,\mu^{*}_i+g_i\left(\frac{\BFq}{\eta^{\alpha}}\right)\eta^{\beta-1}]$ for all $i \in {N}$. Therefore, we have
    \begin{align}
    \lr={}&\underbrace{\sum_{i=1}^n\E\left[g_i\left(\frac{\BFq}{\eta^{\alpha}}\right) \eta^{\beta}\right]\big(\mu^{*}_i G_i'(\mu^{*}_i)+G_i(\mu^{*}_i))\big)-\sum_{j=1}^m\E\left[f_j\left(\frac{\BFq}{\eta^{\alpha}}\right)\eta^{\beta}\right]\big(\lambda^{*}_j F_j'(\lambda^{*}_j)+F_j(\lambda^{*}_j))\big)}_{\mathcal{A}_1}\nonumber \\
    &\underbrace{-\eta^{2\beta-1}\left(\sum_{j=1}^mF_j'(\lambda^{*}_j)\E\left[f^2_j\left(\frac{\BFq}{\eta^{\alpha}}\right)\right]+\sum_{j=1}^m\frac{\lambda^{*}_j}{2}\E\left[F''_j(\tilde{\lambda}_j(\BFq))f^2_j\left(\frac{\BFq}{\eta^{\alpha}}\right)\right]\right)}_{\mathcal{A}_2^{(c)}} \nonumber \\
    &\underbrace{+\eta^{2\beta-1}\left(\sum_{i=1}^nG'_i(\mu^{*}_i)\E\left[g^2_i\left(\frac{\BFq}{\eta^{\alpha}}\right)\right]+\sum_{i=1}^n\frac{\mu^{*}_i}{2}\E\left[G''_i(\tilde{\mu}(\BFq))g^2_i\left(\frac{\BFq}{\eta^{\alpha}}\right)\right]\right)}_{\mathcal{A}_2^{(s)}}\nonumber\\
    &\underbrace{-\frac{\eta^{3\beta-2}}{2}\left(\sum_{j=1}^m\E\left[F''_j(\tilde{\lambda}_j(\BFq))f^3_j\left(\frac{\BFq}{\eta^{\alpha}}\right)\right]-\sum_{i=1}^n\E\left[G''_i(\tilde{\mu}_i(\BFq))g_i^3\left(\frac{\BFq}{\eta^{\alpha}}\right)\right]\right)}_{\mathcal{A}_3}+\underbrace{\E[\inner{\BFs}{\BFq}]}_{\mathcal{A}_4}\label{eq: lowerboundequationtheorem}
    \end{align}

We now bound each of the terms in the above equation individually. Firstly, we have $\mathcal{A}_1\geq 0$ by Claim \ref{claim} and \eqref{eq: boundedqueueexpectedtheorem}. Next, the terms $\mathcal{A}_2^{(s)}$ and $\mathcal{A}_2^{(c)}$ can be simplified using Claim~\ref{claim: uniformconvergence}: for any $\eta > \eta'$, we have
    \begin{align}
      \mathcal{A}_2^{(c)}+\mathcal{A}_2^{(s)} &\geq \eta^{2\beta-1}\left(\sum_{j=1}^m\left(A_j-\frac{\lambda^*_j\tilde{c}}{2}\right)\E\left[f_j^2\left(\frac{\BFq}{\eta^{\alpha}}\right)\right]+\sum_{i=1}^n\left(B_i-\frac{\mu^*_i\tilde{c}}{2}\right)\E\left[g_i^2\left(\frac{\BFq}{\eta^{\alpha}}\right)\right]\right)\nonumber \\
      &\geq \eta^{2\beta-1}\left(\sum_{j=1}^m \frac{A_j}{2}\E\left[f_j^2\left(\frac{\BFq}{\eta^{\alpha}}\right)\right]+\sum_{i=1}^n \frac{B_i}{2}\E\left[g_i^2\left(\frac{\BFq}{\eta^{\alpha}}\right)\right]\right), \nonumber
    \end{align}
where the second inequality uses the definition of $\tilde{c}$.
For the term $\mathcal{A}_3$, we have
    \begin{align*}
        \mathcal{A}_3&=-\frac{\eta^{3\beta-2}}{2}\left(\sum_{j=1}^m\E\left[F''_j(\tilde{\lambda}_j(\BFq))f^3_j\left(\frac{\BFq}{\eta^{\alpha}}\right)\right]-\sum_{i=1}^n\E\left[G''_i(\tilde{\mu}_i(\BFq))g_i^3\left(\frac{\BFq}{\eta^{\alpha}}\right)\right]\right) \\
        &\geq -\eta^{3\beta-2}\frac{\inner{\BFone_{m}}{|F''(\BFlambda^*)|}+\inner{\BFone_n}{|G''(\BFmu^*)|}+2\tilde{c}}{2}\left(\E\left[\bigg|f^3(\frac{\BFq}{\eta^{\alpha}})\bigg|\right]+\E\left[\bigg|g^3(\frac{\BFq}{\eta^{\alpha}})\bigg|\right]\right).
    \end{align*}
Dividing both sides by $\eta^{2\beta -1}$, as $\beta < 1$, we have
    \begin{align*}
        \frac{\mathcal{A}_3}{\eta^{2\beta-1}} \geq &-\eta^{\beta-1}\frac{\inner{\BFone_{m}}{|F''(\BFlambda^*)|}+\inner{\BFone_n}{|G''(\BFmu^*)|}+2\tilde{c}}{2}\left(\sum_{j=1}^m\E\left[\bigg|f_j^3(\frac{\BFq}{\eta^{\alpha}})\bigg|\right]+\sum_{i=1}^n\E\left[\bigg|g_i^3(\frac{\BFq}{\eta^{\alpha}})\bigg|\right]\right)
    \end{align*}
For simplicity, we let $K_1 = (\inner{\BFone_{m}}{|F''(\BFlambda^*)|}+\inner{\BFone_n}{|G''(\BFmu^*)|}+2\tilde{c})/{2}$. Combining everything, for all $\eta>\eta'$ (recall the definition of $\eta'$ in Claim~\ref{claim: uniformconvergence}), we have
\begin{align}
    \lr \geq{}& \eta^{2\beta-1}\left(\sum_{j=1}^m \frac{A_j}{2}\E\left[f_j^2\left(\frac{\BFq}{\eta^{\alpha}}\right)\right]+\sum_{i=1}^n \frac{B_i}{2}\E\left[g_i^2\left(\frac{\BFq}{\eta^{\alpha}}\right)\right]\right)\nonumber\\
    &-\eta^{2\beta-1}\eta^{\beta-1}K_1\left(\sum_{j=1}^m\E\left[\bigg|f_j^3(\frac{\BFq}{\eta^{\alpha}})\bigg|\right]+\sum_{i=1}^n\E\left[\bigg|g_i^3(\frac{\BFq}{\eta^{\alpha}})\bigg|\right]\right)+\E\left[\inner{\BFs}{\BFq}\right]  \label{eq: combiningeverything} \\
    \geq{}& \eta^{2\beta-1}\left(\sum_{j=1}^m \frac{A_j}{2}\E\left[f_j^2\left(\frac{\BFq}{\eta^{\alpha}}\right)\right]+\sum_{i=1}^n \frac{B_i}{2}\E\left[g_i^2\left(\frac{\BFq}{\eta^{\alpha}}\right)\right]\right)\nonumber\\
    &-\eta^{2\beta-1}\frac{\min_{i,j} \{A_j,B_i\}\epsilon}{4(\inner{\BFone_m}{{\BFGamma}^3}+\inner{\BFone_n}{{\BFPsi}^3})}\left(\sum_{j=1}^m\E\left[\bigg|f_j^3(\frac{\BFq}{\eta^{\alpha}})\bigg|\right]+\sum_{i=1}^n\E\left[\bigg|g_i^3(\frac{\BFq}{\eta^{\alpha}})\bigg|\right]\right)+\E\left[\inner{\BFs}{\BFq}\right].  \nonumber
\end{align}
By Lemma~\ref{lemma: boundedf2g2} and Condition \ref{ass: generalpricinggraph}\ref{condition: boundedfg}, for all $\eta>\max\{\eta_1,\eta',\mathcal{M}\}$, we have
\begin{align} 
        \lr\geq{}\ &  \eta^{2\beta-1}\left(\frac{\min_{i,j}\left\{A_j,B_i
      \right\} \epsilon}{2}-\frac{\min_{i,j}\left\{A_j,B_i
      \right\} \epsilon}{4}\right)+\frac{\min_{i \in N, j \in M}\{s_i^{(s)},s_j^{(c)}\}\eta^{1-\beta}\inner{\BFone_m}{\olam}}{2\left(\inner{\BFone_n}{{\BFPsi}}+\inner{\BFone_m}{{\BFGamma}}\right)}\nonumber \\
      &\qquad -0.5\max_{i \in N, j \in M}\{s_i^{(s)},s_j^{(c)}\} \nonumber\\
      {\geq}{}\ &  \inf_{\beta<1}\left\{\eta^{2\beta-1}\frac{\min_{i,j}\left\{A_j,B_i
      \right\} \epsilon}{4}+\frac{\min_{i \in N, j \in M}\{s_i^{(s)},s_j^{(c)}\}\eta^{1-\beta}\inner{\BFone_m}{\olam}}{2\left(\inner{\BFone_n}{{\BFPsi}}+\inner{\BFone_m}{{\BFGamma}}\right)}\right\}-0.5\max_{i \in N, j \in M}\{s_i^{(s)},s_j^{(c)}\}\nonumber \\
     {=}{}\ &\ 
     \eta^{1/3}\frac{3}{2^{2/3}}\left(\frac{\min_{i,j}\left\{A_j,B_i
      \right\} \epsilon}{4}\right)^{1/3}\left(\frac{\min_{i \in N, j \in M}\{s_i^{(s)},s_j^{(c)}\}\inner{\BFone_m}{\olam}}{2\left(\inner{\BFone_n}{{\BFPsi}}+\inner{\BFone_m}{{\BFGamma}}\right)}\right)^{2/3} 
      \\
      &- 0.5\max_{i \in N, j \in M}\{s_i^{(s)},s_j^{(c)}\}, \nonumber
\end{align}
where
the last equality follows as the coefficient of the terms $\eta^{2\beta-1}$ and $\eta^{1-\beta}$ are strictly positive and the infimum is achieved when $2\beta-1=1-\beta$. Let 
\begin{align}
    K := \frac{3}{2^{2/3}}\left(\frac{\min_{i,j}\left\{A_j,B_i
      \right\} \epsilon}{4}\right)^{1/3}\left(\frac{\min_{i \in N, j \in M}\{s_i^{(s)},s_j^{(c)}\}\inner{\BFone_m}{\olam}}{2\left(\inner{\BFone_n}{{\BFPsi}}+\inner{\BFone_m}{{\BFGamma}}\right)}\right)^{2/3}, \label{eq: constant_cube_root}
\end{align}
and $K > 0$ by Claim~\ref{claim: uniformconvergence} so the proof is complete.
$\Halmos$
\endproof
\subsection{Proof of Proposition \ref{prop: twopricesinglelinktwosidedqueue}}
\label{appendix: twopricesinglelinktwosidedqueue}
\begin{lemma}\label{lemma: lowerboundeqtwoprice}
Under the hypothesis of Proposition \ref{prop: twopricesinglelinktwosidedqueue}, if $\alpha+\beta>1$, there exists $\tilde{c}_1>0$, $\tilde{c}_2>0$ such that for all $\eta \geq 1$, it holds that
\begin{align}
& \E[q^{(s)}+q^{(c)}] \geq\max \left\{\frac{1}{2}\eta^{\alpha},\frac{\eta^{1-\beta}\lambda^*}{2(\theta+\phi)}-0.5\right\} \nonumber\\
  &  \E\left[f^2\left(\frac{\BFq}{\eta^{\alpha}}\right)+g^2\left(\frac{\BFq}{\eta^{\alpha}}\right)\right] \geq \tilde{c}_1 \eta^{1-\alpha-\beta} \nonumber\\
   & \E\left[\bigg|f^3\left(\frac{\BFq}{\eta^{\alpha}}\right)\bigg|+\bigg|g^3\left(\frac{\BFq}{\eta^{\alpha}}\right)\bigg|\right] \leq \tilde{c}_2 \eta^{1-\alpha-\beta}, \nonumber
\end{align}
where $f(\BFq/\eta^{\alpha})=-\theta\mathbbm{1}_{\left\{q^{(c)} \geq \eta^{\alpha}\right\}}$ and $g(\BFq/\eta^{\alpha})=-\phi\mathbbm{1}_{\left\{q^{(s)} \geq \eta^{\alpha}\right\}}$.
\end{lemma}
\proof{Proof of Lemma~\ref{lemma: lowerboundeqtwoprice}.}
Recall that the two-price policy for a single-link two-sided queue has the form
\begin{align*}
    \lambda^{\eta}(q)=\begin{cases}
    \eta \lambda^{*} &\textit{if } q \leq \lceil\eta^{\alpha}\rceil \\
    \eta\lambda^{*}-\theta\eta^{\beta} &\textit{otherwise}
    \end{cases} \quad
      \mu^{\eta}(q)=\begin{cases}
    \eta \mu^{*} &\textit{if } q \geq -\lceil\eta^{\alpha}\rceil \\
    \eta\mu^{*}-\phi\eta^{\beta} &\textit{otherwise},
    \end{cases}
\end{align*}
where $q=q^{(c)}-q^{(s)}$. Also, we have $\min\{q^{(c)},q^{(s)}\}=0$ for single-link queues. Consider a birth and death process denoted by $\{\tilde{q}(t) : t \geq 0\}$ with the state space $\mathbb{Z}$ 
and the transition rates given by
\begin{align*}
    \tilde{\lambda}^{\eta}(q)=\begin{cases}
    \eta \lambda^{*} &\textit{if } q \leq \lceil\eta^{\alpha}\rceil \\
    0 &\textit{otherwise}
    \end{cases} \quad
      \tilde{\mu}^{\eta}(q)=\begin{cases}
    \eta \mu^{*} &\textit{if } q \geq -\lceil\eta^{\alpha}\rceil \\
    0 &\textit{otherwise}.
    \end{cases}
\end{align*}
Because $\tilde{\lambda}(q) \leq \lambda(q)$ and $\tilde{\mu}(q) \geq \mu(q)$ for all $q \geq 0$ and $\tilde{\lambda}(q) \geq \lambda(q)$ and $\tilde{\mu}(q) \leq \mu(q)$ for all $q \leq 0$,  we can couple the two systems such that $|\tilde{q}(t)|\leq |q(t)| $ for all $t \geq 1$. Thus, we have $\Pr[|q(t)|\leq \mathcal{K}] \leq \Pr[|\tilde{q}(t)|\leq \mathcal{K}]$ for all $t \in \mathbb{Z}_+$. Since $\{q(t):t \geq 0\}$ is irreducible and positive recurrent by assumption and $\{\tilde{q}(t): t\geq 0\}$ has only one irreducible and positive recurrent class by construction, the stationary distributions exist for both $q$ and $\tilde{q}$. We have
\begin{align*}
    \Pr[|q(\infty)|\leq \mathcal{K}] \leq \Pr[|\tilde{q}(\infty)|\leq \mathcal{K}] \ \forall \mathcal{K} \in \mathbb{R} \ \Rightarrow \E[|q(\infty)|] \geq \E[|\tilde{q}(\infty)|],
\end{align*}
where $q(\infty)$ and $\tilde{q}(\infty)$ are random variables following the stationary distributions of $\{q(t):t \geq 0\}$ and $\{\tilde{q}(t):t \geq 0\}$, respectively. It is easy to verify that $\E[|\tilde{q}(\infty)|]=(\lceil\eta^{\alpha}\rceil(\lceil\eta^{\alpha}\rceil+1))/(2\lceil\eta^{\alpha}\rceil+1) \geq \lceil\eta^{\alpha}\rceil/2$. Thus, we have
\begin{align}
    \E[|q^{(s)}(\infty)+q^{(c)}(\infty)] {=} \E[|q(\infty)|] \geq \E[|\tilde{q}(\infty)|]\geq\frac{\eta^{\alpha}}{2}. \nonumber
\end{align}
Also, by Eq~\eqref{eq: twopricecondition12}, the two-price policy satisfies Condition \ref{ass: generalpricinggraph}\ref{condition: boundedfg} and \ref{ass: generalpricinggraph}\ref{condition: posrecsystem}. By Lemma~\ref{lemma: boundedf2g2}, we have
\begin{align}
    \E[|q(\infty)|] {=}\E[q^{(s)}(\infty)+q^{(c)}(\infty)] &\geq \frac{\eta^{1-\beta}\lambda^{*}}{2(\theta+\phi)}-0.5. \nonumber
\end{align}
where $(q^{(s)}(\infty), q^{(c)}(\infty))$ are random variables following the stationary distributions of $\{q^{(s)}(t),q^{(c)}(t): t \geq 0\}$. The first equality follows as $|q(t)|=q^{(s)}(t)+q^{(c)}(t)\ a.s$.
This proves the first inequality in the lemma.
Now, we will lower bound $\Pr[|q|\geq\eta^{\alpha}]$ for a given two-price policy. Since $\lambda^{*}=\mu^{*}$, the stationary distribution under a given two-price policy denoted by $\{\pi_i\}_{i \in \mathbb{Z}}$ is given by
\begin{align*}
    \pi_{k}&=\begin{cases}
    \pi_0  &\forall k \leq \lceil\eta^{\alpha}\rceil, k \geq -\lceil\eta^{\alpha}\rceil \\
    \pi_0\left(1-\frac{\theta\eta^{\beta}}{\eta\mu^{*}}\right)^{k-\lceil\eta^{\alpha}\rceil}  &\forall k > \lceil\eta^{\alpha}\rceil \\
     \pi_0\left(1-\frac{\phi\eta^{\beta}}{\eta\lambda^{*}}\right)^{-k-\lceil\eta^{\alpha}\rceil}  &\forall k < -\lceil\eta^{\alpha}\rceil. \\
    \end{cases}
\end{align*}
This implies that
\begin{align}
    \Pr[|q(\infty)|\geq\eta^{\alpha}]&=\frac{\eta^{1-\beta}(\mu^*/\theta+\lambda^*/\phi)}{\eta^{1-\beta}(\mu^*/\theta+\lambda^*/\phi)+2\lceil\eta^{\alpha}\rceil-1} \nonumber\\
    &{\geq} \frac{(\mu^*/\theta+\lambda^*/\phi)}{(\mu^*/\theta+\lambda^*/\phi)+3}\eta^{1-\alpha-\beta}\overset{\Delta}{=}\frac{\tilde{c}_1}{\min\{\theta^2,\phi^2\}} \eta^{1-\alpha-\beta} \quad \forall \eta\geq1, \nonumber
\end{align}
where the inequality follows as $1-\alpha-\beta<0$. This implies that
\begin{align*}
    \E\left[f^2\left(\frac{\BFq}{\eta^{\alpha}}\right)\right]+ \E\left[g^2\left(\frac{\BFq}{\eta^{\alpha}}\right)\right]&= \theta^2 \Pr\left[q^{(c)}(\infty)\geq\eta^{\alpha}\right]+\phi^2 \Pr\left[q^{(s)}(\infty)\geq\eta^{\alpha}\right] \\
      &\geq \min\{\theta^2,\phi^2\}\Pr\left[|q(\infty)|\geq\eta^{\alpha}\right] \geq \tilde{c}_1 \eta^{1-\alpha-\beta}.
\end{align*}
Thus, the second inequality of the lemma is proved.
Similarly, we have
\begin{align}
    \Pr[|q(\infty)|\geq\eta^{\alpha}]&\leq\frac{\eta^{1-\alpha-\beta}(\mu^*/\theta+\lambda^*/\phi)}{\eta^{1-\alpha-\beta}(\mu^*/\theta+\lambda^*/\phi)+2-\eta^{-\alpha}}  \nonumber\\
    & {\leq} (\mu^*/\theta+\lambda^*/\phi)\eta^{1-\alpha-\beta}\overset{\Delta}{=}\frac{\tilde{c}_2}{\max\{\theta^3,\phi^3\}} \eta^{1-\alpha-\beta} \quad \forall \eta\geq 1, \nonumber
\end{align}
where 
we replace 
$\lceil\eta^{\alpha}\rceil$ by $\eta^{\alpha}$ in the first inequality, and
the second inequality follows as $\alpha\geq 0$ and $1-\alpha-\beta<0$. This implies that
\begin{align*}
    \E\left[\bigg|f^3\left(\frac{\BFq}{\eta^{\alpha}}\right)\bigg|\right]+ \E\left[\bigg|g^3\left(\frac{\BFq}{\eta^{\alpha}}\right)\bigg|\right]&= \theta^3 \Pr\left[q^{(c)}(\infty)\geq\eta^{\alpha}\right]+\phi^3 \Pr\left[q^{(s)}(\infty)\geq\eta^{\alpha}\right] \\
      &\leq \max\{\theta^3,\phi^3\}\Pr\left[|q(\infty)|\geq\eta^{\alpha}\right] \leq \tilde{c}_2 \eta^{1-\alpha-\beta}.
\end{align*}
This proves the last inequality of the lemma.
$\Halmos$
\endproof
\proof{Proof of Proposition \ref{prop: twopricesinglelinktwosidedqueue}.}
When $\alpha+\beta \leq 1$, by Theorem \ref{theorem: lowerboundmultiplelink}, we know that $\lr \geq K\eta^{1/3}$. So we will only prove Proposition~\ref{prop: twopricesinglelinktwosidedqueue} when $\alpha+\beta>1$.
As the two-price policy is a special case of the general pricing policy given by Eqs~(\ref{eq: pricinggraphgeneralcust}, \ref{eq: pricinggraphgeneralser}), it holds by Eq \eqref{eq: combiningeverything} that
\begin{align}
    \lr \geq{}& \eta^{2\beta-1}\left( \frac{A}{2}\E\left[f^2\left(\frac{\BFq}{\eta^{\alpha}}\right)\right]+ \frac{B}{2}\E\left[g^2\left(\frac{\BFq}{\eta^{\alpha}}\right)\right]\right)\nonumber\\
    &-\eta^{2\beta-1}\eta^{\beta-1}K_1\left(\E\left[\bigg|f^3(\frac{\BFq}{\eta^{\alpha}})\bigg|\right]+\E\left[\bigg|g^3(\frac{\BFq}{\eta^{\alpha}})\bigg|\right]\right)+\E\left[\inner{\BFs}{\BFq}\right] \nonumber \\
    \geq{}& \eta^{\beta-\alpha}\left(\min\{A,B\}\frac{\tilde{c}_1}{2}-\tilde{c}_2 K_1\eta^{\beta-1}\right)+\min_{i \in N, j \in M}\{s^{(s)}_i,s^{(c)}_j\}\max \left\{\frac{1}{2}\eta^{\alpha},\frac{\eta^{1-\beta}\lambda^*}{2(\theta+\phi)}-0.5\right\}. \nonumber
    \end{align}
Letting $\eta>\left(\frac{\min\{A,B\}\tilde{c}_1}{4\tilde{c}_2 K_1}\right)^{1/(\beta-1)}$ and optimizing over all possible values of $\alpha$ and $\beta$, we get
    \begin{align}
        \lr \geq & \inf_{\alpha \geq 0,\ 1 - \beta < \alpha}\left\{\eta^{\beta-\alpha}\min\{A,B\}\frac{\tilde{c}_1}{4}+\min_{i \in N, j \in M}\{s^{(s)}_i,s^{(c)}_j\}\max \left\{\frac{1}{2}\eta^{\alpha},\frac{\eta^{1-\beta}\lambda^*}{2(\theta+\phi)}-0.5\right\}\right\} \nonumber\\
        {=}& \min\left\{\min\{A,B\}\frac{\tilde{c}_1}{4},
        \min_{i \in N, j \in M}\{s^{(s)}_i,s^{(c)}_j\}\max \left\{\frac{1}{2},\frac{\lambda^*}{2(\theta+\phi)}-\frac{1}{2\eta^{1/3}}\right\}
        \right\}\eta^{1/3}, \label{eq: two-price-single-link-lower-bound}
    \end{align}
    where the infimum is approached when $\alpha=1/3$, $\beta = 2/3$. If $\alpha\geq 1/3$, then the holding-cost term is $\Omega(\eta^\alpha)=\Omega(\eta^{1/3})$. If $\alpha<1/3$, then $\alpha+\beta>1$ implies
$\beta-\alpha>1-2\alpha>1/3$, so the revenue-loss term is
$\Omega(\eta^{\beta-\alpha})=\Omega(\eta^{1/3})$.
$\Halmos$
\endproof
\section{Further Analysis on Max-Weight Matching in the Large-System Regime}
\label{appendix: further-analysis-max-weight}
\proof{Proof of Lemma \ref{lemma: positive_chi}.}
By the CRP condition, there exists a constant $\delta>0$ such that for all $J \subsetneq {N}$ and for all $I \subsetneq {N}$, it holds that
\begin{align}
    \sum_{j \in J} \lambda_j^* < \sum_{i : \exists j \in J, (i,j) \in E} \mu_i^*-\delta, \quad \sum_{i \in I} \mu_i^* < \sum_{j: \exists i \in I, (i,j) \in E} \lambda_j^*-\delta. \label{eq: strict_hallscondition}
\end{align}
We claim that there exists
$\BFchi \in \mathbb{R}_+^{m \times n}$ satisfying the following equalities (for $\delta$ small enough):
\begin{subequations} \label{eq: new_chi}
\begin{align}
    \lambda_j^*-\frac{\delta}{n^2}|N(j)|&=\sum_{i=1}^n \chi_{ij}, \quad \forall j \in {M} \\ \mu_i^*-\frac{\delta}{n^2}|N(i)|&=\sum_{j=1}^m \chi_{ij}, \quad \forall i \in {N} \\
    \chi_{ij}&=0, \quad \forall (i,j) \notin E.
\end{align}
\end{subequations}
Such a $\BFchi$ exists because the left-hand sides $(\lambda_j^*-\frac{\delta}{n^2}|N(j)|)_{j=1}^m$ and $(\mu_i^*-\frac{\delta}{n^2}|N(i)|)_{i=1}^n$ satisfy the Hall's condition, which is 
implied by \eqref{eq: strict_hallscondition}. 
Let $\chi^*_{ij}=\frac{\delta}{n^2}+\chi_{ij}$ for all $(i,j) \in E$ and $\chi^*_{ij}=0$ otherwise.
Substituting $\BFchi$ for $\BFchi^*$ in \eqref{eq: new_chi} gives
\begin{subequations} \label{eq: same_as_fluid}
\begin{align}
    \lambda_j^*&=\sum_{i=1}^n \chi_{ij}^*, \quad \forall j \in {M} \\ \mu_i^*&=\sum_{j=1}^m \chi_{ij}^*, \quad \forall i \in {N} \\
    \chi_{ij}^*&=0, \quad \forall (i,j) \notin E.
\end{align}
\end{subequations}
Note that \eqref{eq: same_as_fluid} is the same as the constraints in the fluid problem \eqref{eq:fluid_opt}, so $(\BFlambda^*,\BFmu^*,\BFchi^*)$ is a feasible solution. 
In addition, $\inner{F(\BFlambda^*)}{\BFlambda^*}-\inner{G(\BFmu^*)}{\BFmu^*}$ is the optimal objective function value, so $(\BFlambda^*,\BFmu^*,\BFchi^*)$ is an optimal solution. As $\delta>0$, we have $\chi^*_{ij}>0$ for all $(i,j) \in E$. 
$\Halmos$
\endproof
To prove the state space collapse,
we consider the Lyapunov function 
\begin{equation}\label{eq:def-V}
    U(\BFq)=\sqrt{\inner{\BFone_{2n}}{\BFq^2}-\frac{z^2}{n}}.
\end{equation}
(Recall the definition of $z$ from Eq~\eqref{eq: imbalance}.)
In addition, define the following vectors of size $n$: $\BFq_{\parallel}^{(s)}\overset{\Delta}{=}\frac{1}{n}\inner{\BFone_n}{\BFq^{(s)}}\BFone_n$, $\BFq_{\perp}^{(s)}\overset{\Delta}{=}\BFq^{(s)}-\BFq_{\parallel}^{(s)}$ and $\BFq_{\parallel}^{(c)}\overset{\Delta}{=}\frac{1}{n}\inner{\BFone_n}{\BFq^{(c)}}\BFone_n$, $\BFq_{\perp}^{(c)}\overset{\Delta}{=}\BFq^{(c)}-\BFq_{\parallel}^{(c)}$. We will first prove the following lemmas, which will assist us in proving Proposition \ref{prop: max-weight-optimal-fluid-SSC} and Proposition \ref{prop: max-weight-optimal-two-price-SSC}. Let $S_{\text{mw}}$ be the set of states such that no compatible customer-server pairs are waiting in the queue, i.e.,
\begin{align*}
    S_{\text{mw}} = \left\{\BFq \in \mathbb{Z}_+^{m+n}: q^{(s)}_iq^{(c)}_j = 0 \quad \forall (i, j) \in E\right\}.
\end{align*}
Note that if we employ the max-weight matching policy, then, the state space under the induced Markov chain will always be a subset of $S_{\text{mw}}$.
\begin{lemma} \label{lemma: perp_representation}
 Under any pricing policy and max-weight matching policy, for all $\BFq \in S_{\text{mw}}$, we have
 \begin{align*}
    U(\BFq)=\|\BFq_{\perp}^{(s)}-\BFq_{\perp}^{(c)}\|.
\end{align*}
\end{lemma}
\proof{Proof of Lemma~\ref{lemma: perp_representation}.}
With elementary algebra, we can show
\begin{align*} 
    U^2(\BFq)&=\inner{\BFone_n}{(\BFq^{(s)})^2}-\frac{1}{n}\inner{\BFone_n}{\BFq^{(s)}}^2+\inner{\BFone_n}{(\BFq^{(c)})^2}-\frac{1}{n}\inner{\BFone_n}{\BFq^{(c)}}^2+\frac{2}{n}\inner{\BFone_n}{\BFq^{(s)}}\inner{\BFone_n}{\BFq^{(c)}} \\
    &=\inner{\BFone_n}{(\BFq_{\perp}^{(s)})^2}+\inner{\BFone_n}{(\BFq_{\perp}^{(c)})^2}+\frac{2}{n}\inner{\BFone_n}{\BFq^{(s)}}\inner{\BFone_n}{\BFq^{(c)}} \\
    & = \inner{\BFone_n}{(\BFq_{\perp}^{(s)})^2}+\inner{\BFone_n}{(\BFq_{\perp}^{(c)})^2}-2\inner{\BFq_\perp^{(s)}}{\BFq_\perp^{(c)}} \\
    &=\| \BFq_\perp^{(s)}-\BFq_\perp^{(c)}\|^2,
\end{align*}
where the third equality above holds because
\begin{align*}
    \inner{\BFq_\perp^{(s)}}{\BFq_\perp^{(c)}}&=\sum_{i=1}^n \left(q_i^{(s)}-\frac{1}{n}\inner{\BFone_n}{\BFq^{(s)}}\right)\left(q_i^{(c)}-\frac{1}{n}\inner{\BFone_n}{\BFq^{(c)}}\right) \\
    &=\sum_{i=1}^n q_i^{(s)}q_i^{(c)}-\frac{1}{n}\inner{\BFone_n}{\BFq^{(s)}}\inner{\BFone_n}{\BFq^{(c)}} \\
    &{=}-\frac{1}{n}\inner{\BFone_n}{\BFq^{(s)}}\inner{\BFone_n}{\BFq^{(c)}}.
\end{align*}
The last equality above follows because $E$ contains edges connecting type $i$ server to type $i$ customer by Condition~\ref{condition: perfect_matching}. Thus, for the max-weight matching policy $(\BFq \in S_{\text{mw}})$, either $q_i^{(s)} = 0$ or $q_i^{(c)} = 0$.
\Halmos
\endproof
\begin{lemma} \label{lemma: imp_term_SSC}
Under any pricing policy and max-weight matching policy, for all $\BFq \in S_{\text{mw}}$, we have
\begin{align*}
   T(\BFq) := \sum_{(i,j) \in E} \left(q_i^{(s)}+q_j^{(c)}-\max_{i' : (i',j) \in E} q_{i'}^{(s)}-\max_{j' : (i,j') \in E} q_{j'}^{(c)}\right) \leq -\frac{1}{n}\sqrt{\inner{\BFone_{2n}}{\BFq^2}-\frac{z^2}{n}}.
\end{align*}
\end{lemma}
\proof{Proof of Lemma~\ref{lemma: imp_term_SSC}.}
By the definition of $T(\BFq)$, we have
\begin{align*}
    T(\BFq)= \underbrace{\sum_{(i,j)\in E}\left(q_i^{(s)}-\max_{i' : (i',j) \in E} q_{i'}^{(s)}\right)}_{T_1}+\underbrace{\sum_{(i,j)\in E}\left(q_j^{(c)}-\max_{j' : (i,j') \in E} q_{j'}^{(c)}\right)}_{T_2}.
\end{align*}
As the graph $E$ is connected (implied by the CRP condition), consider a path $P$ from $i_{\min}=\min_{i \in {N}} q_i^{(s)}$ to $i_{\max}=\max_{i \in {N}} q_i^{(s)}$. We have
\begin{align*}
    T_1&=\sum_{(i,j)\in E}\left(q_i^{(s)}-\max_{i' : (i',j) \in E} q_{i'}^{(s)}\right) = \sum_{(i,j)\in P}\left(q_i^{(s)}-\max_{i' : (i',j) \in E} q_{i'}^{(s)}\right)+\sum_{(i,j)\in E \backslash P}\left(q_i^{(s)}-\max_{i' : (i',j) \in E} q_{i'}^{(s)}\right) \\
    &\leq \sum_{(i,j)\in P}\left(q_i^{(s)}-\max_{i' : (i',j) \in E} q_{i'}^{(s)}\right) \leq q_{i_{\min}}^{(s)}-q_{i_{\max}}^{(s)},
\end{align*}
where the last inequality follows from telescoping sum. Similarly, $T_2 \leq q_{j_{\min}}^{(c)}-q_{j_{\max}}^{(c)}$. Thus, we have
\begin{align*}
    T \leq q_{i_{\min}}^{(s)}+q_{j_{\min}}^{(c)}-q_{i_{\max}}^{(s)}-q_{j_{\max}}^{(c)}.
\end{align*}
Note that, $q_{i_{\min}}^{(s)}$ and $q_{j_{\min}}^{(c)}$ cannot both be non-zero by the definition of the max-weight algorithm $(\BFq \in S_{\text{mw}})$. Without loss of generality, we assume  $q_{j_{\min}}^{(c)}=0$. 
We consider two cases:\\
\textit{Case I:} $q_{i_{\min}}^{(s)}=0$. We have
\begin{align*}
    T(\BFq) \leq -\left(q_{i_{\max}}^{(s)}+q_{j_{\max}}^{(c)}\right)\leq -\frac{1}{n}\inner{\BFone_{2n}}{\BFq} \leq -\frac{1}{n} \|\BFq\| \leq -\frac{1}{n} U(\BFq).
\end{align*}
\textit{Case II:} $q_{i_{\min}}^{(s)}>0$. This implies that $q_i^{(s)}>0$ for all $i \in {N}$, so we must have $q_j^{(c)}=0$ for all $j \in {N}$ by the definition of the max-weight algorithm. Thus, we have
\begin{align*}
     T(\BFq) \leq q_{i_{\min}}^{(s)}-q_{i_{\max}}^{(s)} \leq \frac{1}{n}\inner{\BFone_n}{\BFq^{(s)}}-q_{i_{\max}}^{(s)} \leq -\frac{1}{n}\|q_\perp^{(s)}\|=-\frac{1}{n}\|q_\perp^{(s)}-q_\perp^{(c)}\|=-\frac{1}{n}U(\BFq). \Halmos
\end{align*}
\endproof
To prepare for the next lemma, recall the definition 
\[
    V(\BFq)=\inner{\BFone_{2n}}{\BFq^2},\;
    V^{(c)}(\BFq)=\inner{\BFone_{n}}{(\BFq^{(c)})^2},\;
    V^{(s)}(\BFq)=\inner{\BFone_{n}}{(\BFq^{(s)})^2}.
\]
We also define
\begin{equation}
    W_z(\BFq)=\frac{z^2}{n}.
    \label{eq:def-W}
\end{equation}
\begin{lemma} \label{lemma: bounded_drift}
Under any pricing policy and max-weight matching policy, for all $\BFq \in S_{\text{mw}}$, the drift of $U(\BFq)$, defined as $\Delta U(\BFq)=(U(\BFq(k+1))-U(\BFq(k)))\mathbbm{1}_{\{\BFq(k)=\BFq\}}$, satisfies
\begin{align*}
    |\Delta U(\BFq)| \leq 4, \quad
    \Delta U(\BFq) \leq \frac{1}{2U(\BFq)}(\Delta V(\BFq)-\Delta W_z(\BFq)).
\end{align*}
\end{lemma}
\proof{Proof of Lemma~\ref{lemma: bounded_drift}.}
By expanding the drift term, we have
\begin{align}
    \lefteqn{|\Delta U(\BFq)|=|U(\BFq(k+1))-U(\BFq(k))|\mathbbm{1}_{\BFq(k)=\BFq}}\nonumber \\
    \overset{(a)}{=}{}&\big|\| \BFq_\perp^{(s)}(k+1)-\BFq_\perp^{(c)}(k+1)\|-\| \BFq_\perp^{(s)}(k)-\BFq_\perp^{(c)}(k)\|\big|\mathbbm{1}_{\BFq(k)=\BFq}\nonumber \\
    \overset{(b)}{\leq}{}& \|\BFq_\perp^{(s)}(k+1)-\BFq_\perp^{(c)}(k+1)-\BFq_\perp^{(s)}(k)+\BFq_\perp^{(c)}(k)\|\mathbbm{1}_{\BFq(k)=\BFq}\nonumber \\
    ={}&\|\BFq^{(s)}(k+1)-\BFq^{(s)}(k)-\BFq^{(s)}_\parallel(k+1)+\BFq^{(s)}_\parallel(k)-\BFq^{(c)}(k+1)+\BFq^{(c)}(k)+\BFq^{(c)}_\parallel(k+1)-\BFq^{(c)}_\parallel(k)\|\mathbbm{1}_{\BFq(k)=\BFq}\nonumber \\
    \overset{(c)}{\leq}{}& \left(\|\BFq^{(s)}(k+1)-\BFq^{(s)}(k)\|+\|\BFq^{(s)}_\parallel(k+1)-\BFq^{(s)}_\parallel(k)\|+\|\BFq^{(c)}(k+1)-\BFq^{(c)}(k)\|\right.\nonumber\\
    &\left.+\|\BFq^{(c)}_\parallel(k+1)-\BFq^{(c)}_\parallel(k)\|\right)\mathbbm{1}_{\BFq(k)=\BFq} \nonumber\\
    \overset{(d)}{\leq}{}& 4 \|\BFq(k+1)-\BFq(k)\|\mathbbm{1}_{\BFq(k)=\BFq} \overset{(e)}{\leq} 4,
\end{align}
where $(a)$ follows from Lemma \ref{lemma: perp_representation} (using $\BFq \in S_{\text{mw}}$); $(b)$ and $(c)$ follows by the triangle inequality; $(d)$ follows because projection onto a subspace is non-expansive and $\|\BFq^{(l)}(k+1)-\BFq^{(l)}(k)\| \leq \|\BFq(k+1)-\BFq(k)\|$ for $l \in \{1,2\}$; finally, $(e)$ follows as there can be at most one arrival and one matching (under max-weight matching) in one time epoch of the unformized DTMC. This prove the first part of the lemma.
Now we will show the second part of the lemma. Expanding the drift term, we have
\begin{align}
    \Delta U(\BFq)&=\left(\|\BFq_\perp^{(s)}(k+1)-\BFq_\perp^{(c)}(k+1)\|-\|\BFq_\perp^{(s)}(k)-\BFq_\perp^{(c)}(k)\|\right)\mathbbm{1}_{\BFq(k)=\BFq} \nonumber \\
    &=\left(\sqrt{\|\BFq_\perp^{(s)}(k+1)-\BFq_\perp^{(c)}(k+1)\|^2}-\sqrt{\|\BFq_\perp^{(s)}(k)-\BFq_\perp^{(c)}(k)\|^2}\right)\mathbbm{1}_{\BFq(k)=\BFq} \nonumber  \\
    &{\leq} \frac{1}{2\|\BFq_\perp^{(s)}(k)-\BFq_\perp^{(c)}(k)\|}\left(\|\BFq_\perp^{(s)}(k+1)-\BFq_\perp^{(c)}(k+1)\|^2-\|\BFq_\perp^{(s)}(k)-\BFq_\perp^{(c)}(k)\|^2\right)\mathbbm{1}_{\BFq(k)=\BFq} \nonumber \\
    &=\frac{1}{2\|\BFq_\perp^{(s)}(k)-\BFq_\perp^{(c)}(k)\|}\left(\inner{\BFone_{2n}}{\BFq^2(k+1)}-\inner{\BFone_{2n}}{\BFq^2(k)}-\frac{z^2(k+1)}{n}+\frac{z^2(k)}{n}\right)\mathbbm{1}_{\BFq(k)=\BFq} \nonumber \\
    &=\frac{1}{2\|\BFq_\perp^{(s)}-\BFq_\perp^{(c)}\|}\left(\Delta V(\BFq)-\Delta W_z(\BFq)\right), \nonumber
\end{align}
where the inequality follows as $h(x)=\sqrt{x}$ is concave and hence $h(y)-h(x) \leq (y-x)h'(x)=(y-x)/(2\sqrt{x})$.
\Halmos
\endproof
\subsection{Proposition \ref{prop: max-weight-optimal-fluid-SSC}}
We start with a lemma that established the state space collapse.
Recall the definition of the imbalance process $z(k)$ in Eq~\eqref{eq: imbalance}.
\begin{lemma}[State Space Collapse] \label{lemma: SSC_fluid}
Under the fluid pricing policy and max-weight matching policy, for all $\eta \geq 1$, $r \in \mathbb{Z}_+$ there exists a constant $\mathcal{F}_r$ (independent of $\eta$) such that 
\begin{align*}
    \E\left[\left(\inner{\BFone_{2n}}{\BFq^2(\infty)}-\frac{z^2(\infty)}{n}\right)^{r/2}\right] \leq \mathcal{F}_r.
\end{align*}
\end{lemma}
\proof{Proof of Lemma~\ref{lemma: SSC_fluid}.}
First note that the state space of the Markov chain induced by the fluid pricing and max-weight matching policy is a subset of $S_{\text{mw}}$ as max-weight matches all compatible customer-server pairs. By Lemma \ref{lemma: bounded_drift}, the drift of the Lyapunov function $U(\BFq)$ is uniformly bounded. In the proof below, we want to show that the drift is negative outside an appropriately defined finite set. We start by calculating the drift of $V(\BFq)$. By \eqref{eq: drift_const_price}, the drift of $V^{(c)}(\BFq)$ for the uniformized DTMC is
\begin{align*}
   c \E\left[\Delta V^{(c)}(\BFq) \mid \BFq(k)=\BFq\right]
   &\leq \sum_{i=1}^n \mu_i^* + \sum_{j=1}^n \lambda_j^*+2\sum_{(i,j)\in E} \chi_{ij}^{*}\left[q_j^{(c)}\mathbbm{1}_{\left\{q_j^{(c)}< q_{\max}^{\eta}\right\}}-\max_{j':(i,j') \in E}q_{j'}^{(c)}\right].
\end{align*}
Similarly, for $V^{(s)}(\BFq)$ we will have
\begin{align*}
     c \E\left[\Delta V^{(s)}(\BFq)\mid \BFq(k)=\BFq\right]
   &\leq \sum_{i=1}^n \mu_i^* + \sum_{j=1}^n \lambda_j^*+2\sum_{(i,j)\in E} \chi_{ij}^{*}\left[q_i^{(s)}\mathbbm{1}_{\left\{q_i^{(s)}< q_{\max}^{\eta}\right\}}-\max_{i':(i',j) \in E}q_{i'}^{(s)}\right].
\end{align*}
Adding the two inequalities above, we have
\begin{align*}
 \lefteqn{c \E\left[\Delta V(\BFq)\mid \BFq(k)=\BFq\right]} \\
     \leq{}& 2\sum_{i=1}^n \mu_i^* + 2\sum_{j=1}^n \lambda_j^*+2\sum_{(i,j)\in E} \chi_{ij}^{*}\left[q_j^{(c)}\mathbbm{1}_{\left\{q_j^{(c)}< q_{\max}^{\eta}\right\}}+q_i^{(s)}\mathbbm{1}_{\left\{q_i^{(s)}< q_{\max}^{\eta}\right\}}-\max_{j':(i,j') \in E}q_{j'}^{(c)}-\max_{i':(i',j) \in E}q_{i'}^{(s)}\right] \\
     \leq{}& 2\sum_{i=1}^n \mu_i^* + 2\sum_{j=1}^n \lambda_j^*-2q_{\max}^{\eta}\left(\sum_{i=1}^n \mu_i^*\mathbbm{1}_{\{q_i^{(s)}=q_{\max}^\eta\}}+2\sum_{j=1}^n \lambda_j^*\mathbbm{1}_{\{q_j^{(c)}=q_{\max}^\eta\}}\right) \\
     &+2\sum_{(i,j)\in E} \chi_{ij}^{*}\left[q_j^{(c)}+q_i^{(s)}-\max_{j':(i,j') \in E}q_{j'}^{(c)}-\max_{i':(i',j) \in E}q_{i'}^{(s)}\right] \\
     \leq{}&2\sum_{i=1}^n \mu_i^* + 2\sum_{j=1}^n \lambda_j^*-2q_{\max}^{\eta}\left(\sum_{i=1}^n \mu_i^*\mathbbm{1}_{\{q_i^{(s)}=q_{\max}^\eta\}}+\sum_{j=1}^n \lambda_j^*\mathbbm{1}_{\{q_j^{(c)}=q_{\max}^\eta\}}\right)-\frac{2}{n}\left(\min_{{i,j} \in E} \chi_{ij}^*\right)U(\BFq),
\end{align*}
where the last inequality uses Lemma~\ref{lemma: imp_term_SSC}.
The drift of $\Delta W_z(\BFq)$ can be bounded as
\begin{align*}
{c\E\left[\Delta W_z(\BFq)\mid \BFq(k)=\BFq\right]}
    ={}&\frac{1}{n}\E\left[z^2(k+1)-z^2(k)\mid \BFq(k)=\BFq\right] \nonumber \\
    ={}& \frac{1}{n}\left(\sum_{j=1}^n \lambda_j^*\mathbbm{1}_{\{q_j^{(c)}<q_{\max}^\eta\}}\left((z+1)^2-z^2\right)+\sum_{i=1}^n \mu_i^*\mathbbm{1}_{\{q_i^{(s)}<q_{\max}^\eta\}}\left((z-1)^2-z^2\right)\right) \\
    \geq{}& \frac{2z}{n}\left(\sum_{j=1}^n\lambda_j^*\mathbbm{1}_{\{q_j^{(c)}<q_{\max}^\eta\}} -\sum_{i=1}^n \mu_i^*\mathbbm{1}_{\{q_i^{(s)}<q_{\max}^\eta\}}\right) \\
    ={}& \frac{2z}{n}\left(\sum_{i=1}^n \mu_i^* \mathbbm{1}_{\{q_i^{(s)}=q_{\max}^\eta\}}-\sum_{j=1}^n \lambda_j^* \mathbbm{1}_{\{q_j^{(c)}=q_{\max}^\eta\}}\right) \\
    \geq & -2 q_{\max}^\eta\left(\sum_{i=1}^n \mu_i^*\mathbbm{1}_{\{q_i^{(s)}=q_{\max}^\eta\}}+\sum_{j=1}^n \lambda_j^*\mathbbm{1}_{\{q_j^{(c)}=q_{\max}^\eta\}}\right),
\end{align*}
where the third equality uses the fact that $\sum_{i=1}^n \mu^*_i = \sum_{j=1}^n \lambda^*_j$, and the last inequality uses the fact that $|z|/n \leq q^{\eta}_{\max}$ by the definition of the fluid pricing policy.
Combining the inequalities above, by Lemma \ref{lemma: bounded_drift}, 
for any $U(\BFq)> \frac{2n\left(\inner{\BFone_n}{\BFlambda^*}+\inner{\BFone_n}{\BFmu^*}\right)}{\min_{(i,j) \in E} \chi_{ij}^*}$,
we have
\begin{align*}
    \E\left[\Delta U(\BFq)\mid \BFq(k)=\BFq\right]  &\leq \frac{\E\left[\Delta V(\BFq)-\Delta W_z(\BFq)\right]}{2U(\BFq)} \\
    &\leq \frac{2\left(\inner{\BFone_n}{\BFlambda^*}+\inner{\BFone_n}{\BFmu^*}\right)}{c U(\BFq)}-\frac{2}{c n}\left(\min_{(i,j) \in E} \chi_{ij}^*\right)
    < -\frac{1}{c n}\left(\min_{(i,j) \in E} \chi_{ij}^*\right).
\end{align*}
Thus, we have established that the drift of $U(\BFq)$ is uniformly bounded and is negative whenever $U(\BFq)$ is large. 
By \citet[Lemma 1]{atilla_srikant},  we have
\begin{align}
   \E\left[U(\BFq(\infty))^r\right] =  \E\left[\left(\inner{\BFone_{2n}}{\BFq^2(\infty)}-\frac{z^2(\infty)}{n}\right)^{r/2}\right] \leq \mathcal{F}_r \quad \forall \ \eta \geq 1 \quad \forall r \in \mathbb{Z}_+. 
    \Halmos
\end{align}
\endproof
\proof{Proof of Proposition \ref{prop: max-weight-optimal-fluid-SSC}.}
In the steady state, we have $\E\left[\Delta W_z(\BFq(\infty))\right] = 0$. Expanding this equation, we have
\begin{align}
  & \quad \sum_{j=1}^n \lambda_j^* \Pr\left[q_j^{(c)}(\infty)<q_{\max}^\eta\right]+\sum_{i=1}^n \mu_i^* \Pr\left[q_i^{(s)}(\infty)<q_{\max}^\eta\right]\nonumber\\
   &=2\E\left[z(\infty)\left(\sum_{i=1}^n \left(\lambda_i^*\mathbbm{1}_{\{q_i^{(c)}(\infty)=q_{\max}^\eta\}}-\mu_i^*\mathbbm{1}_{\{q_i^{(s)}(\infty)=q_{\max}^\eta\}} \right)\right)\right]  \\
   &\leq 2nq_{\max}^\eta \sum_{i=1}^n\left( \lambda_i^* \Pr\left[q_i^{(c)}(\infty)=q_{\max}^\eta\right]+ \mu_i^* \Pr\left[q_i^{(s)}(\infty)=q_{\max}^\eta\right]\right), 
\end{align}
where the last inequality holds as $|z(\infty)| \leq \langle \BFone_{2n}, \BFq\rangle \leq 2n q_{\max}^\eta$. By the above equation, we have
\begin{align*}
& 2nq_{\max}^\eta \sum_{i=1}^n \left(\lambda_i^*\Pr[q_i^{(c)}(\infty)=q_{\max}^\eta]+\mu_i^*\Pr[q_i^{(s)}(\infty)=q_{\max}^\eta] \right)  \\
\geq & \sum_{j=1}^n \lambda_j^* \Pr\left[q_j^{(c)}(\infty)<q_{\max}^\eta\right]+\sum_{i=1}^n \mu_i^* \Pr\left[q_i^{(s)}(\infty)<q_{\max}^\eta\right]\\
= & \inner{\BFone_n}{\BFlambda^*}+\inner{\BFone_n}{\BFmu^*} - \sum_{i=1}^n \left(\lambda_i^*\Pr[q_i^{(c)}(\infty)=q_{\max}^\eta]+\mu_i^*\Pr[q_i^{(s)}(\infty)=q_{\max}^\eta] \right).
\end{align*}
Rearranging the terms, we get
\begin{equation}\label{eq:any-matching-lower-bound}
    \sum_{i=1}^n \left(\lambda_i^*\Pr[q_i^{(c)}(\infty)=q_{\max}^\eta]+\mu_i^*\Pr[q_i^{(s)}(\infty)=q_{\max}^\eta] \right) \geq \frac{\inner{\BFone_n}{\BFlambda^*}+\inner{\BFone_n}{\BFmu^*}}{2nq_{\max}^\eta+1}. 
\end{equation}
This proves the first part of the proposition.
Note that the proof so far applies to any matching algorithm. In the rest of the proof, we consider the max-weight algorithm specifically (so, the state space is a subset of $S_{\text{mw}}$). Again, by $\E\left[\Delta W_z(\BFq(\infty))\right] = 0$, we get
\begin{align}
& \quad \sum_{j=1}^n \lambda_j^* \Pr\left[q_j^{(c)}(\infty)<q_{\max}^\eta\right]+\sum_{i=1}^n \mu_i^* \Pr\left[q_i^{(s)}(\infty)<q_{\max}^\eta\right]\nonumber\\
   &=2\E\left[z(\infty)\left(\sum_{i=1}^n \left(\lambda_i^*\mathbbm{1}_{\{q_i^{(c)}(\infty)=q_{\max}^\eta\}}-\mu_i^*\mathbbm{1}_{\{q_i^{(s)}(\infty)=q_{\max}^\eta\}} \right)\right)\right]  \\
   & = 2n \E\Bigg[ \underbrace{\sum_{i=1}^n \left(q_{\perp_i}^{(s)}(\infty)-q_{\perp_i}^{(c)}(\infty)-\frac{z(\infty)}{n}\right)\left(\mu_i^*\mathbbm{1}_{\{q_i^{(s)}(\infty)=q_{\max}^\eta\}}-\lambda_i^*\mathbbm{1}_{\{q_i^{(c)}(\infty)=q_{\max}^\eta\}}\right)}_{(\dagger)}\Bigg] \nonumber\\
   & \quad -2 n \E\Bigg[ \underbrace{\sum_{i=1}^n \left(q_{\perp_i}^{(s)}(\infty)-q_{\perp_i}^{(c)}(\infty)\right)\left(\mu_i^*\mathbbm{1}_{\{q_i^{(s)}(\infty)=q_{\max}^\eta\}}-\lambda_i^*\mathbbm{1}_{\{q_i^{(c)}(\infty)=q_{\max}^\eta\}}\right)}_{(\ddagger)}\Bigg] \nonumber \\
   & = 2n  q_{\max}^\eta\sum_{i=1}^n\left( \lambda_i^* \Pr\left[q_i^{(c)}(\infty)=q_{\max}^\eta\right]+ \mu_i^* \Pr\left[q_i^{(s)}(\infty)=q_{\max}^\eta\right]\right) \nonumber\\
   & \quad -2n\E\left[\sum_{i=1}^n \left(q_{\perp_i}^{(s)}(\infty)-q_{\perp_i}^{(c)}(\infty)\right)\left(\mu_i^*\mathbbm{1}_{\{q_i^{(s)}(\infty)=q_{\max}^\eta\}}-\lambda_i^*\mathbbm{1}_{\{q_i^{(c)}(\infty)=q_{\max}^\eta\}}\right)\right]. \label{eq: for-lower-bound-fluid}
\end{align}

To see why the last the last equality holds, note that
\begin{align*}
    q_{\perp_i}^{(s)}(\infty)-q_{\perp_i}^{(c)}(\infty)-\frac{z(\infty)}{n} =&\ q_{i}^{(s)}(\infty) - q_{\parallel_i}^{(s)}(\infty) - q_{i}^{(c)}(\infty) + q_{\parallel_i}^{(c)}(\infty) - \left(q_{\parallel_i}^{(c)}(\infty) - q_{\parallel_i}^{(s)}(\infty) \right) \\
    =&\
    q_{i}^{(s)}(\infty) - q_{i}^{(c)}(\infty),
\end{align*}
and hence the term $(\dagger)$ can be expressed as
\begin{align*}
    (\dagger )&=\sum_{i=1}^n\left(q_i^{(s)}(\infty)-q_i^{(c)}(\infty)\right)\left(\mu_i^*\mathbbm{1}_{\{q_i^{(s)}(\infty)=q_{\max}^\eta\}}-\lambda_i^*\mathbbm{1}_{\{q_i^{(c)}(\infty)=q_{\max}^\eta\}}\right) \\
    &{=}\sum_{i=1}^n\left(\mu_i^{*}q_i^{(s)}(\infty)\mathbbm{1}_{\{q_i^{(s)}(\infty)=q_{\max}^\eta\}}+\lambda_i^*q_i^{(c)}(\infty)\mathbbm{1}_{\{q_i^{(c)}(\infty)=q_{\max}^\eta\}}\right) \\
    &=q_{\max}^\eta\sum_{i=1}^n\left(\mu_i^{*}\mathbbm{1}_{\{q_i^{(s)}(\infty)=q_{\max}^\eta\}}+\lambda_i^*\mathbbm{1}_{\{q_i^{(c)}(\infty)=q_{\max}^\eta\}}\right),
\end{align*}
where the second equality above follows from Condition \ref{condition: perfect_matching} and the fact that $q_i^{(s)}(\infty)\cdot q_i^{(c)}(\infty)=0$.
Furthermore, $q_i^{(s)}(\infty)=q_{\max}^\eta$ implies that $q_{\perp_i}^{(s)}(\infty) \geq 0$, and $q_i^{(c)}(\infty)=q_{\max}^\eta$ implies that $q_{\perp_i}^{(s)}(\infty) \geq 0$. Therefore, the term $(\ddagger)$ is nonnegative. 
By the Cauchy-Schwarz inequality, we have
\begin{align*}
   & \quad \E\left[\sum_{i=1}^n \left(q_{\perp_i}^{(s)}(\infty)-q_{\perp_i}^{(c)}(\infty)\right)\left(\mu_i^*\mathbbm{1}_{\{q_i^{(s)}(\infty)=q_{\max}^\eta\}}-\lambda_i^*\mathbbm{1}_{\{q_i^{(c)}(\infty)=q_{\max}^\eta\}}\right)\right] \\
    &\leq \sqrt{\E\left[\|\BFq_{\perp}^{(s)}(\infty)-\BFq_{\perp}^{(c)}(\infty)\|^2\right]}\sqrt{\E\left[\sum_{i=1}^n \left(\mu_i^*\mathbbm{1}_{\{q_i^{(s)}(\infty)=q_{\max}^\eta\}}-\lambda_i^*\mathbbm{1}_{\{q_i^{(c)}(\infty)=q_{\max}^\eta\}}\right)^2\right]}  \\
    &\leq \sqrt{\mathcal{F}_2} \sqrt{\sum_{i=1}^n\left( (\mu_i^*)^2 \Pr[q_i^{(s)}(\infty)=q_{\max}^\eta]+(\lambda_j^*)^2 \Pr[q_j^{(c)}(\infty)=q_{\max}^\eta]\right)}  \\
    &{\leq} \sqrt{\mathcal{F}_2} \sqrt{\sum_{i=1}^n \lambda_i^* \max_{i \in {N}} \lambda_i^* +\sum_{i=1}^n \mu_i^* \max_{i \in {N}} \mu_i^*} \frac{1}{\sqrt{q_{\max}^\eta}},
\end{align*}
where the second inequality follows by Lemma~\ref{lemma: SSC_fluid} (which requires the matching algorithm to be max-weight) and the fact that $q_i^{(s)}(\infty)\cdot q_i^{(c)}(\infty)=0$, and the third inequality follows by \eqref{proof1}, as we have
\begin{align*}
    \sum_{i=1}^n \lambda_i^* \Pr[q_i^{(c)}(\infty)=q_{\max}^\eta] \leq \frac{1}{q_{\max}^\eta}\sum_{i=1}^n \lambda_i^* \ \Rightarrow\  \sum_{i=1}^n (\lambda_i^*)^2 \Pr[q_i^{(c)}(\infty)=q_{\max}^\eta]\leq \frac{1}{q_{\max}^\eta}\sum_{i=1}^n \lambda_i^* \max_{i \in {N}} \lambda_i^*,
\end{align*}
and similarly for $\Pr[q_i^{(s)}(\infty)=q_{\max}^\eta]$.
Substituting the above inequality in Eq~\eqref{eq: for-lower-bound-fluid}, we get 
\begin{align}
  & \quad q_{\max}^\eta\left(\sum_{j=1}^n \lambda_j^* \Pr\left[q_j^{(c)}(\infty)=q_{\max}^\eta\right]+\sum_{i=1}^n \mu_i^* \Pr\left[q_i^{(s)}(\infty)=q_{\max}^\eta\right]\right)\nonumber \\
   &\leq  \frac{\inner{\BFone_n}{\BFlambda^*}+\inner{\BFone_n}{\BFmu^*}}{2n}+\sqrt{\mathcal{F}_2} \sqrt{\sum_{i=1}^n \lambda_i^* \max_{i \in {N}} \lambda_i^* +\sum_{i=1}^n \mu_i^* \max_{i \in {N}} \mu_i^*} \frac{1}{\sqrt{q_{\max}^\eta}}. \label{eq: max-weight-bound}
\end{align}
Because $q_{\max}^\eta \rightarrow \infty$ as $\eta \rightarrow \infty$, the second term in Eq~\eqref{eq: max-weight-bound} vanishes. Combined with Eq~\eqref{eq:any-matching-lower-bound}, we obtain the second part of the proposition.
\Halmos
\endproof
\subsection{Proposition \ref{prop: max-weight-optimal-two-price-SSC}}
Similar to Lemma~\ref{lemma: SSC_fluid}, we establish the state space collapse for the two-price policy.
\begin{lemma}[State Space Collapse]\label{lemma: SSC_two-price}
Under the two-price policy and max-weight matching, if 
$\tau_{\max}^\eta \sigma^\eta \leq \eta$, then
for all $\eta \geq 1$, $r \in \mathbb{Z}_+$ there exists constants $\mathcal{T}_r$ independent of $\eta$ such that 
\begin{align*}
    \E\left[\left(\inner{\BFone_{2n}}{\BFq^2(\infty)}-\frac{z^2(\infty)}{n}\right)^{r/2}\right] \leq \mathcal{T}_r.
\end{align*}
\end{lemma}
\proof{Proof of Lemma~\ref{lemma: SSC_two-price}.}
First note that the state space of the Markov chain induced by the two-pricing and max-weight matching policy is $S_{\text{mw}}$ as max-weight matches all compatible customer-server pairs. By Lemma \ref{lemma: bounded_drift} we already know that the drift of the Lyapunov function $U(\BFq)$ is uniformly bounded. Now, we will show that the drift of $U(\BFq)$ is negative outside a finite set. The drift of $V(\BFq)$ for the uniformized DTMC is (by Eq \eqref{eq: SSC_MW}) 
\begin{align*}
   & \quad c \E[\Delta V^{(c)}(\BFq)\mid \BFq(k)=\BFq]\\
   &=\inner{\BFone_m}{\olam}+\inner{\BFone_n}{\omu}+2\sum_{i,j \in E}\chi^{*}_{ij} \left( q_j^{(c)}\mathbbm{1}_{\left\{\max_{i':(i',j)\in E}q_{i'}^{(s)}=0\right\}}-\max_{j':(i,j')\in E}q_{j'}^{(c)} \right)-2\frac{\sigma^{\eta}}{\eta}\sum_{j=1}^n  q_j^{(c)}\mathbbm{1}_{\left\{q_j^{(c)} >\tau_{\max}^{\eta}\right\}} \\
   &\leq \inner{\BFone_n}{\BFlambda^*}+\inner{\BFone_n}{\BFmu^*}+2n {\sigma^\eta \tau^\eta_{\max}}/{\eta}+2 \sum_{(i,j) \in E} \chi_{ij}^* \left(q_j^{(c)}-\max_{j' : (i,j') \in E} q_{j'}^{(c)}\right)-2\frac{\sigma^\eta}{\eta}\inner{\BFone_n}{\BFq^{(c)}} \\
   &\leq \inner{\BFone_n}{\BFlambda^*}+\inner{\BFone_n}{\BFmu^*}+2n +2 \sum_{(i,j) \in E} \chi_{ij}^* \left(q_j^{(c)}-\max_{j' : (i,j') \in E} q_{j'}^{(c)}\right)-2\frac{\sigma^\eta}{\eta}\inner{\BFone_n}{\BFq^{(c)}}, 
\end{align*}
where 
the last inequality uses the assumption that $\tau_{\max}^\eta \sigma^\eta \leq \eta$.
Similarly, by Eq~\eqref{eq: driftv1}, we have
\begin{align*}
    c\E[\Delta V^{(s)}(\BFq)\mid \BFq(k)=\BFq]\leq \inner{\BFone_n}{\BFlambda^*}+\inner{\BFone_n}{\BFmu^*}+2n+2 \sum_{(i,j) \in E} \chi_{ij}^* \left(q_i^{(s)}-\max_{i' : (i',j) \in E} q_{i'}^{(s)}\right)-2\frac{\sigma^\eta}{\eta}\inner{\BFone_n}{\BFq^{(s)}}.
\end{align*}
Combining the two inequalities above, since $\Delta V(\BFq) = \Delta V^{(s)}(\BFq) + \Delta V^{(c)}(\BFq)$, we have
\begin{align}
    c\E[\Delta V(\BFq)\mid \BFq(k)=\BFq] \leq & 2\inner{\BFone_n}{\BFlambda^*}+2\inner{\BFone_n}{\BFmu^*}+4n  \nonumber\\
  & +2 \sum_{(i,j) \in E} \chi_{ij}^* \left(q_i^{(s)}+q_j^{(c)}-\max_{i' : (i',j) \in E} q_{i'}^{(s)}-\max_{j' : (i,j') \in E} q_{j'}^{(c)}\right)-2\frac{\sigma^\eta}{\eta}\inner{\BFone_{2n}}{\BFq}. \label{eq: drift_q2}
\end{align}
Next, the drift of $W_z(\BFq)$ is bounded by
\begin{align}
    \lefteqn{c\E\left[\Delta W_z(\BFq)\mid \BFq(k)=\BFq\right]}\nonumber\\
    ={}&\frac{1}{n}\E\left[z^2(k+1)-z^2(k)\mid \BFq(k)=\BFq\right] \nonumber \\
    ={}&\frac{1}{n}\left(\sum_{j=1}^n \left(\lambda_j^*-\frac{\sigma^\eta}{\eta}\mathbbm{1}_{\{q_j^{(c)}>\tau_{\max}^\eta\}}\right)\left((z+1)^2-z^2\right)+\sum_{i=1}^n \left(\mu_i^*-\frac{\sigma^\eta}{\eta}\mathbbm{1}_{\{q_i^{(s)}>\tau_{\max}^\eta\}}\right)\left((z-1)^2-z^2\right)\right) \nonumber \\
    ={}&\frac{1}{n}\left(\inner{\BFone_n}{\BFlambda^*}+\inner{\BFone_n}{\BFmu^*}-\frac{\sigma^\eta}{\eta}\left(\sum_{j=1}^n \mathbbm{1}_{\{q_j^{(c)}>\tau_{\max}^\eta\}}+\sum_{i=1}^n \mathbbm{1}_{\{q_i^{(s)}>\tau_{\max}^\eta\}}\right)\right) \nonumber \\
    &+\frac{2\sigma^\eta z}{n\eta}\left(\sum_{i=1}^n \mathbbm{1}_{\{q_i^{(s)}>\tau_{\max}^\eta\}}-\sum_{j=1}^n \mathbbm{1}_{\{q_j^{(c)}>\tau_{\max}^\eta\}}\right)  \label{eq: drift_z}\\
    \geq{}& \frac{1}{n}\left(\inner{\BFone_n}{\BFlambda^*}+\inner{\BFone_n}{\BFmu^*}\right)-\frac{2\sigma^\eta}{\eta} -\frac{2\sigma^\eta }{\eta}\inner{\BFone_{2n}}{\BFq}. \nonumber
\end{align}
where the third equality holds because $\sum_{i=1}^n \mu^*_i = \sum_{j=1}^n \lambda^*_j$ and the last inequality holds because $|z|\leq \inner{\BFone_{2n}}{\BFq}$.
Let $B=2(\inner{\BFone_n}{\BFlambda^*}+\inner{\BFone_n}{\BFmu^*})+4n+2$. 
Combining Eq \eqref{eq: drift_q2} and Eq~\eqref{eq: drift_z},
for any 
$U(\BFq) \geq {Bn}/({2}\min_{(i,j) \in E} \chi_{ij}^*)$,
we have
\begin{align*}
    \E[\Delta U(\BFq)\mid \BFq(k)=\BFq] &\leq \frac{\E[\Delta V(\BFq)-\Delta W_z(\BFq)]}{2U(\BFq)} \\
    &\leq \frac{1}{2cU(\BFq)}\left(B+2\sum_{(i,j) \in E} \chi_{ij}^* \left(q_i^{(s)}+q_j^{(c)}-\max_{i' : (i',j) \in E} q_{i'}^{(s)}-\max_{j' : (i,j') \in E} q_{j'}^{(c)}\right)\right), \\
    &{\leq} \frac{B}{2cU(\BFq)}-\frac{2}{cn}\min_{(i,j) \in E} \chi_{ij}^* \leq -\frac{1}{cn}\min_{(i,j) \in E} \chi_{ij}^*.
\end{align*}
where 
the first inequality is by Lemma \ref{lemma: bounded_drift} and
the third inequality is by Lemma~\ref{lemma: imp_term_SSC}.
This proves that
the drift is negative outside a finite set.
By \citet[Lemma 1]{atilla_srikant}, we have
\begin{align*}
    \E\left[U(\BFq(\infty))^r\right] = \E\left[\left(\inner{\BFone_{2n}}{\BFq^2(\infty)}-\frac{z^2(\infty)}{n}\right)^{r/2}\right] \leq \mathcal{T}_r, \quad \forall \ \eta \geq 1. \Halmos 
\end{align*}
\endproof
\begin{lemma} \label{lemma: finite_expectation}
Under the two-price policy and max-weight matching policy, for all $\eta \geq 1$, we have $\E[z^2(\infty)]<\infty$.
\end{lemma}
\proof{Proof of Lemma~\ref{lemma: finite_expectation}.}
We use $\sqrt{V(\BFq)} = \|\BFq\|$ as the Lyapunov function for this proof. The one-step drift is
\begin{align*}
    |\Delta \sqrt{V(\BFq)}| = \big| \| \BFq(k+1)\|-\|\BFq(k)\| \big| \mathbbm{1}_{\{\BFq(k)=\BFq\}} \leq \|\BFq(k+1)-\BFq(k)\| \mathbbm{1}_{\{\BFq(k)=\BFq\}} \leq 1.
\end{align*}
For any $\|\BFq\|\geq {\eta(2\inner{\BFone_n}{\BFlambda^*}+2\inner{\BFone_n}{\BFmu^*}+4n)}/{\sigma^\eta}$, we have
\begin{align*}
   \lefteqn{\quad \E\left[\Delta \sqrt{V(\BFq)}\mid \BFq(k)=\BFq\right]}\\
   &=\E\left[\|\BFq(k+1)\|-\|\BFq(k)\| \big| \BFq(k)=\BFq\right] \\
    &\leq \frac{1}{2\|\BFq\|}\E\left[\|\BFq(k+1)\|^2-\|\BFq(k)\|^2 \big| \BFq(k)=\BFq\right]\\
    & {=} \frac{1}{2c\|\BFq\|}\Bigl(2\inner{\BFone_n}{\BFlambda^*}+2\inner{\BFone_n}{\BFmu^*}+4n\\
    & \qquad \quad+2\sum_{(i,j) \in E} \chi_{ij}^* \left(q_i^{(s)}+q_j^{(c)}-\max_{i' : (i',j) \in E} q_{i'}^{(s)}-\max_{j' : (i,j') \in E} q_{j'}^{(c)}\right)-2\frac{\sigma^\eta}{\eta}\inner{\BFone_{2n}}{\BFq}\Bigr) \\
    &\leq \frac{1}{2c\|\BFq\|}\left(2\inner{\BFone_n}{\BFlambda^*}+2\inner{\BFone_n}{\BFmu^*}+4n-2\frac{\sigma^\eta}{\eta} \inner{\BFone_{2n}}{\BFq}\right)\\
    & \leq \frac{\inner{\BFone_n}{\BFlambda^*}+\inner{\BFone_n}{\BFmu^*}+2n}{c\|\BFq\|}-\frac{\sigma^\eta}{c\eta}<-\frac{\sigma^\eta}{2c\eta},
\end{align*}
where the second equality follows from Eq~\eqref{eq: drift_q2}.
So the drift is negative outside a finite set.
By \citet[Lemma 1]{atilla_srikant}, we have $\E[\|\BFq\|^2]<\infty$ and hence $\E[z^2]<\infty$. \Halmos
\endproof
\proof{Proof of Proposition \ref{prop: max-weight-optimal-two-price-SSC}.}
In the steady state, it holds that $\E[\Delta W_z(\BFq(\infty))] = 0$.
By \eqref{eq: drift_z}, we expand this equation to get
\begin{align}
& \quad 2\sigma^\eta\E\left[z(\infty)\left(\sum_{j=1}^n \mathbbm{1}_{\{q_j^{(c)}(\infty)>\tau_{\max}^\eta\}}-\sum_{i=1}^n \mathbbm{1}_{\{q_i^{(s)}(\infty)>\tau_{\max}^\eta\}}\right)\right] \nonumber \\
    &= \eta\inner{\BFone_n}{\BFlambda^*}+\eta\inner{\BFone_n}{\BFmu^*}-\sigma^\eta\left(\sum_{j=1}^n \Pr\left[q_j^{(c)}(\infty)>\tau_{\max}^\eta\right]+\sum_{i=1}^n \Pr\left[q_i^{(s)}(\infty)>\tau_{\max}^\eta\right]\right) \label{eq: tight_queue_length} \\
    & \geq \eta\inner{\BFone_n}{\BFlambda^*}+\eta\inner{\BFone_n}{\BFmu^*}-2\sigma^\eta n. \nonumber
\end{align}
Dividing both sides by $2\sigma^\eta$, we get
\begin{align}
    \E\left[z(\infty)\left(\sum_{j=1}^n \mathbbm{1}_{\{q_j^{(c)}(\infty)>\tau_{\max}^\eta\}}-\sum_{i=1}^n \mathbbm{1}_{\{q_i^{(s)}(\infty)>\tau_{\max}^\eta\}}\right)\right]\geq \frac{\eta}{2\sigma^\eta}\left(\inner{\BFone_n}{\BFlambda^*}+\inner{\BFone_n}{\BFmu^*}\right)-n.
    \label{eq: two-price-delay-lower-bound}
\end{align}
Recall that $z(\infty) = \sum_{j=1}^n q^{(c)}_j(\infty) - \sum_{i=1}^n q^{(s)}_i(\infty)$, so 
\[
z(\infty)\left(\sum_{j=1}^n \mathbbm{1}_{\{q_j^{(c)}(\infty)>\tau_{\max}^\eta\}}-\sum_{i=1}^n \mathbbm{1}_{\{q_i^{(s)}(\infty)>\tau_{\max}^\eta\}}\right)
\leq n \sum_{j=1}^n q^{(c)}_j(\infty) + n \sum_{i=1}^n q^{(s)}_i(\infty) = n \inner{\BFone_{2n}}{\BFq(\infty)}.
\]
Substituting the above inequality in Eq~\eqref{eq: two-price-delay-lower-bound} and dividing both sides by $n$, we get
\begin{equation}
    \E[\inner{\BFone_{2n}}{\BFq(\infty)}] \geq \frac{\eta}{2\sigma^\eta n}\left(\inner{\BFone_n}{\BFlambda^*}+\inner{\BFone_n}{\BFmu^*}\right)-1.
    \label{eq: two-price-delay-lower-bound-2}
\end{equation}
This proves the first part of the proposition. Now, under the max-weight matching algorithm, by Lemma \ref{lemma: finite_expectation}, we have $\E[W_z(\BFq)]<\infty$. So Eq~\eqref{eq: two-price-delay-lower-bound-2} holds for max-weight matching algorithm. 
Next, we consider the max-weight matching algorithm specifically.
Using \eqref{eq: tight_queue_length} again,
we have
\begin{align}
    {}& \E\left[z(\infty)\left(\sum_{j=1}^n \mathbbm{1}_{\{q_j^{(c)}(\infty)>\tau_{\max}^\eta\}}-\sum_{i=1}^n \mathbbm{1}_{\{q_i^{(s)}(\infty)>\tau_{\max}^\eta\}}\right)\right]\leq \frac{\eta}{2\sigma^\eta}\left(\inner{\BFone_n}{\BFlambda^*}+\inner{\BFone_n}{\BFmu^*}\right)+\frac{n}{2}. 
    \nonumber
\end{align}
Dividing both sides by $n$ and splitting the left-hand side, we get
\begin{align}
    & \E\left[\sum_{i=1}^n\left(\frac{z(\infty)}{n}+q_i^{(s)}(\infty)-q_i^{(c)}(\infty)\right)\left(\mathbbm{1}_{\{q_i^{(c)}(\infty)>\tau_{\max}^\eta\}}-\mathbbm{1}_{\{q_i^{(s)}(\infty)>\tau_{\max}^\eta\}}\right)\right] \nonumber \\
    & + \E\left[\sum_{i=1}^n\left(q_i^{(s)}(\infty)-q_i^{(c)}(\infty)\right)\left(\mathbbm{1}_{\{q_i^{(s)}(\infty)>\tau_{\max}^\eta\}}-\mathbbm{1}_{\{q_i^{(c)}(\infty)>\tau_{\max}^\eta\}}\right)\right] 
    \leq \frac{\eta}{2n\sigma^\eta}\left(\inner{\BFone_n}{\BFlambda^*}+\inner{\BFone_n}{\BFmu^*}\right)+\frac{1}{2}. \label{eq: two-price-delay-opt-upper-bound}
\end{align}
We bound the two expectation term separately.
The first term in Eq~\eqref{eq: two-price-delay-opt-upper-bound} can be bounded by
\begin{align*}
    &\E \left[\sum_{i=1}^n \left(\frac{z(\infty)}{n}+q_i^{(s)}(\infty)-q_i^{(c)}(\infty)\right)\left(\mathbbm{1}_{\left\{q_i^{(s)}>\tau_{\max}^\eta\right\}}-\mathbbm{1}_{\left\{q_i^{(c)}>\tau_{\max}^\eta\right\}}\right)\right] \\
    \geq{}& -\E\left[\sum_{i=1}^n \left|q_i^{(s)}(\infty)-q_i^{(c)}(\infty)+\frac{z(\infty)}{n}\right|\right] \\
    ={}& -\E\left[\|\BFq_{\perp}^{(s)}(\infty)-\BFq_{\perp}^{(c)}(\infty)\|_1\right] \\
    \geq{}& -\sqrt{n} \E \left[\|\BFq_{\perp}^{(s)}(\infty)-\BFq_{\perp}^{(c)}(\infty)\|_2\right] \\
    \geq{}& - \sqrt{n \E\left[\|\BFq_{\perp}^{(s)}(\infty)-\BFq_{\perp}^{(c)}(\infty)\|_2^2\right]}  \\
    ={}& - \sqrt{n\E\left[\inner{\BFone_{2n}}{\BFq^2}-\frac{z^2}{n}\right]} \\
    \geq{}& - \sqrt{n\mathcal{T}_2}, 
\end{align*}
where the first equality uses the definition of $\BFq_{\perp}$, the second inequality uses the fact that $\|\BFx\|_1 \leq \sqrt{n}\|\BFx\|_2$ for any $\BFx \in \mathbb{R}^n$,
the second inequality uses Lemma~\ref{lemma: perp_representation}, and the final inequality uses Lemma~\ref{lemma: SSC_two-price} (which requires the matching algorithm to be max-weight) with $r=2$.
The second term on the left-hand side of Eq~\eqref{eq: two-price-delay-opt-upper-bound} is equal to
\begin{align*}
    &\sum_{i=1}^n \E\left[ (q_i^{(s)}(\infty)-q_i^{(c)}(\infty))\left(\mathbbm{1}_{\{q_i^{(s)}(\infty)>\tau_{\max}^\eta\}}-\mathbbm{1}_{\{q_j^{(c)}(\infty)>\tau_{\max}^\eta\}}\right) \right]\\
    ={}& \sum_{i=1}^n\E\left[ q_i^{(s)}(\infty)\mathbbm{1}_{\{q_i^{(s)}(\infty)>\tau_{\max}^\eta\}}+q_i^{(c)}(\infty)\mathbbm{1}_{\{q_j^{(c)}(\infty)>\tau_{\max}^\eta\}}-q_i^{(s)}(\infty)\mathbbm{1}_{\{q_j^{(c)}(\infty)>\tau_{\max}^\eta\}}-q_i^{(c)}(\infty)\mathbbm{1}_{\{q_i^{(s)}(\infty)>\tau_{\max}^\eta\}} \right]\\
    ={}&\sum_{i=1}^n\E\left[ q_i^{(s)}(\infty)\mathbbm{1}_{\{q_i^{(s)}(\infty)>\tau_{\max}^\eta\}}+q_i^{(c)}(\infty)\mathbbm{1}_{\{q_j^{(c)}(\infty)>\tau_{\max}^\eta\}}\right]\\
    ={}& \sum_{i=1}^n \E \left[ q_i^{(s)}(\infty) +q_i^{(c)}(\infty) \right]-\sum_{i=1}^n \E\left[  q_i^{(s)}(\infty)\mathbbm{1}_{\{q_i^{(s)}(\infty)\leq\tau_{\max}^\eta\}}\right]
    - \sum_{i=1}^n \E\left[ q_i^{(c)}(\infty)\mathbbm{1}_{\{q_i^{(c)}(\infty)\leq\tau_{\max}^\eta\}}\right]\\
    \geq{}&\sum_{i=1}^n \E\left[q_i^{(s)}(\infty) + q_i^{(c)}(\infty)\right]- 2n \tau_{\max}^\eta\\
     ={}&\E\left[\inner{\BFone_{2n}}{\BFq(\infty)}\right] -2n \tau_{\max}^\eta.
\end{align*}
Note the the last two terms in the first equality are equal to zero because $q_i^{(s)}(\infty)q_i^{(c)}(\infty)=0$ for all $i \in [n]$. 
Combining the two inequalities above and using Eq~\eqref{eq: two-price-delay-opt-upper-bound}, we have
\begin{align}
& \E\left[\inner{\BFone_{2n}}{\BFq(\infty)}\right] \leq \frac{\eta}{2n\sigma^\eta}\left(\inner{\BFone_n}{\BFlambda^*}+\inner{\BFone_n}{\BFmu^*}\right)+\frac{1}{2}+\sqrt{n\mathcal{T}_2}+2n\tau_{\max}^\eta, \label{eq: queue_upper_bound_tight}
\end{align}
Combining the upper bound Eq~\eqref{eq: queue_upper_bound_tight} with the lower bound Eq~\eqref{eq: two-price-delay-lower-bound-2}, under the max-weight algorithm, we have
\begin{align*}
    \frac{\eta}{2n\sigma^\eta}\left(\inner{\BFone_n}{\BFlambda^*}+\inner{\BFone_n}{\BFmu^*}\right)-\frac{1}{2} \leq \E[\inner{\BFone_{2n}}{\BFq(\infty)}] \leq \frac{\eta}{2n\sigma^\eta}\left(\inner{\BFone_n}{\BFlambda^*}+\inner{\BFone_n}{\BFmu^*}\right)+\frac{1}{2}+ \sqrt{n\mathcal{T}_2}+2n\tau_{\max}^\eta.
\end{align*}
As $\eta \rightarrow \infty$, by assumption, we have
$\sigma^\eta/\eta \rightarrow 0$ and $\sigma^\eta\tau_{\max}^\eta/\eta \rightarrow 0$; therefore,
\begin{align*}
    \lim_{\eta \rightarrow \infty}\frac{\sigma^\eta}{\eta}\E[\inner{\BFone_{2n}}{\BFq(\infty)}] =  \frac{1}{2n}\left(\inner{\BFone_n}{\BFlambda^*}+\inner{\BFone_n}{\BFmu^*}\right).\Halmos
\end{align*}
\endproof
\subsection{Theorem \ref{theo: fluid_max_weight_SSC}}
\begin{lemma} \label{lemma: profit-loss-max-weight-fluid}
Under the fluid pricing policy with max-weight matching, for $q_{\max}^\eta=\sqrt{\eta/n}$, we have
\begin{align*}
    \limsup_{\eta \rightarrow \infty}\frac{\lr}{\eta^{1/2}} \leq n^{1/2}\left(\frac{\inner{\BFone_n}{\BFlambda^*}+\inner{\BFone_n}{\BFmu^*}}{2n} \max_{j \in {N}} F_j(\lambda_j^*)+2\max_{i \in N, j\in M}\{s_i^{(s)},s_j^{(c)}\}\right).
\end{align*}
\end{lemma}
\proof{Proof of Lemma~\ref{lemma: profit-loss-max-weight-fluid}.}
Recall that the loss $L^\eta$ \eqref{eq: rl} is composed of the holding cost term and the revenue loss term.
The holding cost term depends on the queue length vector $\BFq$. According to the definition of the fluid pricing policy, we have $\BFq \leq q_{\max}^{\eta}\BFone_{2n}\ a.s$. Thus, it is trivially true that $\inner{\BFs}{\E[\BFq(\infty)]} \leq 2\max_{i \in N, j\in M}\{s_i^{(s)},s_j^{(c)}\}nq_{\max}^{\eta}$. 
For the revenue loss term, we apply the upper bound in Eq~\eqref{eq: max-weight-bound} through the inequality
\begin{align*}
    &\quad \eta\inner{F(\olam)}{(\olam\circ\E[\mathbf{I}^{(c)}(q_{\max}^{\eta})])}-\eta\inner{G(\omu)}{(\omu\circ\E[\mathbf{I}^{(s)}(q_{\max}^{\eta})]} \\
    & \leq \eta \max_{j \in {N}}  F_j(\lambda_j^*) \left(\sum_{j=1}^n \lambda_j^* \Pr\left[q_j^{(c)}(\infty)=q_{\max}^\eta\right]+ \sum_{i=1}^n \mu_i^* \Pr\left[q_i^{(s)}(\infty)=q_{\max}^\eta\right]\right).
\end{align*}
Therefore, the loss $\lr$ \eqref{eq: rl} is bounded by
\begin{align*}
    \frac{\lr}{\eta^{1/2}} \leq{}& \max_{j \in {N}} F_j(\lambda_j^*)\left( \eta^{1/2}\frac{\inner{\BFone_n}{\BFlambda^*}+\inner{\BFone_n}{\BFmu^*}}{2nq_{\max}^\eta}+ \sqrt{\mathcal{F}_2} \sqrt{\sum_{i=1}^n \lambda_i^* \max_{i \in {N}} \lambda_i^* +\sum_{i=1}^n \mu_i^* \max_{i \in {N}} \mu_i^*} \frac{\eta^{1/2}}{(q_{\max}^\eta)^{3/2}}\right)\\
    &+2\max_{i \in N, j\in M}\{s_i^{(s)},s_j^{(c)}\}\frac{nq_{\max}^{\eta}}{\eta^{1/2}},
\end{align*}
where the constant $\mathcal{F}_2$ is defined in Lemma~\ref{lemma: SSC_fluid} and is independent of $\eta$.
Setting $q_{\max}^\eta=\sqrt{\eta/n}$, we note that the term involving $\mathcal{F}_2$ converges to 0 as $\eta \to\infty$.
Therefore, we get
\begin{align}
    \limsup_{\eta \rightarrow \infty}\frac{\lr}{\eta^{1/2}} \leq{}& \frac{1}{n^{1/2}}\frac{\inner{\BFone_n}{\BFlambda^*}+\inner{\BFone_n}{\BFmu^*}}{2} \max_{j \in {N}} F_j(\lambda_j^*)+2n^{1/2}\max_{i \in N, j\in M}\{s_i^{(s)},s_j^{(c)}\}\nonumber\\
    ={}&n^{1/2}\left(\frac{\inner{\BFone_n}{\BFlambda^*}+\inner{\BFone_n}{\BFmu^*}}{2n} \max_{j \in {N}} F_j(\lambda_j^*)+2\max_{i \in N, j\in M}\{s_i^{(s)},s_j^{(c)}\}\right). \nonumber \Halmos
\end{align}
\endproof
\begin{lemma} \label{lemma: profit-loss-random-fluid}
Under the fluid pricing policy with randomized matching (Algorithm \ref{alg: random_non_empty}), for $q_{\max}^\eta=\gamma \sqrt{\eta}$, we have
\begin{align*}
 \lr \leq \left(\frac{\sum_{i=1}^n \mu_i^* + \sum_{j=1}^m \lambda_j^* }{2\gamma} \max_{j \in {M}} \{F_j(\lambda_j^*)\}+\inner{\BFone_{n+m}}{\BFs}\gamma\right)\eta^{1/2}.
\end{align*}
\end{lemma}
\proof{Proof of Lemma~\ref{lemma: profit-loss-random-fluid}.}
We use $V^{(c)}(\BFq) = \inner{\BFone_{m}}{(\BFq^{(c)})^2}$ as the Lyapunov function. By Eq \eqref{eq: drift_q2_fluid}, the one-step drift of $V^{(c)}(\BFq)$ is bounded by
\begin{align*}
    \eta c\E[\Delta V^{(c)}(\BFq)]&\leq \sum_{i=1}^n \etam_i + \sum_{j=1}^m \etal_j+2 \E \left[\sum_{j=1}^m\etal_j q_j^{(c)}\left(1-\sum_{i=1}^n y_{ij}^{}\right)\mathbbm{1}_{\left\{q_j^{(c)}< q_{\max}^{\eta}\right\}}-\sum_{i=1}^n\etam_i\sum_{j=1}^mq_j^{(c)}y_{ij}^{}\right] \\
    &\leq \sum_{i=1}^n \etam_i + \sum_{j=1}^m \etal_j+2 \E \left[\sum_{j=1}^m\etal_j q_j^{(c)}\mathbbm{1}_{\left\{q_j^{(c)}< q_{\max}^{\eta}\right\}}-\sum_{i=1}^n\etam_i\sum_{j=1}^mq_j^{(c)}y_{ij}^{}\right] \\
    &= \sum_{i=1}^n \etam_i + \sum_{j=1}^m \etal_j+2 \E \left[\sum_{j=1}^m\etal_j q_j^{(c)} -\sum_{i=1}^n\etam_i\sum_{j=1}^mq_j^{(c)}y_{ij}^{} - \sum_{j=1}^m\etal_j q_{\max}^\eta I_j^{(c)}(q_{\max}^\eta)\right] \\
    &{\leq}\sum_{i=1}^n \etam_i + \sum_{j=1}^m \etal_j+2\E \left[\sum_{j=1}^m\etal_j q_j^{(c)}  -\sum_{i=1}^n\etam_i\sum_{j=1}^m q_j^{(c)}\frac{\chi^*_{ij}}{\mu_i^*} - \sum_{j=1}^m\etal_j q_{\max}^\eta I_j^{(c)}(q_{\max}^\eta)\right] \\
    &=\sum_{i=1}^n \etam_i + \sum_{j=1}^m \etal_j-2\E\left[\sum_{j=1}^m\etal_j q_{\max}^\eta I_j^{(c)}(q_{\max}^\eta)\right].
\end{align*}
The third inequality above follows from the definition of the randomized matching algorithm (Algorithm~\ref{alg: random_non_empty}): when $q_j^{(c)}>0$, if 
a type $i$ server arrives,
the conditional probability that the server is matched to a type $j$ customer is at least $\chi_{ij}^*/\mu_i^*$. The last equality holds because $\lambda^*_j = \sum_{i=1}^n \chi^*_{ij}$ for all $j\in M$ by Eq~\eqref{ball}.
In the steady state, we have $\E[\Delta V^{(c)}(\BFq)]=0$. Rearranging the term in the above inequality, we get
\[
   2\sum_{j=1}^m\etal_j q_{\max}^\eta \mathbb{e}\E\left[I_j^{(c)}(q_{\max}^\eta) \right] \leq  \sum_{i=1}^n \etam_i + \sum_{j=1}^m \etal_j,
\]
which implies that
\[
   \sum_{j=1}^m\etal_j F_j(\lambda_j^*) \E\left[I_j^{(c)}(q_{\max}^\eta) \right] \leq \frac{\sum_{i=1}^n \etam_i + \sum_{j=1}^m \etal_j }{2q_{\max}^\eta} \max_{j \in {M}} \{F_j(\lambda_j^*)\}.
\]
Setting $q_{\max}^\eta=\gamma\eta^{1/2}$, we have
\begin{align*}
   L^{\eta}&=\sum_{j=1}^m\etal_j F_j(\lambda_j^*) \E\left[I_j^{(c)}(q_{\max}^\eta) \right]+\E[\inner{\BFs}{\BFq}]\\
   &\leq \frac{\sum_{i=1}^n \etam_i + \sum_{j=1}^m \etal_j }{2q_{\max}^\eta} \max_{j \in {M}} \{F_j(\lambda_j^*)\}+\inner{\BFone_{n+m}}{\BFs}q_{\max}^\eta \\
   &= \eta^{1/2} \frac{\sum_{i=1}^n \mu_i^* + \sum_{j=1}^m \lambda_j^* }{2\gamma} \max_{j \in {M}} \{F_j(\lambda_j^*)\}+\inner{\BFone_{n+m}}{\BFs}\gamma\eta^{1/2}.
   \Halmos
\end{align*}
\endproof
\proof{Proof of Theorem \ref{theo: fluid_max_weight_SSC}.}
In Lemma \ref{lemma: profit-loss-max-weight-fluid}, we have shown that 
$\lim\sup_{\eta\to\infty} L^{\eta}/\eta^{1/2} = O({n}^{1/2})$
for the max-weight matching algorithm. By Lemma \ref{lemma: profit-loss-random-fluid}, we know that the profit loss is $O(\eta^{1/2})$ for the randomized matching policy.
To complete the proof, we want to show $\lim\inf_{\eta\to\infty} L^{\eta}/\eta^{1/2} = \Omega(n)$ for the randomized matching policy.
\begin{figure}[!htb]
    \FIGURE
    {\begin{tabular}{c}
        \begin{tikzpicture}[scale=0.7]
\draw[black, very thick] (0,0) -- (2,0) -- (2,1) -- (0,1);
\node[black,very thick] at (0.25,0.5) {2};
\draw[black, very thick] (0,1.5) -- (2,1.5) -- (2,2.5) -- (0,2.5);
\node[black,very thick] at (0.25,2) {1};
\fill[black] (1,-0.3) circle (0.05);
\fill[black] (1,-0.5) circle (0.05);
\fill[black] (1,-0.7) circle (0.05);
\draw[black, very thick] (0,-1) -- (2,-1) -- (2,-2) -- (0,-2);
\node[black,very thick] at (0.25,-1.5) {n};
\draw[black,very thick] (8,0) -- (6,0) -- (6,1) -- (8,1);
\node[black,very thick] at (7.75,0.5) {2};
\draw[black, very thick] (8,1.5) -- (6,1.5) -- (6,2.5) -- (8,2.5);
\node[black,very thick] at (7.75,2) {1};
\fill[black] (7,-0.3) circle (0.05);
\fill[black] (7,-0.5) circle (0.05);
\fill[black] (7,-0.7) circle (0.05);
\draw[black, very thick] (8,-1) -- (6,-1) -- (6,-2) -- (8,-2);
\node[black,very thick] at (7.75,-1.5) {n};
\draw[black,thick]  (2.75, 2.1) edge[<->]  (5.25, 2.1);
\draw[black,thick]  (2.75, 1.9) edge[<->]  (5.25, -1.4);
\draw[black,thick]  (2.75, 0.6) edge[<->]  (5.25, 1.9);
\draw[black,thick]  (2.75, 0.4) edge[<->]  (5.25, 0.4);
\draw[dotted,black,thick]  (2.75, -1.2) edge[<->]  (5.25, 0.2);
\fill[black] (4,-0.3) circle (0.05);
\fill[black] (4,-0.5) circle (0.05);
\fill[black] (4,-0.7) circle (0.05);
\draw[black,thick]  (2.75, -1.5) edge[<->]  (5.25, -1.5);
\node[black, align=center] at (1,3) {\footnotesize Servers};
\node[black, align=center] at (4,3) {\footnotesize \shortstack{Compatible\\Matching}};
\node[black, align=center] at (7,3) {\footnotesize Customers};
\end{tikzpicture}
    \end{tabular}}
    {\centering{The ``long chain'' graph used in the proof.}
    \label{fig: long_chain}}
    {}
\end{figure}
To this end, we consider the following instance. Consider a bipartite graph with customer type $j$ is connected with server type $j$ and type $j+1$ modulo $n$ (see Figure~\ref{fig: long_chain}).
We assume that the demand curve is $F_j(x)=2-0.5x$ for all $j \in {N}$ and the supply curve is $G_i(x)=0.5x$ for all $i \in {N}$.
It is easy to verify that $\BFlambda^*=\BFmu^*=\BFone_n$, $\chi_{ii}=1$ for all $i \in {N}$ is an optimal solution to the fluid problem. It is also easy to see that the graph and the fluid solution satisfy both Condition \ref{condition: crp} and Condition \ref{condition: perfect_matching}. As $\chi_{ij}=0$ for all $i \neq j$, the system behaves like $n$ independent two-sided queues under the randomized matching algorithm.
For each single-link two-sided queue (with loss $\lr/n$),
by Eq~\eqref{eq: sl-rl} in Proposition \ref{prop: singlelinktwosidedqueue}, the loss under any fluid pricing policy is lower bounded by
\begin{align*}
    \frac{\lr}{n} & \geq  \frac{1}{2q_{\max}^{\eta}+1}\eta+\min\{s^{(s)},s^{(c)}\}\frac{q_{\max}^{\eta}(q_{\max}^{\eta}+1)}{2q_{\max}^{\eta}+1}.
\end{align*}
To minimize the order of $\eta$ on the right-hand side, we set $q^{\eta}_{\max} = \gamma\eta^{1/2}$.
Therefore, using the AM-GM inequality $(a+b)/2\geq \sqrt{ab}$, the total loss for this system is bounded by
 \begin{align*}  
    \liminf_{\eta \rightarrow \infty}\frac{\lr}{\eta^{1/2}}&\geq n\left(\frac{1}{2\gamma}+\frac{\min\{s^{(s)},s^{(c)}\}\gamma}{2}\right) \geq n \sqrt{\min\{s^{(s)},s^{(c)}\}}. \Halmos
\end{align*}
\endproof
\subsection{Theorem \ref{theo: two_price_maxweight_SSC}}
\begin{lemma} \label{lemma: profit-loss-two-price-max-weight}
Under the two-price and max-weight matching policy, for $\BFtheta=\BFone_n$ and $\BFphi=\BFone_n$, $\sigma^\eta=n^{-1/3}\eta^{2/3}$ and $\tau_{\max}^\eta=o(\eta^{1/3})$, we have
\begin{align*}
\limsup_{\eta \rightarrow \infty} \frac{\lr}{\eta^{1/3}} 
    \leq& \left(\sum_{i \in {N}} \left(\frac{\mu_i^*G_i''(\mu_i^*)}{2}+G_i'(\mu_i^*)\right)-\sum_{j \in {N}} \left(\frac{\lambda_j^*F_j''(\lambda_j^*)}{2}+F_j'(\lambda_j^*)\right)\right.\nonumber\\
    &\left.+\frac{\max_{i \in N,j\in M}\{s_i^{(s)},s_j^{(c)}\}}{2}\left(\inner{\BFone_n}{\BFlambda^*}+\inner{\BFone_n}{\BFmu^*}\right)\right)n^{-2/3}=O(n^{1/3}).
\end{align*}
\end{lemma}
\proof{Proof of Lemma~\ref{lemma: profit-loss-two-price-max-weight}.}
By \eqref{eq: theo2_profit_loss}, the profit loss under the two-price policy is
\begin{align}
    \lr={}&\eta\sum_{j \in {M}} \left(F_j(\lambda_j^{*})\lambda_j^{*}-F_j\left(\lambda_j^{*}-\theta_j\frac{\sigma^{\eta}}{\eta}\right)\left(\lambda_j^{*}-\theta_j\frac{\sigma^{\eta}}{\eta}\right)\right) \Pr[q^{(c)}_j > \tau_{\max}^{\eta}] \nonumber \\
    &-\eta\sum_{i \in {N}} \left(G_i(\mu_i^{*})\mu_i^{*}-G_i\left(\mu_i^{*}-\phi_i\frac{\sigma^{\eta}}{\eta}\right)\left(\mu_i^{*}-\phi_i\frac{\sigma^{\eta}}{\eta}\right)\right) \Pr[q^{(s)}_i > \tau_{\max}^{\eta}]+\inner{\BFs}{\E[\BFq(\infty)]}
    \nonumber\\
    ={}& \sum_{j \in {M}} \left(F_j'(\lambda_j^{*})\lambda_j^{*}+F_j(\lambda_j^{*})\right)\theta_j\sigma^{\eta} \Pr[q^{(c)}_j>\tau_{\max}^{\eta}]
  -\sum_{i \in {N}} \left(G_i'(\mu_i^{*})\mu_i^{*}+G_i(\mu_i^{*})\right)\phi_i\sigma^{\eta} \Pr[q^{(s)}_i > \tau_{\max}^{\eta}]\nonumber \\
  &-\left(\sum_{j \in {N}} \left(\frac{\lambda_j^*F_j''(\lambda_j^*)}{2}+F_j'(\lambda_j^*)\right)\theta_j^2\Pr[q_j^{(c)}>\tau_{\max}^\eta] -\sum_{i \in {N}} \left(\frac{\mu_i^*G_i''(\mu_i^*)}{2}+G_i'(\mu_i^*)\right)\phi_i^2\Pr[q_i^{(s)}>\tau_{\max}^\eta]\right)\frac{(\sigma^\eta)^2}{\eta} \nonumber\\
  &+o\left(\frac{(\sigma^\eta)^2}{\eta}\right)+  \inner{\BFs}{\E[\BFq(\infty)]}, \label{eq: profit-loss-taylor-thm}
\end{align}
where the second equality uses Taylor's theorem.
By Lemma~\ref{lemma: positive_chi}, there exists a fluid optimal solution such that $\chi^*_{ij}>0$ for all $(i,j) \in E$. Therefore, by Claim~\ref{claim}, it holds that
\[
\sum_{j \in {M}} \left(F_j'(\lambda_j^{*})\lambda_j^{*}+F_j(\lambda_j^{*})\right)\theta_j\sigma^{\eta} \Pr[q^{(c)}_j>\tau_{\max}^{\eta}] - \sum_{i \in {N}} \left(G_i'(\mu_i^{*})\mu_i^{*}+G_i(\mu_i^{*})\right)\phi_i\sigma^{\eta} \Pr[q^{(s)}_i > \tau_{\max}^{\eta}] = 0.
\]
By Eq~\eqref{eq: queue_upper_bound_tight}, we have
\[
\inner{\BFs}{\E[\BFq(\infty)]} \leq 
\max_{i \in N,j\in M}\{s_i^{(s)},s_j^{(c)}\}\left(\frac{\eta}{2n\sigma^\eta}\left(\inner{\BFone_n}{\BFlambda^*}+\inner{\BFone_n}{\BFmu^*}\right)+\frac{1}{2}+\sqrt{n\mathcal{T}_2}+2n\tau_{\max}^\eta\right).
\]
Substituting the above two equations in Eq~\eqref{eq: profit-loss-taylor-thm} and setting $\theta_j = \phi_i = 1$, $\sigma^\eta = n^{-1/3}\eta^{2/3}$ and $\tau^{\eta}_{\max} = o(\eta^{1/3})$, we get
\begin{align*}
\lr{\leq}
  {}\ &\  \left(\sum_{i \in {N}} \left(\frac{\mu_i^*G_i''(\mu_i^*)}{2}+G_i'(\mu_i^*)\right)-\sum_{j \in {N}} \left(\frac{\lambda_j^*F_j''(\lambda_j^*)}{2}+F_j'(\lambda_j^*)\right)\right)\eta^{1/3}n^{-2/3}\\
  &+\max_{i \in N,j\in M}\{s_i^{(s)},s_j^{(c)}\}\left(\frac{\eta^{1/3}}{2n^{2/3}}\left(\inner{\BFone_n}{\BFlambda^*}+\inner{\BFone_n}{\BFmu^*}\right)+\frac{1}{2}+ \sqrt{n\mathcal{T}_2}+o(\eta^{1/3})\right)+o(\eta^{1/3}).
\end{align*}
Thus, we have
\begin{align*}
    \limsup_{\eta \rightarrow \infty} \frac{\lr}{\eta^{1/3}} 
    \leq{}& \left(\sum_{i \in {N}} \left(\frac{\mu_i^*G_i''(\mu_i^*)}{2}+G_i'(\mu_i^*)\right)-\sum_{j \in {N}} \left(\frac{\lambda_j^*F_j''(\lambda_j^*)}{2}+F_j'(\lambda_j^*)\right)\right.\nonumber\\
    &\left.+\frac{\max_{i \in N,j\in M}\{s_i^{(s)},s_j^{(c)}\}}{2}\left(\inner{\BFone_n}{\BFlambda^*}+\inner{\BFone_n}{\BFmu^*}\right)\right)n^{-2/3}=O(n^{1/3}). \Halmos
\end{align*}
\endproof
\begin{lemma} \label{lemma: profit-loss-two-price-random}
Under the two-price and randomized matching policy (Algorithm \ref{alg: random_non_empty}), for $\sigma^\eta=\eta^{2/3}$ and $\tau_{\max}^\eta=\gamma \eta^{1/3}$, we have
\begin{align*}
 \lr=O(\eta^{1/3}).
\end{align*}
\end{lemma}
\proof{Proof of Lemma~\ref{lemma: profit-loss-two-price-random}.}
We use $V^{(c)}(\BFq) = \inner{\BFone_{m}}{(\BFq^{(c)})^2}$ as the Lyapunov function. By Eq \eqref{eq: drift_q2_fluid}, the one-step drift of $V^{(c)}(\BFq)$ is bounded by
\begin{align}
     \eta c \E\left[\Delta V^{(c)}(\BFq)\right]={}&\E\Bigg[\sum_{j=1}^m \left(\etal_j-\theta_j\sigma^{\eta}\mathbbm{1}_{\left\{q_j^{(c)}>\tau_{\max}^{\eta}\right\}}\right)\left(1-\sum_{i=1}^n y_{ij}^{}\right)\left(1+2q_j^{(c)}\right) \nonumber \\
    & \; +\sum_{i=1}^n \left(\etam_i-\phi_i\sigma^{\eta}\mathbbm{1}_{\left\{q_i^{(s)}>\tau_{\max}^{\eta}\right\}}\right)\sum_{j=1}^m y_{ij}^{}\left(1-2q_j^{(c)}\right)\Bigg] \nonumber \\
    \leq{}& \E\Bigg[\sum_{j=1}^m \left(\etal_j-\theta_j\sigma^{\eta}\mathbbm{1}_{\left\{q_j^{(c)}>\tau_{\max}^{\eta}\right\}}\right)\left(1+2q_j^{(c)}\right) \nonumber \\
    &\; +\sum_{i=1}^n \left(\etam_i-\phi_i\sigma^{\eta}\mathbbm{1}_{\left\{q_i^{(s)}>\tau_{\max}^{\eta}\right\}}\right)\sum_{j=1}^m y_{ij}^{}\left(1-2q_j^{(c)}\right)\Bigg] \nonumber \\
    \leq{}& \eta\left(\inner{\BFone_{n}}{\BFmu^*}+\inner{\BFone_m}{\BFlambda^*}\right)+2\E\Bigg[\eta\left(\sum_{j=1}^m \lambda_j^* q_j^{(c)}-\sum_{i=1}^n \mu_i^*\sum_{j=1}^m y_{ij}^{}q_j^{(c)}\right)\nonumber \\
    &-\sigma^\eta \sum_{j=1}^m \theta_jq_j^{(c)}\mathbbm{1}_{\left\{q_j^{(c)}>\tau_{\max}^\eta\right\}} 
    +\sigma^\eta \sum_{i=1}^n \phi_i \mathbbm{1}_{\left\{q_i^{(s)}>\tau_{\max}^\eta\right\}}\sum_{j=1}^m y_{ij}^{}q_j^{(c)}\Bigg] \nonumber \\
    \leq{}& \eta\left(\inner{\BFone_{n}}{\BFmu^*}+\inner{\BFone_m}{\BFlambda^*}\right)+2\E\Bigg[\eta\left(\sum_{j=1}^m \lambda_j^* q_j^{(c)}-\sum_{i=1}^n \sum_{j=1}^m \chi_{ij}^{*}q_j^{(c)}\right)\nonumber \\
    &-\sigma^\eta \sum_{j=1}^m \theta_jq_j^{(c)}\mathbbm{1}_{\left\{q_j^{(c)}>\tau_{\max}^\eta\right\}} 
    +\sigma^\eta \sum_{i=1}^n \phi_i \mathbbm{1}_{\left\{q_i^{(s)}>\tau_{\max}^\eta\right\}}\sum_{j=1}^m y_{ij}q_j^{(c)}\Bigg] \nonumber \\
    {=}{}&\eta\left(\inner{\BFone_{n}}{\BFmu^*}+\inner{\BFone_m}{\BFlambda^*}\right)-2\sigma^\eta\E\left[ \sum_{j=1}^m \theta_jq_j^{(c)}\mathbbm{1}_{\left\{q_j^{(c)}>\tau_{\max}^\eta\right\}}- \sum_{i=1}^n \phi_i \mathbbm{1}_{\left\{q_i^{(s)}>\tau_{\max}^\eta\right\}}\sum_{j=1}^m y_{ij}q_j^{(c)}
    \right]\nonumber \\
    {=}{}&\eta\left(\inner{\BFone_{n}}{\BFmu^*}+\inner{\BFone_m}{\BFlambda^*}\right)-2\sigma^\eta \E\left[\sum_{j=1}^m \theta_jq_j^{(c)}\mathbbm{1}_{\left\{q_j^{(c)}>\tau_{\max}^\eta\right\}}\right], \nonumber
\end{align}
where the second inequality uses $\sum_{j=1}^m \lambda_j^* = \sum_{i=1}^n \mu_i^*$. The third inequality
follows by the definition of Algorithm \ref{alg: random_non_empty}: when $q_j^{(c)}>0$, if 
a type $i$ server arrives,
the conditional probability that the server is matched to a type $j$ customer is at least $\chi_{ij}^*/\mu_i^*$. The second equality uses the fact that $\lambda^*_j = \sum_{i=1}^n \chi^*_{ij}$ for all $j\in M$ by Eq~\eqref{ball}.
The last equality follows as for all $(i,j) \in E$ and $\BFq \in S$, we have $q_i^{(s)}q_j^{(c)}=0$.
Similarly, we can calculate the drift of $V^{(s)}(\BFq)$ to get
\begin{align*}
    \eta c\E[\Delta V^{(s)}(\BFq)] \leq \eta\left(\inner{\BFone_{n}}{\BFmu^*}+\inner{\BFone_m}{\BFlambda^*}\right)-2\sigma^\eta \E\left[\sum_{i=1}^n \phi_iq_i^{(s)}\mathbbm{1}_{\left\{q_i^{(s)}>\tau_{\max}^\eta\right\}}\right].
\end{align*}
Combining the two inequalities above, since $V(\BFq) = V^{(c)}(\BFq) + V^{(s)}(\BFq)$, we have
\begin{align*}
    \eta c\E[\Delta V(\BFq)] \leq 2\eta\left(\inner{\BFone_{n}}{\BFmu^*}+\inner{\BFone_m}{\BFlambda^*}\right)-2\sigma^\eta \sum_{j=1}^m \theta_jq_j^{(c)}\mathbbm{1}_{\left\{q_j^{(c)}>\tau_{\max}^\eta\right\}}-2\sigma^\eta \sum_{i=1}^n \phi_iq_i^{(s)}\mathbbm{1}_{\left\{q_i^{(s)}>\tau_{\max}^\eta\right\}}.
\end{align*}
Thus, the drift is negative outside the finite set $\mathcal{B}$ defined as
\begin{align*}
    \mathcal{B}\overset{\Delta}{=}\left\{\BFq: q_j^{(c)} \leq \max\left\{\tau_{\max}^\eta,\frac{2\eta}{\theta_j \sigma^\eta}\left(\inner{\BFone_{n}}{\BFmu^*}+\inner{\BFone_m}{\BFlambda^*}\right)\right\},q_i^{(s)} \leq \max\left\{\tau_{\max}^\eta,\frac{2\eta}{\phi_i \sigma^\eta}\left(\inner{\BFone_{n}}{\BFmu^*}+\inner{\BFone_m}{\BFlambda^*}\right)\right\}\right\}.
\end{align*}
By the Foster-Lyapunov theorem, the uniformized DTMC is positive recurrent. Moreover, using the moment bound theorem as a corollay of the Foster-Lyapunov theorem  \citep{HajekComm}, we get
\begin{align*}
  \E \left[\inner{\BFtheta}{\BFq^{(c)}}\right]+\E \left[\inner{\BFphi}{\BFq^{(s)}}\right]  \leq &
    \tau_{\max}^{\eta}\left(\sum_{j = 1}^{m} \theta_j  \Pr[q_j^{(c)}\leq \tau_{\max}^{\eta}]+ \sum_{i =1}^{n} \phi_i \Pr[q^{(s)}_i \leq \tau_{\max}^{\eta}]\right) \nonumber\\
    &+ \frac{\eta}{\sigma^{\eta}}\left(\inner{\BFone_n}{\omu} + \inner{\BFone_m}{\olam}\right). 
\end{align*}
Given the above bound on the expected queue length, we can bound the profit loss to be $O(\eta^{1/3})$ using the identical proof in Theorem \ref{theorem: twoprice}, so we omit the details. 
\Halmos
\endproof
\proof{Proof of Theorem \ref{theo: two_price_maxweight_SSC}.}
By Lemma \ref{lemma: profit-loss-two-price-max-weight} we have the required bound for profit loss under max-weight matching. By Lemma \ref{lemma: profit-loss-two-price-random}, we know that the profit loss is $O(\eta^{1/3})$ under randomized matching policy. 
To complete the proof, we will show that $\lim\inf_{\eta\to\infty} \lr / \eta^{1/3} = \Omega(n)$ for the two-price randomized matching algorithm.
We use the same graph and fluid solution as in the proof of Theorem \ref{theo: two_price_maxweight_SSC} (see Figure~\ref{fig: long_chain}) and show that the limiting scaled profit loss is of the order $n$. Note that as $\chi_{ij}=0$ for all $i \neq j$, the system behaves like $n$ independent two-sided queue under the randomized matching algorithm. For each single-link two-sided queue under any two-price pricing policy (which has loss $\lr/n$), 
by \eqref{eq: two-price-single-link-lower-bound},
the loss for
for any $\sigma^\eta = \gamma\eta^{2/3}$ is bounded by
\begin{align}
    \frac{\lr}{n} &\geq \left(\min_{i \in N, j\in M}\{A_j,B_i\}\frac{\gamma K_2}{4}+\min_{i \in N, j\in M}\{s_i^{(s)},s_j^{(c)}\}\max_{i \in N, j\in M} \left\{\frac{1}{2},\frac{\lambda_j^*}{2(\theta_j+\phi_i)\gamma}-0.5\right\}\right)\eta^{1/3} \\
    &\geq \frac{\min_{i \in N, j\in M}\{s_i^{(s)},s_j^{(c)}\} }{2} \eta^{1/3}.
\nonumber
\end{align}
(Recall that $A_j$ and $B_i$ represent the second derivative of $-F(\lambda_j)\lambda_j$ and $G(\mu_i)\mu_i$, respectively.)
Therefore, we have
\[
    \lim\inf_{\eta\to\infty} \frac{\lr}{\eta^{1/3}}
    \geq  \frac{\min_{i \in N, j\in M}\{s_i^{(s)},s_j^{(c)}\} }{2}n.  \Halmos 
\]
\endproof
\subsection{Theorem \ref{theo: lower_bound_n}}
\begin{lemma} \label{lemma: bias_complete_graph}
Let $\BFq_1, \BFq_2 \in \mathbb{Z}_+^{m+n}$ be such that $\inner{\BFone_m}{\BFq^{(c)}_1}-\inner{\BFone_n}{\BFq^{(s)}_1}=\inner{\BFone_m}{\BFq^{(c)}_2}-\inner{\BFone_n}{\BFq^{(s)}_2}$, then
\begin{align}
    h(\BFq_1) &\leq h(\BFq_2) \quad \textit{if } \inner{\BFone_m}{\BFq^{(c)}_1} \geq \inner{\BFone_m}{\BFq^{(c)}_2} \label{eq: bias_complete_graph_unequal}, \\
    h(\BFq_1) &= h(\BFq_2) \quad \textit{if } \inner{\BFone_m}{\BFq^{(c)}_1} = \inner{\BFone_m}{\BFq^{(c)}_2}, \label{eq: bias_complet_graph_equal}
\end{align}
where $h(\cdot)$ is the optimal bias function of the MDP defined in Eq~\eqref{MDP}.
\end{lemma}
\proof{Proof of Lemma~\ref{lemma: bias_complete_graph}.}
Similar to Lemma~\ref{lemma: monotonic}, As the relative value iteration may not converge in general for average cost MDP with countable state space, we work with the discounted MDP with the following Bellman equation with discount factor $\alpha \in (0, 1)$:
\begin{align*}
    V_{\alpha}(q) = \frac{1}{c}\max_{(\BFlambda,\BFmu,\BFx) \in Z(\BFq)}\left\{\sum_{j=1}^n F(\lambda_j)\lambda_j -\sum_{i=1}^n G(\mu_i)\mu_i -\inner{\BFone_{2n}}{\BFq}+\alpha\sum_{j=1}^n \lambda_j V_{\alpha}(\BFq+e_j^{(c)}-\BFx)\right. \\
 \left.+\alpha\sum_{i=1}^n \mu_i V_{\alpha}(\BFq+e_i^{(s)}-\BFx) +\alpha\left(c-\sum_{j=1}^n \lambda_j - \sum_{i=1}^n \mu_i\right)V_{\alpha}(\BFq) \right\}
\end{align*}
We consider the value iteration method to compute $V_\alpha(\BFq)$.
The value iteration starts with an arbitrary initial value $V_{\alpha, 0}(\BFq)$. In each iteration $k=1,2,\cdots$, for all $\BFq \in \mathbb{Z}_+^{m+n}$, the value iteration algorithm updates the value function by
\begin{align}
    V_{\alpha, k+1}(\BFq)=\frac{1}{c}\max_{(\BFlambda,\BFmu,\BFx) \in Z(\BFq)}\left\{\sum_{j=1}^n F(\lambda_j)\lambda_j -\sum_{i=1}^n G(\mu_i)\mu_i -\inner{\BFone_{2n}}{\BFq}+\alpha\sum_{j=1}^n \lambda_j V_{\alpha, k}(\BFq+e_j^{(c)}-\BFx)\right. \nonumber\\
 \left.+\alpha\sum_{i=1}^n \mu_i V_{\alpha, k}(\BFq+e_i^{(s)}-\BFx) +\alpha\left(c-\sum_{j=1}^n \lambda_j - \sum_{i=1}^n \mu_i\right)V_{\alpha, k}(\BFq) \right\} \label{eq:value_iteration_discounted_complete_graph}
\end{align}


Now we show using induction that \eqref{eq: bias_complete_graph_unequal} and \eqref{eq: bias_complet_graph_equal} holds for $V_{\alpha, k}$ for all $k \geq 0$.

\underline{Base Case:} 
Pick $V_{\alpha, 0}(\BFq)=0$ for all $\BFq \in \mathbb{Z}_+^{m+n}$. Thus, $V_{\alpha, 0}(\cdot)$ trivially satisfies Eqs \eqref{eq: bias_complete_graph_unequal}-\eqref{eq: bias_complet_graph_equal}.

\underline{Induction Step:} Suppose Eqs \eqref{eq: bias_complete_graph_unequal}-\eqref{eq: bias_complet_graph_equal} hold for $V_{\alpha, k}(\cdot)$, namely, after $k$ steps of value iteration.
Let $(\BFlambda^\star,\BFmu^\star,\BFx^\star_1) \in Z(\BFq_1)$ be the solution attaining the maximum on the right hand side of \eqref{eq:value_iteration_discounted_complete_graph}. Then, for $\BFq_1 \in \mathbb{Z}_+^{m+n}$, we have
\begin{align*}
   c V_{\alpha, k+1}(\BFq_1)=\sum_{j=1}^n F(\lambda_j^\star)\lambda_j^\star-\sum_{i=1}^n G(\mu_i^\star)\mu_i^\star-\inner{\BFone_{2n}}{\BFq_1}+\alpha\sum_{j=1}^n \lambda_j^\star V_{\alpha, k}(\BFq_1+e_j^{(c)}-\BFx^\star_1) \\
    +\alpha\sum_{i=1}^n \mu_i^\star V_{\alpha, k}(\BFq_1+e_i^{(s)}-\BFx^\star_1)+\alpha\left(c-\sum_{j=1}^n \lambda_j^\star - \sum_{i=1}^n \mu_i^\star\right)V_{\alpha, k}(\BFq_1)
\end{align*}
Consider another state $\BFq_2 \in \mathbb{Z}_+^{m+n}$ such that $\inner{\BFone_m}{\BFq^{(c)}_1}-\inner{\BFone_n}{\BFq^{(s)}_1}=\inner{\BFone_m}{\BFq^{(c)}_2}-\inner{\BFone_n}{\BFq^{(s)}_2}$ and $\inner{\BFone_m}{\BFq^{(c)}_1} \geq \inner{\BFone_m}{\BFq^{(c)}_2}$. As the compatibility graph forms a complete graph, we can construct $\BFx^\star_2 \in X(\BFq_2)$ such that $\inner{\BFone_m}{\BFq^{(c)}_1}-\inner{\BFone_n}{\BFx^{\star,(c)}_1} \geq \inner{\BFone_m}{\BFq^{(c)}_2}-\inner{\BFone_n}{\BFx^{\star,(c)}_2}$. Thus, by the induction hypothesis, we have
\begin{align*}
   c V_{\alpha, k+1}(\BFq_1)\leq {}&\sum_{j=1}^n F(\lambda_j^\star)\lambda_j^\star-\sum_{i=1}^n G(\mu_i^\star)\mu_i^\star-\inner{\BFone_{2n}}{\BFq_2}+\alpha\sum_{j=1}^n \lambda_j^\star V_{\alpha, k}(\BFq_2+e_j^{(c)}-\BFx^\star_2) \\
    &+\alpha\sum_{i=1}^n \mu_i^\star V_{\alpha, k}(\BFq_2+e_i^{(s)}-\BFx^\star_2)+\alpha\left(c-\sum_{j=1}^n \lambda_j^\star - \sum_{i=1}^n \mu_i^\star\right)V_{\alpha, k}(\BFq_2)  \\
    \leq{}& c V_{\alpha, k+1}(\BFq_2)
\end{align*}
This proves \eqref{eq: bias_complete_graph_unequal} for $V_{\alpha, k+1}$. Now, \eqref{eq: bias_complet_graph_equal} also follows by the same argument with $\BFq_1$ and $\BFq_2$ flipped. 

As $k\to\infty$, $\Delta V_{\alpha, k}(\BFq)$  converges to $\Delta V_{\alpha}^{*}(\BFq)$ (e.g., see Theorem 3.4 of \citet{feinberg2018convergence}).
Since the structural properties \eqref{eq: bias_complete_graph_unequal}-\eqref{eq: bias_complet_graph_equal} is preserved in the limit, $\Delta V_{\alpha}^{*}(\BFq)$ satisfies \eqref{eq: bias_complete_graph_unequal}-\eqref{eq: bias_complet_graph_equal}. Lastly, by \citet[Theorem 3.2]{cavazos1989weak}, there exists a subsequence of discount factors $\{\alpha_n\}$ and reference states $\{\BFq_n\}$ such that 
$$
V_{\alpha_n}^\star(\BFq) - V_{\alpha_n}^\star(\BFq_n) \to h^\star(\BFq) \quad \text{as} \quad \alpha_n \to 1,
$$
where, the convergence above is  pointwise. Additive normalization and pointwise convergence preserve the structural properties \eqref{eq: bias_complete_graph_unequal}-\eqref{eq: bias_complet_graph_equal}; thus, the resulting bias satisfies structural properties \eqref{eq: bias_complete_graph_unequal}-\eqref{eq: bias_complet_graph_equal}. This completes the proof.
$\Halmos$
\endproof
\proof{Proof of Proposition \ref{prop: optimal_policy_complete_graph}.}
Let $h(\cdot)$ be the bias function of the Bellman equation \eqref{MDP} and $(\BFlambda^\star,\BFmu^\star,\BFx^\star)$ be the optimal policy for the state $\BFq \in \mathbb{Z}_+^{2n}$. We have
\begin{align*}
& \max_{(\BFlambda,\BFmu,\BFx) \in Z(\BFq)}\left\{\sum_{j=1}^n \frac{\lambda_j}{c}\left(F(\lambda_j)+h(\BFq+\BFe_j^{(c)}-\BFx)-h(\BFq)\right)-\sum_{i=1}^n \frac{\mu_i}{c}\left(G(\mu_i)-h(\BFq+\BFe_i^{(s)}-\BFx)+h(\BFq)\right)\right\} \\
      ={}& \sum_{j=1}^n\frac{\lambda_j^\star}{c} \left(F(\lambda_j^\star)+h(\BFq+\BFe_j^{(c)}-\BFx^\star)-h(\BFq)\right)-\sum_{i=1}^n \frac{\mu_i^\star}{c}\left(G(\mu_i^\star)-h(\BFq+\BFe_i^{(s)}-\BFx^\star)+h(\BFq)\right) \\
    {\leq}{}& \sum_{j=1}^n \left(\frac{\inner{\BFone_n}{\BFlambda^\star}}{nc} F\left(\frac{\inner{\BFone_n}{\BFlambda^\star}}{n}\right)+\frac{\lambda_j^\star}{c}h(\BFq+\BFe_j^{(c)}-\BFx^\star)-\frac{\lambda_j^\star}{c} h(\BFq)\right)\\
    &-\sum_{i=1}^n \left(\frac{\inner{\BFone_n}{\BFmu^\star}}{nc} G\left(\frac{\inner{\BFone_n}{\BFmu^\star}}{n}\right)-\frac{\mu_i^\star}{c} h(\BFq+\BFe_i^{(s)}-\BFx^\star)+\frac{\mu_i^\star}{c} h(\BFq)\right) \\
    {=}{}& \sum_{j=1}^n \frac{\inner{\BFone_n}{\BFlambda^\star}}{nc} \left( F\left(\frac{\inner{\BFone_n}{\BFlambda^\star}}{n}\right)+h(\BFq+\BFe_j^{(c)}-\BFx^\star)- h(\BFq)\right)\\
    &-\sum_{i=1}^n \frac{\inner{\BFone_n}{\BFmu^\star}}{nc}\left( G\left(\frac{\inner{\BFone_n}{\BFmu^\star}}{n}\right)-h(\BFq+\BFe_i^{(s)}-\BFx^\star)+ h(\BFq)\right)\\
    \leq& \max_{(\BFlambda,\BFmu,\BFx) \in Z(\BFq)}\left\{\sum_{j=1}^n \frac{\lambda_j}{c}\left(F(\lambda_j)+h(\BFq+\BFe_j^{(c)}-\BFx)-h(\BFq)\right)-\sum_{i=1}^n \frac{\mu_i}{c}\left(G(\mu_i)-h(\BFq+\BFe_i^{(s)}-\BFx)+h(\BFq)\right)\right\},
\end{align*}
where the first inequality follows by Jensen's inequality (by Assumptions \ref{ass: monotonic} and \ref{ass: concavity}), the second equality holds because
$h(\BFq+\BFe_j^{(c)}-\BFx^\star)$ is equal for all $j\in M$ and $h(\BFq+\BFe_i^{(s)}-\BFx^\star)$ is equal for all $i\in N$ by
Lemma~\ref{lemma: bias_complete_graph},
and the second inequality holds by definition of the maximum operator.
Therefore, all the inequalities above must hold as equality.
Because the function $F(\lambda)\lambda$ is strictly concave and the function $G(\mu)\mu$ is strictly convex, we must have $\lambda^\star_j = \frac{\inner{\BFone_n}{\BFlambda^\star}}{n}$ for all $j\in M$ and $\mu^\star_i = \frac{\inner{\BFone_n}{\BFmu^\star}}{n}$ for all $i\in N$.
This completes the proof of property (a).

To prove property (b), suppose two states $\BFq_1$ and $\BFq_2$ satisfy
$\inner{\BFone_n}{\BFq_1^{(c)}}=\inner{\BFone_n}{\BFq_2^{(c)}}$ and $\inner{\BFone_n}{\BFq_1^{(s)}}=\inner{\BFone_n}{\BFq_2^{(s)}}$.
Suppose $(\BFlambda,\BFmu,\BFx_1) \in Z(\BFq_1)$. Because the instance assumes a complete graph, we can construct $\BFx_2 \in X(\BFq_2)$ such that $\inner{\BFone_n}{\BFx^{(c)}_1}=\inner{\BFone_n}{\BFx^{(c)}_2}$ and $\inner{\BFone_n}{\BFx^{(s)}_1}=\inner{\BFone_n}{\BFx^{(s)}_2}$. By Lemma~\ref{lemma: bias_complete_graph}, we have
\begin{align*}
     {}& \sum_{j=1}^n \frac{\lambda_j}{c}\left(F(\lambda_j)+h(\BFq_1+\BFe_j^{(c)}-\BFx_1)-h(\BFq_1)\right)-\sum_{i=1}^n \frac{\mu_i}{c}\left(G(\mu_i)-h(\BFq_1+\BFe_i^{(s)}-\BFx_1)+h(\BFq_1)\right) \\
     ={}& \sum_{j=1}^n \frac{\lambda_j}{c}\left(F(\lambda_j)+h(\BFq_2+\BFe_j^{(c)}-\BFx_2)-h(\BFq_2)\right)-\sum_{i=1}^n \frac{\mu_i}{c}\left(G(\mu_i)-h(\BFq_2+\BFe_i^{(s)}-\BFx_2)+h(\BFq_2)\right),
\end{align*}
which implies that
\begin{align*}
    {}& \max_{(\BFlambda,\BFmu,\BFx) \in Z(\BFq_1)}\left\{\sum_{j=1}^n \frac{\lambda_j}{c}\left(F(\lambda_j)+h(\BFq_1+\BFe_j^{(c)}-\BFx)-h(\BFq_1)\right)-\sum_{i=1}^n \frac{\mu_i}{c}\left(G(\mu_i)-h(\BFq_1+\BFe_i^{(s)}-\BFx)+h(\BFq_1)\right)\right\} \\
    ={}& \max_{(\BFlambda,\BFmu,\BFx) \in Z(\BFq_2)}\left\{\sum_{j=1}^n \frac{\lambda_j}{c}\left(F(\lambda_j)+h(\BFq_2+\BFe_j^{(c)}-\BFx)-h(\BFq_2)\right)-\sum_{i=1}^n \frac{\mu_i}{c}\left(G(\mu_i)-h(\BFq_2+\BFe_i^{(s)}-\BFx)+h(\BFq_2)\right)\right\}.
\end{align*}
Since $\inner{\BFone_{2n}}{\BFq_1} = \inner{\BFone_{2n}}{\BFq_2}$, 
if $(\BFlambda^\star,\BFmu^\star)$ is a maximizer for
\begin{align*}
    &\frac{\gamma}{c}+\frac{\inner{\BFone_{2n}}{\BFq_1}}{c}\\
    ={}& \max_{(\BFlambda,\BFmu,\BFx) \in Z(\BFq_1)}\left\{\sum_{j=1}^n \frac{\lambda_j}{c}\left(F(\lambda_j)+h(\BFq_1+\BFe_j^{(c)}-\BFx)-h(\BFq_1)\right)-\sum_{i=1}^n \frac{\mu_i}{c}\left(G(\mu_i)-h(\BFq_1+\BFe_i^{(s)}-\BFx)+h(\BFq_1)\right)\right\},
\end{align*}
then it is also a maximizer for
\begin{align*}
    &\frac{\gamma}{c}+\frac{\inner{\BFone_{2n}}{\BFq_2}}{c}\\
    ={}& \max_{(\BFlambda,\BFmu,\BFx) \in Z(\BFq_2)}\left\{\sum_{j=1}^n \frac{\lambda_j}{c}\left(F(\lambda_j)+h(\BFq_2+\BFe_j^{(c)}-\BFx)-h(\BFq_2)\right)-\sum_{i=1}^n \frac{\mu_i}{c}\left(G(\mu_i)-h(\BFq_2+\BFe_i^{(s)}-\BFx)+h(\BFq_2)\right)\right\}.
\end{align*}
This means that the optimal prices for $\BFq_1$ and $\BFq_2$ are equal, which completes the second part of the proposition. 

To prove property (c), assume that $(\BFlambda^\star,\BFmu^\star,\BFx^\star)$ is the optimal action for some state $\BFq$ and $\inner{\BFone_{2n}}{\BFx^\star} < 2 \min \{\inner{\BFone_n}{\BFq^{(c)}},\inner{\BFone_n}{\BFq^{(s)}}\}$. Let $\BFx^{\star \star} \in X(\BFq)$ be such that $\inner{\BFone_{2n}}{\BFx^{\star\star}} = 2 \min \{\inner{\BFone_n}{\BFq^{(c)}},\inner{\BFone_n}{\BFq^{(s)}}\}$.  
We have
\begin{align*}
    &\frac{\gamma}{c}+\frac{\inner{\BFone_{2n}}{\BFq}}{c} \\
    {=}{}& \max_{(\BFlambda,\BFmu,\BFx) \in Z(\BFq)}\left\{\sum_{j=1}^n \frac{\lambda_j}{c}\left(F(\lambda_j)+h(\BFq+\BFe_j^{(c)}-\BFx)-h(\BFq)\right)-\sum_{i=1}^n \frac{\mu_i}{c}\left(G(\mu_i)-h(\BFq+\BFe_i^{(s)}-\BFx)+h(\BFq)\right)\right\} \\
    ={}& \sum_{j=1}^n \frac{\lambda_j^\star}{c}\left(F(\lambda_j^\star)+h(\BFq+\BFe_j^{(c)}-\BFx^\star)-h(\BFq)\right)-\sum_{i=1}^n \frac{\mu_i^\star}{c}\left(G(\mu_i^\star)-h(\BFq+\BFe_i^{(s)}-\BFx^\star)+h(\BFq)\right) \\
    \leq{}& \sum_{j=1}^n \frac{\lambda_j^\star}{c}\left(F(\lambda_j^\star)+h(\BFq+\BFe_j^{(c)}-\BFx^{\star\star})-h(\BFq)\right)-\sum_{i=1}^n \frac{\mu_i^\star}{c}\left(G(\mu_i^\star)-h(\BFq+\BFe_i^{(s)}-\BFx^{\star\star})+h(\BFq)\right), 
\end{align*}
where the last inequality follows by Eq \eqref{eq: bias_complete_graph_unequal} in Lemma \ref{lemma: bias_complete_graph}. Since $\BFx^{\star\star} \in X(\BFq)$, this proves that
$(\BFlambda^\star,\BFmu^\star,\BFx^{\star\star})$ must also be an optimal action.
$\Halmos$
\endproof
\begin{lemma} \label{lemma: n_lower_bound}
If Condition \ref{cond: n_dep} is satisfied, there exists a constant $\mathcal{M}>0$ that does not depend on $\eta$ such that
\begin{align*}
    \mathbb{E}\left[f^2\left(\frac{z}{\eta^\alpha}\right)+g^2\left(\frac{z}{\eta^\alpha}\right)\right] \geq \frac{\delta^2 \Gamma^2}{2} \quad \forall \eta>\mathcal{M},
\end{align*}
where $z$ is defined in \eqref{eq: imbalance}.
\end{lemma}
\proof{Proof of Lemma \ref{lemma: n_lower_bound}.}
By Condition \ref{cond: n_dep}\ref{cond: stability_n}, we can lower bound the required expectation by
\begin{align*}
    \mathbb{E}\left[f^2\left(\frac{z}{\eta^\alpha}\right)+g^2\left(\frac{z}{\eta^\alpha}\right)\right] \geq \delta^2 \Gamma^2 \mathbb{P}\left[|z(\infty)|>\kappa\eta^\alpha\right].
\end{align*}
Considering the coupling defined in Lemma \ref{lemma: boundedf2g2} and using the lower bound given by \eqref{eq: lb_couping_imbalance}, we have
\begin{align*}
    \Pr[|z(\infty)|>\kappa\eta^\alpha] \geq \left(\frac{\lambda^\star-2\Gamma\eta^{\beta-1}}{\lambda^\star+2\Gamma\eta^{\beta-1}}\right)^{\kappa\eta^\alpha+1}.
\end{align*}
Taking the limit as $\eta \rightarrow \infty$ and using Condition \ref{cond: n_dep}\ref{cond: scaling_n}, the right-hand side converges to
\begin{align*}
    \lim_{\eta \rightarrow \infty} \left(\frac{\lambda^\star-2\Gamma\eta^{\beta-1}}{\lambda^\star+2\Gamma\eta^{\beta-1}}\right)^{\kappa\eta^\alpha+1} = 1.
\end{align*}
Thus, there exists $\mathcal{M}>0$ (which does not depend on $\eta$) such that 
\begin{align*}
    \Pr[|z(\infty)|>\kappa\eta^\alpha] \geq \frac{1}{2} \quad \forall \eta>\mathcal{M}.
\end{align*}
This completes the proof.
$\Halmos$
\endproof
\proof{Proof of Theorem \ref{theo: lower_bound_n}.}
Consider an auxiliary system where we match all compatible customer-server pairs. In this system, either only customers or servers are waiting. The pricing decisions in this auxiliary system is exactly the same (as it only depends on the imbalance due to Condition~\ref{cond: n_dep}), so, the revenue and cost is the same, while the queue lengths are always lower. So, the profit-loss in the auxiliary system lower bounds the profit-loss for the original system. Noting that Condition~\ref{cond: n_dep} for the auxiliary system implies Condition~\ref{ass: generalpricinggraph}, by Theorem \ref{theorem: lowerboundmultiplelink}, we have
\begin{align*}
    L^\eta \geq K \eta^{1/3},
\end{align*}
where $K$ is given by Eq \eqref{eq: constant_cube_root}. In particular, by substituting ${\Gamma}_j={\Psi}_i=\Gamma$ for all $i \in N, j \in M$, $\epsilon=n\delta^2 \Gamma^2/2$, $\lambda_j^\star=\lambda^\star$ for all $j \in M$, $\BFs=\BFone$, and $A_j=A$ for all $j \in M$, $B_i=B$ for all $i \in N$, we get
\begin{align*}
    K = n^{1/3}\frac{3}{2^{2/3}}\left(\frac{\min\left\{A,B
      \right\} \delta^2}{8}\right)^{1/3}\left(\frac{\lambda^*}{4}\right)^{2/3}
\end{align*}
Thus, we have $K\geq\Omega(n^{1/3})$, which concludes the proof.
$\Halmos$
\endproof
\subsection{Relaxing Condition \ref{condition: perfect_matching}}
\label{subsec: relax-perfect-matching}
In this section, we show that Condition \ref{condition: perfect_matching} can be relaxed without loss of generality. Suppose we are given a graph $G(N_s \cup N_c,E)$ and the fluid solution $(\BFlambda^*,\BFmu^*)$ such that the CRP condition is satisfied but Condition \ref{condition: perfect_matching} is violated. We will construct another graph $G'(N_s' \cup N_c',E')$ such that Condition \ref{condition: perfect_matching} is satisfied and use max-weight matching and fluid price policy/two-price policy for this new graph. Then, we will conduct a similar analysis as in Theorem \ref{theo: fluid_max_weight_SSC} and Theorem \ref{theo: two_price_maxweight_SSC} to get the limiting bound on the scaled loss in profit.
\subsubsection{Fluid price policy and max-weight matching policy.}
The idea of the graph construction is to split each customer/server type into multiple types so that we can eventually achieve the same number of customer and server types.
We assume that there exists a $\delta>0$ and such that the fluid solution $(\BFlambda^*,\BFmu^*)$ are integral multiple of $\delta$. (In practice, if the model parameters are rational numbers, this can always be achieved.) In addition, let $\lambda_{\max}$ be an upper bound of the entries in $(\BFlambda^*,\BFmu^*)$. 
Specifically, for each customer type $j \in M$, we split arrivals into ${\lambda_j^*}/{\delta}$ separate queues.
Let $q_{jl}^{(c)}$ for all $l \in [{\lambda_j^*}/{\delta}]$ denote the queue length of the $l^{th}$ replication. Suppose the arrival rate for (original) type $j$ customer is $\lambda_j^\eta(\BFq)=\eta\lambda_j^*-\delta\sum_{i=1}^{\lambda_j^*/\delta}\mathbbm{1}_{\{q_{jl}^{(c)}=q_{\max}^\eta\}}$. Then, the arrival rate of each split queue is $\lambda_{jl}^\eta(\BFq)=\eta\delta\mathbbm{1}_{\{q_{jl}^{(c)}<q_{\max}^\eta\}}$. 
We split the server types similarly. 
This leads to a modified graph $G'(N_s' \cup N_c',E')$ with $N_c=\{jl: l \in [\lambda_j/\delta],j \in {M}\}$, $|N_c'|=\inner{\BFone_m}{\BFlambda^*}/\delta$ and $N_s=\{im: m \in [\mu_i/\delta],i \in {N}\}$, $|N_s'|=\inner{\BFone_n}{\BFmu^*}/\delta$. The edges $E'$ of the modified graph satisfy the following:  $(i,j) \in E$ if and only for all $l \in [\lambda_{j}/\delta]$ and for all $m \in [\mu_{i}/\delta]$, we have $(im,jl) \in E'$.
We continue to use the max-weight matching algorithm under the modified graph $G'(N_s' \cup N_c',E')$.
Next, we show that Condition~\ref{condition: perfect_matching} holds for the modified graph. Since $\inner{\BFone_m}{\BFlambda^*}=\inner{\BFone_n}{\BFmu^*}$, we have $|N_c'|=|N_s'|=\inner{\BFone_m}{\BFlambda^*}/\delta \leq n\lambda_{\max}/\delta$. In addition, for any $J \subsetneq N_c'$, we have
\begin{align}
    \delta|J|&=\sum_{jl \in J} \delta \leq \sum_{j\in {M}: \exists l, jl \in J} \lambda_j^* {<} \sum_{i: \exists j, \exists l,(i,j) \in E, jl \in J} \mu_i^*=\sum_{i \in {N}: \exists jl\in J,(i,j) \in E} \sum_{m \in [\mu_i^*/\delta]}\delta \nonumber \\
    &{=}  \sum_{im \in N_c': \exists jl \in J, (im,jl) \in E' }\delta=\delta|N(J)|, \label{eq: hallscondition}
\end{align}
where the second inequality follows from the CRP condition and the third equality follows as 
\begin{align*}
    \left\{im: i \in {M}, m \in [\mu_i^*/\delta], \exists jl \in J,(i,j) \in E\right\} \ \iff\ \{im \in N_c': \exists jl \in J, (im,jl) \in E'\}.
\end{align*}
Eq~\eqref{eq: hallscondition} is Hall's condition for all the customer types. We can similarly verify Hall's condition for any subset of servers. This implies that there exists a perfect matching in $G'(N_s' \cup N_c',E')$. Thus, Condition \ref{condition: perfect_matching} is satisfied for this modified graph. 
Applying Proposition \ref{prop: max-weight-optimal-fluid-SSC} to the modified graph, we have
\begin{align*}
    \lim_{\eta\to\infty} q_{\max}^\eta\left( \sum_{j=1}^m\sum_{l=1}^{\lambda_j^*/\delta} \lambda_j^* \Pr\left[q_{jl}^{(c)}=q_{\max}^\eta\right]+\sum_{i=1}^n \sum_{m=1}^{\mu_i^*/\delta}\mu_i^*\Pr\left[q_{im}^{(s)}=q_{\max}^\eta\right]\right) =
    \frac{\inner{\BFone_n}{\BFlambda^*}+\inner{\BFone_n}{\BFmu^*}}{2n}.
\end{align*}
In addition, it is trivially true that
\begin{align*}
    \E[\inner{\BFone_{n+m}}{\BFq}] \leq \frac{\inner{\BFone_m}{\BFmu^*}}{\delta}q_{\max}^\eta \leq n \frac{\lambda_{\max}}{\delta} q_{\max}^\eta.
\end{align*}
The above two results together with $q_{\max}^\eta=n^{-1/2}\eta^{1/2}$ will produce $\limsup_{\eta \rightarrow \infty} \lr/\sqrt{\eta}=O(\sqrt{n})$.
\subsubsection{Two-price policy and max-weight matching policy}
As in previous subsection, 
we split customer and server types to satisfy Condition~\ref{condition: perfect_matching}.
For customer type $j$, let $q_{jl}^{(c)}$ for all $l \in \left[{\lambda_{j}^*}/{\delta}\right]$ denote the queue length for the $l^{th}$ replication. The arrival rate for the $l^{th}$ replication is $\lambda_{jl}^\eta(\BFq)=\eta\delta-\sigma^\eta\mathbbm{1}_{\left\{q_{jl}^{(c)}>\tau_{\max}^\eta\right\}}$. Similarly, we replicate the servers as well. 
We apply the two-price max-weight matching policy for the modified graph $G'(N_s' \cup N_c',E')$. By Proposition \ref{prop: max-weight-optimal-two-price-SSC},  we have
\begin{align*}
   \lim_{\eta \to \infty} \frac{\sigma^\eta}{\eta}\E[\inner{\BFone_{2n}}{\BFq}]= \delta\frac{\inner{\BFone_m}{\BFlambda^*}+\inner{\BFone_n}{\BFmu^*}}{2\inner{\BFone_m}{\BFlambda^*}}.
\end{align*}
Next, we analyze the profit loss similarly as in the proof of Theorem \ref{theo: two_price_maxweight_SSC}. By using the definition of the pricing policy and Taylor's theorem, we get
\begin{align*}
    \lr={}&\eta\left(\sum_{j=1}^m F_j(\lambda_j^*)\lambda_j^*-\E\left[\sum_{j=1}^m \left(\lambda_j^*-\frac{\sigma^\eta}{\eta}\sum_{l=1}^{\lambda_j^*/\delta}\mathbbm{1}_{\{q_{jl}>\tau_{\max}^\eta\}}\right)F_j\left(\lambda_j^*-\frac{\sigma^\eta}{\eta}\sum_{l=1}^{\lambda_j^*/\delta}\mathbbm{1}_{\{q_{jl}>\tau_{\max}^\eta\}}\right)\right]\right) \nonumber \\
    &-\eta\left(\sum_{i=1}^n G_i(\mu_i^*)\mu_i^*-\E\left[\sum_{i=1}^n \left(\mu_i^*-\frac{\sigma^\eta}{\eta}\sum_{m=1}^{\mu_i^*/\delta}\mathbbm{1}_{\{q_{im}>\tau_{\max}^\eta\}}\right)G_i\left(\mu_i^*-\frac{\sigma^\eta}{\eta}\sum_{m=1}^{\mu_i^*/\delta}\mathbbm{1}_{\{q_{im}>\tau_{\max}^\eta\}}\right)\right]\right) \\
    ={}& \sigma^\eta \left(\sum_{j=1}^m \left(\lambda_j^*F_j'(\lambda_j^*)+F_j(\lambda_j^*)\right)\E\left[\sum_{l=1}^{\lambda_j^*/\delta}\mathbbm{1}_{\{q_{jl}>\tau_{\max}^\eta\}}\right]-\sum_{i=1}^m \left(\mu_i^*G_i'(\mu_i^*)+G_i(\mu_i^*)\right)\E\left[\sum_{m=1}^{\mu_i^*/\delta}\mathbbm{1}_{\{q_{im}>\tau_{\max}^\eta\}}\right]\right) \nonumber \\
    &-\frac{(\sigma^\eta)^2}{\eta}\sum_{j=1}^m\left( \lambda_j^*F''_j(\lambda_j^*)\E\left[\sum_{l=1}^{\lambda_j^*/\delta}\mathbbm{1}_{\{q_{jl}>\tau_{\max}^\eta\}}\right]+F_j'(\lambda_j^*)\E\left[\left(\sum_{l=1}^{\lambda_j^*/\delta}\mathbbm{1}_{\{q_{jl}>\tau_{\max}^\eta\}}\right)^2\right]\right) \\
    &+\frac{(\sigma^\eta)^2}{\eta}\sum_{i=1}^n\left( \mu_i^*G''_i(\mu_i^*)\E\left[\sum_{m=1}^{\mu_i^*/\delta}\mathbbm{1}_{\{q_{im}>\tau_{\max}^\eta\}}\right]+G_i'(\mu_i^*)\E\left[\left(\sum_{m=1}^{\mu_i^*/\delta}\mathbbm{1}_{\{q_{im}>\tau_{\max}^\eta\}}\right)^2\right]\right)+o(\eta^{1/3}) \\
    {=}{}&-\frac{(\sigma^\eta)^2}{\eta}\sum_{j=1}^m\left( \lambda_j^*F''_j(\lambda_j^*)\E\left[\sum_{l=1}^{\lambda_j^*/\delta}\mathbbm{1}_{\{q_{jl}>\tau_{\max}^\eta\}}\right]+F_j'(\lambda_j^*)\E\left[\left(\sum_{l=1}^{\lambda_j^*/\delta}\mathbbm{1}_{\{q_{jl}>\tau_{\max}^\eta\}}\right)^2\right]\right) \\
    &+\frac{(\sigma^\eta)^2}{\eta}\sum_{i=1}^n\left( \mu_i^*G''_i(\mu_i^*)\E\left[\sum_{m=1}^{\mu_i^*/\delta}\mathbbm{1}_{\{q_{im}>\tau_{\max}^\eta\}}\right]+G_i'(\mu_i^*)\E\left[\left(\sum_{m=1}^{\mu_i^*/\delta}\mathbbm{1}_{\{q_{im}>\tau_{\max}^\eta\}}\right)^2\right]\right)+o(\eta^{1/3}) \\
    \leq{}& \frac{(\sigma^\eta)^2\lambda_{\max}}{\eta\delta}\left(\sum_{i=1}^n \mu_i^*G''_i(\mu_i^*)+G_i'(\mu_i^*)-\sum_{j=1}^m \lambda_j^*F''_j(\lambda_j^*)+F_j'(\lambda_j^*)\right)+o(\eta^{1/3}),
\end{align*}
where the third equality holds by claim \ref{claim}:
as the CRP condition is satisfied, there exists a fluid solution such that $\chi_{ij}^*>0$ for all $(i,j) \in E$. Setting $\sigma^\eta=\eta^{2/3}n^{-1/3}$, we get 
\begin{align*}
    \limsup_{\eta \rightarrow \infty} \frac{\lr}{\eta^{1/3}}\leq n^{1/3}\left(\lambda_{\max}\frac{\sum_{i=1}^n \mu_i^*G''_i(\mu_i^*)+G_i'(\mu_i^*)-\sum_{j=1}^m \lambda_j^*F''_j(\lambda_j^*)+F_j'(\lambda_j^*)}{n\delta}+\delta\max_{i \in N, j \in M}\{s_i^{(s)},s_j^{(c)}\}\right). 
\end{align*}
\end{APPENDICES}
\end{document}